# $L^2 - L^\infty$ DECAY ESTIMATES AND INVISCID LIMITS FOR GLOBAL SMOOTH SOLUTIONS TO THE COMPRESSIBLE NAVIER-STOKES-RIESZ SYSTEM

AGNIESZKA ŚWIERCZEWSKA-GWIAZDA, YUAN XU, AND JUNHAO ZHANG

Abstract. We study the Cauchy problem in $\mathbb{R}^3$ for the repulsive compressible Navier-Stokes-Riesz system with Riesz exponent $0 < s < 1$ and viscosity $0 < \varepsilon \leq 1$, where the Riesz interaction $\nabla(-\Delta)^{-s}(\rho - \bar{\rho})$ is a generalization of the Coulomb interaction for electrons. For small perturbations of a constant equilibrium, with the solenoidal component of the initial velocity of order $\mathcal{O}(\varepsilon)$, we prove the global existence and uniqueness of smooth solutions. We derive time-decay estimates in $L^2$ norms and $L^\infty$ norms that capture both uniform-in-$\varepsilon$ dispersive behavior and viscosity-dependent dissipation. We further establish a global-in-time inviscid limit to the irrotational global solution of the compressible Euler-Riesz system whose initial data consist of the same density and the curl-free component of the velocity, with an explicit convergence rate in $W^{k,p}$ norms. The proof combines viscosity-adapted dispersive estimates, normal-form analysis and nonlinear energy estimates with control of both negative and positive Sobolev norms.

## 1. Introduction and Main Results

1.1. **Introduction.** The motion of a charged compressible flow with its self-consistent force can be described by the following compressible Navier-Stokes-Riesz system:

$$\partial_t \rho^\varepsilon + \mathrm{div}\,(\rho^\varepsilon u^\varepsilon) = 0, \tag{1.1}$$

$$\partial_t(\rho^\varepsilon u^\varepsilon) + \mathrm{div}\,(\rho^\varepsilon u^\varepsilon \otimes u^\varepsilon) + \nabla P(\rho^\varepsilon) + \rho^\varepsilon \nabla \phi^\varepsilon = \varepsilon \mathcal{L} u^\varepsilon, \tag{1.2}$$

$$(-\Delta)^s \phi^\varepsilon = \rho^\varepsilon - \bar{\rho}, \tag{1.3}$$

for $(t, x) \in (0, \infty) \times \mathbb{R}^3$. Here $\rho^\varepsilon \geq 0$, $u^\varepsilon \in \mathbb{R}^3$, and $\phi^\varepsilon \in \mathbb{R}$ denote the density, the velocity, and the self-consistent potential, respectively. The pressure $P$ is given by

$$P(\rho^\varepsilon) = K(\rho^\varepsilon)^\gamma,$$

where $K > 0$ is the specific gas constant and $\gamma > 1$ is the adiabatic exponent. The viscous term $\mathcal{L}u^\varepsilon$ is represented as

$$\mathcal{L}u^\varepsilon = \mu \Delta u^\varepsilon + (\mu + \nu)\nabla \mathrm{div}\, u^\varepsilon,$$

where $\mu > 0$ and $2\mu + 3\nu \geq 0$. The potential $\phi^\varepsilon$, governed by a nonlocal fractional diffusion process (1.3), is often called the Riesz potential. Here we consider $s \in (0, 1)$, then the fractional Laplacian $(-\Delta)^s$ is defined by (1.12) with $m(\xi) := |\xi|^{2s}$, and $\bar{\rho} > 0$ is a given constant background state. The Riesz potential has been widely applied to plasma physics, star or galaxy dynamics, and swarming models in mathematical biology; see, for instance, [4, 5, 6, 27].

In this paper, we consider the Cauchy problem for (1.1)-(1.3) with the initial condition and the far field condition,

$$(\rho^\varepsilon, u^\varepsilon)(0, x) = (\rho_0^\varepsilon, u_0^\varepsilon)(x), \quad \lim_{|x| \to \infty} (\rho^\varepsilon, u^\varepsilon)(t, x) = (\bar{\rho}, 0), \tag{1.4}$$

where $(\rho_0^\varepsilon, u_0^\varepsilon)$ is the given initial data. We also denote $|\nabla|^s\phi_0^\varepsilon := |\nabla|^{-s}(\rho_0^\varepsilon - \bar\rho)$ with $|\nabla|^s := (-\Delta)^{\frac{s}{2}}$. Our first goal is to investigate the global well-posedness of the problem. In particular, for $s = 1$, which reduces the system (1.1)-(1.3) to the electron compressible Navier-Stokes-Poisson system, there has been an extensive literature concerning the problem, and the interaction potential is also called the Coulomb potential. Therefore, we recall the related results in this aspect. For $\bar\rho \in L^1 \cap L^\infty$, Ducomet, Feireisl, Petzeltová, and Straškraba in [13] showed the global existence of weak solutions with the adiabatic index $\gamma > 3/2$. Although they only studied the flow with self-gravitation (the attractive case), their proof, as observed in [2], works also in the case with the Coulomb force (the repulsive case). Later, Bella [2] further studied the long-time behavior of the weak solutions and proved that, in the attractive case, any weak solution converges to a unique weak solution to the corresponding stationary problem. For $\bar\rho$ that does not vanish at infinity, much work has been done around the steady equilibrium when $\bar\rho \equiv$ const. Li, Matsumura, and Zhang [24] obtained the global existence of the strong solution provided $(\rho_0 - \bar\rho, u_0) \in L^1 \cap H^N$ with $N \geq 4$. Moreover, by carrying out a spectral analysis of the semigroup associated with the linearized system, they proved the following optimal time-decay rates,

$$\begin{aligned} &\|(\rho - \bar\rho)(t)\|_{L^2} \lesssim (1+t)^{-\frac{3}{4}}, \quad \|(\rho - \bar\rho)(t)\|_{L^\infty} \lesssim (1+t)^{-\frac{3}{2}}, \\ &\|u(t)\|_{L^2} \lesssim (1+t)^{-\frac{1}{4}}, \quad \|u(t)\|_{L^\infty} \lesssim (1+t)^{-1}. \end{aligned} \tag{1.5}$$

Intuitively, the electric potential may enhance the dissipation of density such that the density exhibits a faster time-decay rate than the velocity. The natural energy of the system contains an electric-field contribution, so finite-energy perturbations require $\nabla\phi_0 \in L^2$, which further implies $\Lambda^{-1}(\rho_0 - 1) \in L^2$ by the Poisson equation. Motivated by this energy structure, in [35], Wang used negative Sobolev spaces. He established the time-decay rates via a refined pure energy method in $H^N \cap \dot H^{-s}$ with $N \geq 3$ and $s \in [0, 3/2)$. Moreover, by the Sobolev embedding $L^q \subset \dot H^{-s}$ with $q \in (1, 2]$,

$$\begin{aligned} &\|(\rho - \bar\rho)(t)\|_{L^2} \lesssim (1+t)^{-\frac{3}{2}(\frac{1}{q}-\frac{1}{2})-\frac{1}{2}}, \quad \|(\rho - \bar\rho)(t)\|_{L^\infty} \lesssim (1+t)^{-\frac{3}{2q}-\frac{1}{2}}, \\ &\|u(t)\|_{L^2} \lesssim (1+t)^{-\frac{3}{2}(\frac{1}{q}-\frac{1}{2})}, \quad \|u(t)\|_{L^\infty} \lesssim (1+t)^{-\frac{3}{2q}}, \end{aligned}$$

provided $(\rho_0, \nabla\phi_0, u_0) \in L^q$. For more related results, we refer the readers to [34, 36] and the references therein. There are very few results concerning the compressible Navier-Stokes-Riesz system. Recently, Wang and Zhang [37] proved the global existence of strong solutions to the generalized compressible Navier-Stokes-Poisson system with a Lévy diffusion process and derived optimal time-decay rates. Their result can be applied to the compressible Navier-Stokes-Riesz system for $s \in (1/2, 1)$, with which they obtained the global strong solution to (1.1)-(1.3) satisfying

$$\|(\rho - \bar\rho)(t)\|_{L^2} \lesssim (1+t)^{-\frac{3}{4}}, \quad \|u(t)\|_{L^2} \lesssim (1+t)^{-\frac{3-2s}{4}}. \tag{1.6}$$

However, it is still unknown whether the global strong solution to the compressible Navier-Stokes-Riesz system exists for $s \in (0, 1/2)$. In this paper, we will prove that the global existence result holds for the full range of $s \in (0, 1)$ and we will also recover the $L^\infty$ time-decay estimate for the solution.

In (1.1)-(1.3), if we take $\varepsilon = 0$, then the system is reduced to the compressible Euler-Riesz system. The literature for the compressible Euler-Poisson system ($s = 1$) for charged particles is also rather complete. Let us first review the related results. The existence of global large solutions with spherical symmetry was studied in several situations. In [33], the weak entropy solution was constructed in the BV space, while in [8, 9], the finite-energy weak solution was constructed. In contrast, it is generally believed that the global smooth solution does not exist. Indeed, similar to the result [31] for the compressible Euler system, Perthame [28] established a blow-up result for spherically symmetric initial data with large initial energy. However, for irrotational flows, by exploiting the decay property of the "Klein-Gordon effect", Guo [15] proved the global existence of smooth solutions with small initial data. Later, Guo's result was extended by Germain, Masmoudi, and Pausader in [21] by removing the neutrality assumption $\int_{\mathbb{R}^3} \rho_0 - \bar{\rho}\,\mathrm{d}x = 0$. The study was also extended to $\mathbb{R}^2$ by Ionescu, Pausader, Li, and Wu in [22, 26]. Besides, Blanc, Danchin, Ducomet, and Necasova [3] proved the existence of global solutions with Sobolev regularity for small initial density and the initial velocity that is close to a reference velocity $\tilde{u}$ such that the spectrum of $D\tilde{u}$ is bounded away from zero, which is called the "dispersive" spectral condition. Recently, for the compressible Euler-Riesz system, Carrillo, Charles, Chen, and Yuan [7] constructed the finite-energy weak solution with spherical symmetry for large initial data. In addition, Choi and Jeong [10] established the local-in-time existence and uniqueness of classical solutions and presented sufficient conditions on the finite-time blow-up. Moreover, by using the idea of "dispersive" spectral condition on the initial velocity, the unique global solution was constructed in [14] by Danchin and Ducomet, as well as in [11] by Choi, Jung, and Lee. Very recently, Choi, Jung, and Lee [12] proved the global existence of smooth solutions for small irrotational initial perturbations. Therefore, it is natural to ask whether solutions of the compressible Navier-Stokes-Riesz system converge to solutions of the compressible Euler-Riesz system in the inviscid limit. Our second goal in this paper is to give this question an affirmative answer.

We introduce the Leray projector $\mathcal{P}$ onto divergence-free vector fields and $\mathcal{P}^{\perp}$ onto curl-free vector fields

$$\mathcal{P} := \mathrm{Id} - \mathcal{P}^{\perp} \quad \text{and} \quad \mathcal{P}^{\perp} := \nabla\Delta^{-1}\mathrm{div}\,.$$

Our main result can be summarized as:

**Theorem 1.1.** *If the initial density perturbation $\rho_0^{\varepsilon} - \bar{\rho}$ and the potential part of the initial velocity $\mathcal{P}^{\perp}u_0^{\varepsilon}$ is sufficiently small, and if the incompressible part of the initial velocity $\mathcal{P}u_0^{\varepsilon}$ is small relative to $\varepsilon$, then there exists a unique global smooth solution to the system* (1.1)-(1.3) *with* (1.4). *Additionally, if the initial perturbation in density and the potential part of the initial velocity are bounded, and the incompressible part is $\mathcal{O}(\varepsilon)$ in $\dot{H}^{-\vartheta}$ with $0 < \vartheta \le 1/2$, then the $L^{\infty}$ norm of the solution decays at the rate $\min\left\{\varepsilon, (1+t)^{-\frac{\vartheta}{2}}\right\} + (1+t)^{-\beta_s}$, while for the $L^2$ norm of the solution, the density perturbation decays at the rate $\min\left\{\varepsilon, (1+t)^{-\frac{\vartheta}{2}}\right\} + (1+\varepsilon^{\Upsilon}t)^{-\frac{s}{2}-\frac{\vartheta}{2}}$ and the velocity perturbation decays at the rate $(1+\varepsilon^{\Upsilon}t)^{-\frac{\vartheta}{2}}$. Here, $\beta_s$ is defined by*

$$\beta_s := \begin{cases} \dfrac{3}{2}\left(1 - \dfrac{2}{p}\right) & \text{when } 0 < s \le 1/2 \text{ and } 6 < p < 12, \\[2mm] \dfrac{4}{3}\left(1 - \dfrac{2}{p}\right) & \text{when } 1/2 < s < 1 \text{ and } 8 < p < 16, \end{cases}$$

*and $\Upsilon$ is a constant depending on $\vartheta$. Moreover, in the inviscid-limit sense, the solution converges to the global smooth solution of the compressible Euler-Riesz system with initial data $(\widetilde{\rho}_0, \widetilde{u}_0) := (\rho_0^{\varepsilon}, \mathcal{P}^{\perp}u_0^{\varepsilon})$ at the $L^p$-rate $\mathcal{O}(\varepsilon^{\lambda_s})$ globally in time where $\lambda_s$ is defined by*

$$\lambda_s := \begin{cases} \dfrac{3(p-2)}{5p-6} & \text{when } 0 < s \le 1/2 \text{ and } 6 < p < 12, \\[2mm] \dfrac{4(p-2)}{7p-8} & \text{when } 1/2 < s < 1 \text{ and } 8 < p < 16. \end{cases}$$

**Remark 1.2.** *The precise statement of the result can be found in Theorem 1.9 and Theorem 1.11. It is worth noting that our analysis yields an additional $L^2$ decay rate. This reveals the dissipation of energy caused by the viscosity, which makes a distinct difference from the compressible Euler-Riesz (Poisson) system (*[15, 20, 22, 26, 29, 12]*). Our result also extends the result in* [37] *for the Navier-Stokes-Riesz system* (1.1)–(1.3) *with full range $s \in (0,1)$.*

In order to study the time-asymptotic stability and vanishing viscosity limit in a unified framework, we are inspired by Rousset and Sun's result in [29]. To be more specific, we first construct an approximate solution $(\rho, u)$ by the following system with irrotational initial data,

$$\begin{cases} \partial_t \rho + \operatorname{div}(\rho u) = 0, \\ \partial_t u + u \cdot \nabla u + \dfrac{1}{\gamma - 1} \nabla(\rho^{\gamma-1}) + \nabla \phi = \varepsilon \mathcal{L} u, \\ (-\Delta)^s \phi = \rho - \bar{\rho}, \\ (\rho, u)(0, x) = (\rho_0^\varepsilon, \mathcal{P}^\perp u_0^\varepsilon)(x). \end{cases} \tag{1.7}$$

Clearly, taking $\varepsilon = 0$ in (1.7), we obtain the compressible Euler-Riesz system. Therefore (1.7) is a viscous approximation of the Euler-Riesz system. In addition, this system also preserves the irrotationality of the flow if it does initially. Indeed, taking the curl in the second equation of system (1.7) gives the equation for $\omega \overset{\text{def}}{=} \operatorname{curl} u$

$$\begin{cases} \partial_t \omega + \omega \operatorname{div} u + (u \cdot \nabla)\omega - (\omega \cdot \nabla) u - \varepsilon \mathcal{L} \omega = 0, \\ \omega(0, x) = \omega_0(x) = 0. \end{cases}$$

By the standard energy estimate and Grönwall's inequality, we obtain

$$\|\omega(t)\|_{L^2}^2 \leq e^{C \int_0^t \|\nabla u(s)\|_{L^\infty} \, \mathrm{d}s} \|\omega_0\|_{L^2}^2 = 0,$$

which implies that

$$\nabla \times u = 0.$$

Moreover, since the artificial viscous term $\varepsilon \mathcal{L} u$ is added, the system is still dissipative. For (1.7), we obtain the following result:

**Theorem 1.3.** *For $0 < s < 1$, let $p$ be chosen according to*

$$\begin{cases} 6 < p < 12 & \text{when } 0 < s \leq 1/2, \\ 8 < p < 16 & \text{when } 1/2 < s < 1. \end{cases} \tag{1.8}$$

*Then there exists $\delta_0 > 0$ such that for any initial data satisfying*

$$\sup_{\varepsilon \in (0,1]} \Bigg( \Big\| \Big( \rho_0^\varepsilon - \bar{\rho}, |\nabla|^s \phi_0^\varepsilon, \mathcal{P}^\perp u_0^\varepsilon \Big) \Big\|_{\dot{H}^{-(1-s)} \cap H^N} + \Big\| \Big( \rho_0^\varepsilon - \bar{\rho}, |\nabla|^s \phi_0^\varepsilon, \mathcal{P}^\perp u_0^\varepsilon \Big) \Big\|_{L^1} + \Big\| \Big( \rho_0^\varepsilon - \bar{\rho}, |\nabla|^s \phi_0^\varepsilon, \mathcal{P}^\perp u_0^\varepsilon \Big) \Big\|_{W^{N_0 + 3(1-\frac{2}{p}), \frac{p}{p-1}}} \Bigg) \leq \delta_0,$$

*with $N_0 \geq 6$ and $N \geq N_0 + 14$, the system* (1.7) *admits a unique global solution such that $(\rho - \bar{\rho}, |\nabla|^s \phi, u) \in C([0, \infty); H^N)$. Moreover, the following time-decay estimate holds uniformly for $0 < \varepsilon \leq 1$,*

$$\|(\rho - \bar{\rho}, |\nabla|^s \phi, u)(t)\|_{W^{N_0, p}} \leq C \delta_0 (1 + t)^{-\beta_s}, \quad \forall\, t \geq 0,$$

*where the decay exponent $\beta_s$ is defined as*

$$\beta_s := \begin{cases} \dfrac{3}{2} \left( 1 - \dfrac{2}{p} \right) & \text{when } 0 < s \leq 1/2, \\ \dfrac{4}{3} \left( 1 - \dfrac{2}{p} \right) & \text{when } 1/2 < s < 1. \end{cases} \tag{1.9}$$

*In addition, if we assume that for any $0<\vartheta<3/2$, $\sup_{\varepsilon\in(0,1]}\||\nabla|^{-\vartheta}(\rho_0^\varepsilon-\bar\rho,|\nabla|^s\phi_0^\varepsilon,\mathcal{P}^\perp u_0^\varepsilon)\|_{L^2}$ is bounded, then the solution satisfies the following $\varepsilon$-dependent time-decay estimates*

$$\begin{aligned}&\|(\rho-\bar\rho)(t)\|_{L^2}\le C(1+\varepsilon^\Upsilon t)^{-\frac{s}{2}-\frac{\vartheta}{2}},\quad \|\nabla(\rho-\bar\rho)(t)\|_{H^{N-1}}\le C(1+\varepsilon^\Upsilon t)^{-\frac12-\frac{\vartheta}{2}},\\ &\|\nabla^k(|\nabla|^s\phi,u)(t)\|_{H^{N-k}}\le C(1+\varepsilon^\Upsilon t)^{-\frac{k}{2}-\frac{\vartheta}{2}},\quad \text{for } k=0,1,\end{aligned}$$

*where the exponent $\Upsilon\ge1$ corresponds to $\Upsilon_1$ (defined in (2.91)) for $0<s\le1/2$, and to $\Upsilon_2$ (defined in (2.96)) for $1/2<s<1$.*

**Remark 1.4.** *It should be noted that the finiteness of the $L^1$ norm of initial data implies that the perturbation is electrically neutral. However, for the case $1/2<s<1$, this assumption can be dropped. In fact, the $L^1$ boundedness is used to control the contribution of the viscous term. For the precise details, we refer to (2.80).*

**Remark 1.5.** *Although it seems reasonable to expect that the Riesz potential serves as a good approximation of the Coulomb as $s\to1$, the time-decay rate for $s\in(1/2,1)$ is slower than the case $s=1$. More interestingly, the time-decay rate for $s\in(0,1/2]$ is the same as the case $s=1$, we refer the reader to the estimates (2.12), (2.42). This phenomenon was also observed by Choi, Jung, and Lee in [12] for the compressible Euler-Riesz system.*

**Remark 1.6.** *The regularity requirements for initial data in $W^{N_0,p}$ and $H^N$ heavily depend on the estimates of the phase function $\Phi$ as shown in Lemma A.2. Similar estimates were also obtained in [12]. By carrying out a similar argument, we obtain a slightly improved result on $\Phi$, which weakens the assumptions on the regularity of initial data.*

To study the large-time behavior, we introduce a new unknown as follows

$$V:=\left(\frac{\sqrt{-\Delta+(-\Delta)^{1-s}}}{|\nabla|}(\rho-\bar\rho),\frac{\mathrm{div}}{|\nabla|}u\right),$$

then by Duhamel's principle, we obtain the following mild formulation for $V$

$$V=e^{-tA(D)}V_0+\int_0^t e^{-(t-\tau)A(D)}B(V)(\tau)\,\mathrm{d}\tau,$$

where $A$ is a linear operator defined in (2.4) and $B$ is a nonlinear term defined in (2.5). It is expected that the large-time behavior of the solution mainly depends on $e^{-tA}$. After computing the eigenvalues of $A$

$$\lambda_\pm(\xi):=-\varepsilon|\xi|^2\pm\mathrm{i}\sqrt{|\xi|^2+|\xi|^{2-2s}-\varepsilon^2|\xi|^4},$$

we conclude that, on the one hand, the system should exhibit a $L^p$ ($p>2$) time-decay estimate due to the dispersion if $\varepsilon|\xi|^2\le2\kappa_0$. On the other hand, the system should exhibit a $L^2$ time-decay estimate due to the dissipation if $\varepsilon|\xi|^2\ge2\kappa_0$. Here, $\kappa_0$ is a constant to be chosen small but independent of $\varepsilon$. Therefore, we consider two different regions, i.e., the dispersion-dominated domain $\Omega^L$ and the dissipation-dominated domain $\Omega^H$. In $\Omega^L$, we prove Proposition 2.1 and Corollary 2.7, which give the linear time-decay estimates of $e^{\mathrm{i}tb}$ uniformly in $\varepsilon$ with

$$b(\xi):=\sqrt{|\xi|^2+|\xi|^{2-2s}-\varepsilon^2|\xi|^4}$$

Since the Coulomb potential gives

$$b(\xi):=\sqrt{1+|\xi|^2-\varepsilon^2|\xi|^4},$$

the Riesz potential behaves much closer to the wave dispersion $\omega(\xi):=|\xi|$ than the Klein-Gordon dispersion $\omega(\xi):=\sqrt{1+|\xi|^2}$ as $s\in(0,1)$. The proof of Proposition 2.1 relies on the stationary phase analysis and the properties of $b(\xi)$. The reason for the difference between the case $s\in(0,1/2]$ and the case $s \in (1/2, 1)$ is that the phase function determined by $b$ has a degeneracy when $s \in (1/2, 1)$. The detailed analysis of $b$ can be found in the Appendix A.

To construct the global solution as well as its $L^2 - L^\infty$ time-decay rates, we derive a series of estimates in Sobolev spaces.

First, we apply nonlinear energy estimates to the system and obtain the estimates in $H^N$ norms. Since we are interested in the case for a general adiabatic constant $\gamma > 1$, we introduce a new unknown as follows

$$\mathfrak{U} := \left( \frac{2}{\gamma - 1} (\rho^{\frac{\gamma-1}{2}} - \bar{\rho}^{\frac{\gamma-1}{2}}), u \right).$$

This change of variables transforms the system (1.7) into a symmetric hyperbolic-parabolic system (2.46). The basic $L^2$ energy estimate shows that we only need to focus on the following term

$$I_3 := \int_{\mathbb{R}^3} \nabla^k \phi \, \nabla^k \mathrm{div}\, u \, \mathrm{d}x.$$

For $s \in [1/2, 1)$, $I_3$ can be bounded directly using integration by parts and the Sobolev interpolation $\dot{H}^{k+s} \cap \dot{H}^{k+2s} \subset \dot{H}^{k+1}$. When $s \in (0, 1/2)$, the situation becomes delicate. Inspired by [12], our key observation is the following: We first notice by the continuity equation $(1.7)_1$

$$-\mathrm{div}\, u = \frac{\partial_t \rho + u \cdot \nabla \rho}{\rho},$$

then applying $|\nabla|^{-s}$ to the above equation and substituting $(1.7)_3$, we obtain

$$\begin{aligned} -\mathrm{div}\, \nabla^k |\nabla|^{-s} u = & \frac{\partial_t \nabla^k |\nabla|^{-s} \rho}{\rho} + \frac{u \cdot \nabla \nabla^k |\nabla|^s \phi}{\rho} \\ & + \left[ \nabla^k |\nabla|^{-s}, \frac{1}{\rho} \right] \mathrm{div}\, (\rho u) + \left[ \nabla^k |\nabla|^{-s}, \frac{u}{\rho} \right] \cdot \nabla \rho. \end{aligned}$$

Multiplying the above equation by $\nabla^k |\nabla|^s \phi$ and using integration by parts, the left-hand side of the equation gives $I_3$, while on the right-hand side, the first term gives

$$\frac{1}{2} \frac{d}{dt} \int_{\mathbb{R}^3} \frac{1}{\rho} |\nabla^k |\nabla|^s \phi|^2 \, \mathrm{d}x,$$

which is a modified energy, the second term can be bounded by

$$\int_{\mathbb{R}^3} \mathrm{div} \left( \frac{u}{\rho} \right) |\nabla^k |\nabla|^s \phi|^2 \, \mathrm{d}x \lesssim \|\nabla(\varrho, u)\|_{L^\infty} \|\nabla^k |\nabla|^s \phi\|_{L^2}^2.$$

The rest of the terms can be bounded by using Lemma B.4.

For the estimate in $W^{k,p}$ norm, the dissipation-dominated part can be directly controlled by applying the linear estimates of $e^{-tA}$ to the mild formulation (2.9). For the dispersion-dominated part, the situation becomes more complicated. By Duhamel's formula, we write

$$V^L = e^{-tA(D)} V^L(0) + \int_0^t e^{-(t-\tau)A(D)} \chi_{\varepsilon, \kappa_0}(D) \left( \mathcal{B}(V^H, V) + \mathcal{B}(V^L, V^H) + \mathcal{B}(V^L, V^L) + \mathcal{Q}(V) \right) \mathrm{d}\tau,$$

where the superscript $L$ denotes the dispersion-dominated part and $H$ denotes dissipation-dominated part. For the precise definition of $V^L$, $V^H$, $\mathcal{B}$, and $\mathcal{Q}$, we refer the reader to (2.70), (2.6), and (2.7), respectively. Due to the fact that the linear estimates in the dispersion-dominated region are not sufficient to control the quadratic nonlinearity $\mathcal{B}(V^L, V^L)$, a standard remedy is to apply Shatah's normal form method and the space-time resonance method based on Shatah [30], Germain, Masmoudi, and Shatah [18, 19]. We rewrite the term involving $\mathcal{B}(V^L, V^L)$ as the sum of $K_{jk}$, see (2.77). Then, we study the bilinear operator $T_{\frac{m}{\Phi_{jk}}}$ defined by (1.13) where $m$ and $\Phi_{jk}$ are given by (2.74) and (2.76). Different from the Coulomb potential case where "time resonance set" is empty in the dispersive-dominated region, $\Phi_{jk}$ may present a significant zero set $\{\eta = 0\} \cup \{\xi - \eta = 0\}$. To resolve the singularity induced by $\Phi_{jk}^{-1}$, we are inspired by [20, 12]. Let us take $\Phi_{11}$ for illustration, by investigating the behavior of $\Phi_{11}$ near the zero sets, we find that $\Phi_{11}(\xi, \eta) \simeq |\eta|^{1-s}$ as $\eta \simeq 0$ and $\Phi_{11}(\xi, \eta) \simeq |\xi - \eta|^{1-s}$ as $\eta \simeq \xi$ (see Lemma A.2), therefore we introduce the following modified kernel:

$$\mathfrak{M}_{jk} := \frac{|\xi|^{1-s}|\eta|^{1-s}|\xi-\eta|^{1-s}}{\langle \xi - \eta \rangle^{2\lambda} \langle \eta \rangle^{2\lambda} \Phi_{jk}} \widetilde{\chi}_{\varepsilon,\kappa_0}(\xi - \eta) \widetilde{\chi}_{\varepsilon,\kappa_0}(\eta) \widetilde{\chi}_{\varepsilon,\kappa_0}(\xi) \frac{a(\xi)}{2\mathrm{i}b(\xi)},$$

then the estimates for the bilinear operator $T_{\frac{m}{\Phi_{jk}}}$ are obtained by employing the estimates of the solution in the negative Sobolev space $\dot{H}^{-(1-s)}$ and $\mathcal{M}^{\ell,\infty}_{\xi,\eta} \cap \mathcal{M}^{\ell,\infty}_{\eta,\xi}$ estimates of the kernel $\mathfrak{M}$. The estimate in $W^{k,p}$ norm shows the optimal $L^\infty$ time-decay rate in the sense of inviscid limit.

To complete the argument, we prove some estimates in negative Sobolev spaces. As indicated in the above paragraph, we first obtain the uniform estimate in $\dot{H}^{-(1-s)}$ norm by direct computation using (2.9). Then, in order to derive the $L^2$ time-decay rate, different from the $W^{k,p}$ estimate, we have to investigate the parabolic structure of the compressible Navier-Stokes-Riesz system. Inspired by [35], we further focus on uniform estimates in $\dot{H}^{-\vartheta}$ for $0 < \vartheta < 3/2$. Since the $L^\infty$ time-decay rate is not strong enough to close the estimates for the whole range of $\vartheta$, we design an iterative procedure by using additional $L^2$ time-decay estimates to improve the estimates step by step.

We notice that the range of the index $p$ for the function space of the initial data is strongly influenced by the time-decay rate of the solution and the contribution of error terms created by the artificial viscosity.

Now, we further construct $(n, v)$ as the correction to $(\rho, u)$ through the following system:

$$\begin{cases} \partial_t n + \operatorname{div}(\rho v + nu + nv) = 0, \\ \partial_t v + u \cdot \nabla v + v \cdot (\nabla u + \nabla v) + \dfrac{1}{\gamma - 1} \nabla \Big[ (\rho + n)^{\gamma-1} - \rho^{\gamma-1} \Big] + \nabla \psi - \varepsilon \mathcal{L} v = \varepsilon \left( \dfrac{1}{\rho + n} - 1 \right) (\mathcal{L} v + \mathcal{L} u), \\ (-\Delta)^s \psi = n, \\ (n, v)(0, x) = (0, \mathcal{P} u_0^\varepsilon)(x). \end{cases} \tag{1.10}$$

For (1.10), we prove the following result:

**Theorem 1.7.** *Consider $(\rho, u)$ and $\delta_0$ given by Theorem 1.3. Let $3 \le N_1 \le N_0 - 3$, if $\|\mathcal{P} u_0^\varepsilon\|_{H^{N_1}} \le \delta_0 \varepsilon$, then the system (1.10) has a unique solution in $C([0, \infty), H^{N_1})$ with the uniform estimate*

$$\sup_{0 \le t < +\infty} \|(n, |\nabla|^s \psi, v)(t)\|_{H^{N_1}} \le C \delta_0 \varepsilon.$$

*Moreover, if we assume in addition that for any $0 < \vartheta \le 1/2$, $\sup_{\varepsilon \in (0,1]} \varepsilon^{-1} \| |\nabla|^{-\vartheta} \mathcal{P} u_0^\varepsilon \|_{L^2}$ is bounded, then the solution satisfies the following time-decay estimate*

$$\|(n, |\nabla|^s \psi, v)(t)\|_{H^{N_1}} \le C \min\left\{ \varepsilon, (1+t)^{-\frac{\vartheta}{2}} \right\}.$$

**Remark 1.8.** *It is desired to extend the result for $1/2 < \vartheta < 3/2$. However, due to a technical issue, which shall be explained later, our analysis breaks down in such a scenario. Specifically, the lack of sharper decay estimates for the higher-order derivatives of the correction system (1.10) (similar to (2.51) in Proposition 2.8) prevents us from iteratively enhancing the $L^2$ decay rate through a bootstrap argument similar to Theorem 1.3.*

To obtain the uniform estimate for $(n, v)$, we will apply similar nonlinear energy estimates as in the analysis of $(\rho, u)$. Taking the contribution of the approximation solution into account, we introduce a new variable as follows

$$\mathfrak{V} := \left( \frac{2}{\gamma - 1} \left[ (\rho + n)^{\frac{\gamma-1}{2}} - \rho^{\frac{\gamma-1}{2}} \right], v \right).$$

With this specific transformation, the new unknown $\mathfrak{V}$ again satisfies a symmetric hyperbolic-parabolic system (3.1).

Different from the Lyapunov inequalities for energy of the approximation solution, here, we lack a similar estimate as (2.51) for purely high-order energy of the correction solution, which precludes the faster decay framework of Theorem 1.3. Specifically, to control coupling terms involving the interaction between the approximation solution and the correction solution such as $\int_{\mathbb{R}^3} v \cdot \nabla \nabla^k \mathfrak{c} \nabla^k \mathfrak{n} \, \mathrm{d}x$ ($k \geq 1$), since the limit profile is $\mathcal{O}(1)$, this term cannot be absorbed by the $\varepsilon$-scaled dissipation. Alternatively, if one expects to bound it via Gronwall's inequality, it turns out that by exploiting the decay of $\|\nabla^{k+1}\mathfrak{c}\|_{L^q}$ ($q > 6$), one needs to utilize $\|v\|_{L^r}$ ($r < 3$), which inevitably introduces the zero-order norm $\|v\|_{L^2}$ after applying Sobolev embeddings.

The restriction $0 < \vartheta \leq 1/2$ stems from the uniform estimate in negative Sobolev spaces. Indeed, since we need to control $\|(\mathfrak{n}, v)\|_{L^{\frac{3}{\vartheta}}}$ and $\||\nabla|^s \psi\|_{\dot{H}^{\frac{3}{2}-\vartheta}}$ in linear terms, we require $3/\vartheta \geq 6$, $3/2 - \vartheta \geq 1$, which ultimately imposes $\vartheta \leq 1/2$.

Having obtained Theorem 1.3 and Theorem 1.7, the solution to the compressible Navier-Stokes-Riesz system (1.1)-(1.3) with (1.4) can be obtained by

$$(\rho^\varepsilon, u^\varepsilon) = (\rho, u) + (n, v),$$

and we conclude the precise statement of Theorem 1.1 from Theorem 1.3 and Theorem 1.7:

**Theorem 1.9.** *For $0 < s < 1$, let $p$ be chosen according to*

$$\begin{cases} 6 < p < 12 & \text{when } 0 < s \leq 1/2, \\ 8 < p < 16 & \text{when } 1/2 < s < 1. \end{cases}$$

*Then there exists $\delta_0 > 0$ such that for any initial data satisfying*

$$\sup_{\varepsilon \in (0,1]} \left( \left\| \left( \rho_0^\varepsilon - \bar{\rho}, |\nabla|^s \phi_0^\varepsilon, \mathcal{P}^\perp u_0^\varepsilon \right) \right\|_{\dot{H}^{-(1-s)} \cap H^N} + \left\| \left( \rho_0^\varepsilon - \bar{\rho}, |\nabla|^s \phi_0^\varepsilon, \mathcal{P}^\perp u_0^\varepsilon \right) \right\|_{L^1} \right.$$
$$\left. + \left\| \left( \rho_0^\varepsilon - \bar{\rho}, |\nabla|^s \phi_0^\varepsilon, \mathcal{P}^\perp u_0^\varepsilon \right) \right\|_{W^{N_0 + 3(1-\frac{2}{p}), \frac{p}{p-1}}} \right) + \varepsilon^{-1} \|\mathcal{P} u_0^\varepsilon\|_{H^{N_1}} \leq \delta_0,$$

*with $N_0 \geq 6$, $N \geq N_0 + 14$ and $3 \leq N_1 \leq N_0 - 3$, there exists a unique global solution to the system* (1.1)-(1.3) *with* (1.4) *such that $(\rho^\varepsilon - \bar{\rho}, |\nabla|^s \phi^\varepsilon, u^\varepsilon) \in C([0,\infty), H^{N_1})$. In addition, if we assume that for any $0 < \vartheta \leq 1/2$, $\sup_{\varepsilon \in (0,1]} \{ \||\nabla|^{-\vartheta} (\rho_0^\varepsilon - \bar{\rho}, |\nabla|^s \phi_0^\varepsilon, \mathcal{P}^\perp u_0^\varepsilon)\|_{L^2} + \varepsilon^{-1} \||\nabla|^{-\vartheta} \mathcal{P} u_0^\varepsilon\|_{L^2} \}$ is bounded, then the solution satisfies the following time-decay estimates*

$$\|(\rho^\varepsilon - \bar{\rho}, |\nabla|^s \phi^\varepsilon, u^\varepsilon)\|_{W^{N_1 - 2, \infty}} \leq C \Big( \min \Big\{ \varepsilon, (1+t)^{-\frac{\vartheta}{2}} \Big\} + (1+t)^{-\beta_s} \Big),$$

*and*

$$\|\rho^\varepsilon - \bar{\rho}\|_{H^{N_1}} \leq C \Big( \min \Big\{ \varepsilon, (1+t)^{-\frac{\vartheta}{2}} \Big\} + (1 + \varepsilon^\Upsilon t)^{-\frac{s}{2} - \frac{\vartheta}{2}} \Big),$$
$$\|(|\nabla|^s \phi^\varepsilon, u^\varepsilon)\|_{H^{N_1}} \leq C (1 + \varepsilon^\Upsilon t)^{-\frac{\vartheta}{2}},$$

*where the exponent $\Upsilon \geq 1$ corresponds to $\Upsilon_1$ (defined in* (2.91)*) for $0 < s \leq 1/2$, and to $\Upsilon_2$ (defined in* (2.96)*) for $1/2 < s < 1$, and $\beta_s$ is given by* (1.9).

**Remark 1.10.** *Since our result is valid for $\varepsilon = 1$, we extend the existence result of* [37] *to the full range of $0 < s < 1$. While the time-decay rate established in the present work is not sharp (the above $L^2$ time-decay rate is slower than* (1.6)*) due to the approach of construction for the solution, our analysis should be useful to complement the optimal time-decay result in* [37] *for the case $0 < s \leq 1/2$.*

Our last result justifies the inviscid limit from the compressible Navier-Stokes-Riesz system to the compressible Euler-Riesz system.

**Theorem 1.11.** *Under the assumption of Theorem 1.9, let* $(\widetilde{\rho}, \widetilde{u})$ *be the global solution to the compressible Euler-Riesz system with the initial data* $(\widetilde{\rho}_0, \widetilde{u}_0) := (\rho_0^\varepsilon, \mathcal{P}^\perp u_0^\varepsilon)$*. Then, for any integer* $3 \le N_1 \le N_0 - 3$, $3 \le N_0' \le N_1$ *and* $N_2 \le N_0' - 3/2(1-2/p)$*, we have*

$$\sup_{0 \le t < +\infty} \|(\rho^\varepsilon - \widetilde{\rho}, |\nabla|^s(\phi^\varepsilon - \widetilde{\phi}), u^\varepsilon - \widetilde{u})(t)\|_{W^{N_2,p}} \le C\varepsilon^{\lambda_s},$$

*where* $C$ *is a constant independent of* $t$ *and* $\varepsilon$ *and* $\lambda_s$ *is defined by*

$$\lambda_s := \begin{cases} \dfrac{3(p-2)}{5p-6} & \text{when } 0 < s \le 1/2, \\ \dfrac{4(p-2)}{7p-8} & \text{when } 1/2 < s < 1. \end{cases} \tag{1.11}$$

**Remark 1.12.** *For the initial data* $(\widetilde{\rho}_0, \widetilde{u}_0) := (\rho_0^\varepsilon, \mathcal{P}^\perp u_0^\varepsilon)$ *in the above theorem, it has been proved by* [12] *that the Cauchy problem for the compressible Euler-Riesz system admits a global solution* $(\tilde{\rho}, \tilde{u})$*. We also refer the reader to Lemma C.1 in the appendix.*

The rest of the paper is organized as follows: In Section 2, we study the viscous approximation system (1.7). First, we derive the linear decay estimates of the semigroup $e^{-tA}$, then by performing estimates in various Sobolev spaces, we prove Theorem 1.3. In Section 3, we study the correction system (1.10). In order to prove Theorem 1.7, some Lyapunov-type differential inequalities are established. In Section 4, we rigorously justify the inviscid limit from the compressible Navier-Stokes-Riesz system to the compressible Euler-Riesz system as $\varepsilon \to 0$. By establishing uniform energy estimates for the error system, we obtain an explicit convergence rate for all time $t \ge 0$. In Appendix A, we study the properties of the imaginary part of the eigenvalues of the semigroup $e^{-tA}$, together with lower and derivative estimates for the phase function $\Phi(\xi, \eta)$. In Appendix B, we collect useful inequalities that are frequently used in this paper, and establish the estimates for the kernel $\mathfrak{M}$ together with the corresponding bilinear estimates. Finally, in Appendix C, we recall the global well-posedness result for the compressible Euler-Riesz system, which is used in the inviscid limit.

*Notations.* Throughout this paper, we denote $\langle \cdot \rangle := \sqrt{1 + |\cdot|^2}$, $a^+$ a constant which is larger and arbitrarily close to $a$ and $\delta$ a sufficiently small positive constant. We write $A \lesssim B$ ($A \gtrsim B$) for $A \le CB$ ($A \ge C^{-1}B$) for $C > 0$ a generic number that can be chosen independent of $\varepsilon \in (0, 1]$ and $t > 0$. We also write $A \simeq B$ if $A \lesssim B$ and $A \gtrsim B$. We denote $\mathcal{F}$ and $\mathcal{F}^{-1}$ the Fourier transform and the inverse Fourier transform, respectively,

$$\mathcal{F}(f) := \int_{\mathbb{R}^d} f(x) e^{-\mathrm{i}x\cdot\xi}\,\mathrm{d}x, \quad \mathcal{F}^{-1}(f) := \frac{1}{(2\pi)^d} \int_{\mathbb{R}^d} f(\xi) e^{\mathrm{i}x\cdot\xi}\,\mathrm{d}\xi.$$

We define the homogeneous Sobolev space $\dot{H}^s$ ($s \in \mathbb{R}$) of all $f$ for which $\|f\|_{\dot{H}^s}$ is finite, where

$$\|f\|_{\dot{H}^s} := \||\xi|^s \widehat{f}\|_{L^2}.$$

We will also denote $\widehat{f} := \mathcal{F}(f)$. Then, we define $m(D)$ the Fourier multiplier by

$$m(D)f := \mathcal{F}^{-1}(m(\cdot)\widehat{f}(\cdot)). \tag{1.12}$$

Moreover, for a given function $m(\zeta, \eta)$, we define the bilinear operator $T_m(f, g)$ as:

$$T_m(f, g) := \mathcal{F}_\xi^{-1}\left(\int m(\xi, \eta)\widehat{f}(\xi - \eta)\widehat{g}(\eta)\,\mathrm{d}\eta\right). \tag{1.13}$$

In our analysis, we will need the Littlewood-Paley decomposition. We define:

$$\text{radially symmetric functions } \chi(\xi), \varphi_0(\xi) \in C_c^\infty \text{ with } \chi, \varphi_0 = \begin{cases} 1, & \text{on } \{\xi : |\xi| \le 1\}, \\ 0, & \text{in } \{\xi : |\xi| > 2\}, \end{cases} \tag{1.14}$$

and

$$\widetilde{\chi}(\xi) := \chi(\xi/2), \tag{1.15}$$

$$\chi_{\varepsilon,\kappa_0}(\xi) := \chi\big(\sqrt{\frac{\varepsilon}{\kappa_0}}\xi\big), \tag{1.16}$$

$$\varphi_{\eta,\delta}(\xi) := \varphi_0(\frac{\xi - \eta}{\delta}), \tag{1.17}$$

$$\varphi(\frac{\xi}{2^j}) := \varphi_0(\frac{\xi}{2^j}) - \varphi_0(\frac{\xi}{2^{j-1}}), \quad j \in \mathbb{Z}. \tag{1.18}$$

Then we define the inhomogeneous Besov space $B^s_{p,r}$ $(p, r \ge 1, s \in \mathbb{R})$ of all $f$ for which $\|f\|_{B^s_{p,r}}$ is finite, where

$$\|f\|_{B^s_{p,r}} := \left( \|\varphi_0(D)f\|^r_{L^p} + \sum_{j=1}^{\infty} 2^{jsr} \|\varphi(\frac{D}{2^j})f\|^r_{L^p} \right)^{\frac{1}{r}}. \tag{1.19}$$

In addition, we introduce the following multiplier norm

$$\|\mathfrak{M}\|_{\mathcal{M}^{k,\infty}_{\xi,\eta}} := \sum_{j \in \mathbb{Z}} \|\varphi(\frac{D_\eta}{2^j})\mathfrak{M}\|_{L^\infty_\xi \dot{H}^k_\eta}.$$

From now on, we set $\bar{\rho} = 1$, $\mu = 1$, $\nu = 0$, $K\gamma = 1$ for simplicity.

## 2. Viscous approximation of the compressible Euler-Riesz System

To study the system (1.7), we introduce the following new unknown functions

$$h := \frac{a(D)}{|\nabla|}\varrho, \quad q := \frac{\mathrm{div}}{|\nabla|}u, \quad \varrho := \rho - 1, \tag{2.1}$$

where

$$a(D) := \sqrt{-\Delta + (-\Delta)^{1-s}} \tag{2.2}$$

and we set $V := (h, q)^{\mathsf{T}}$, then a direct computation shows that $V$ satisfies the system

$$\partial_t V + A(D)V = B(V), \tag{2.3}$$

where $A(D)$ is the linear operator

$$A(D) := \begin{pmatrix} 0 & a(D) \\ -a(D) & -2\varepsilon\Delta \end{pmatrix}, \tag{2.4}$$

and $B(V)$ is the nonlinear term

$$B(V) := \mathcal{B}(V, V) + \mathcal{Q}(V), \tag{2.5}$$

$$\mathcal{B}(V, V) := \begin{pmatrix} -a(D)\mathcal{R}^* \left(\frac{|\nabla|}{a(D)} h \mathcal{R}q\right) \\ \mathcal{R}^*\nabla \left[\frac{1}{2}|\mathcal{R}q|^2 + \frac{\gamma-2}{2}\left|\frac{|\nabla|}{a(D)}h\right|^2\right] \end{pmatrix}, \tag{2.6}$$

$$\mathcal{Q}(V) := \begin{pmatrix} 0 \\ \mathcal{R}^*\nabla \left[\frac{1}{\gamma-1}(1 + \frac{|\nabla|}{a(D)}h)^{\gamma-1} - \frac{1}{\gamma-1} - \frac{|\nabla|}{a(D)}h - \frac{\gamma-2}{2}|\frac{|\nabla|}{a(D)}h|^2\right] \end{pmatrix}, \tag{2.7}$$

with the Riesz transform and its adjoint,

$$\mathcal{R} := \frac{\nabla}{|\nabla|}, \quad \mathcal{R}^* := -\frac{\mathrm{div}}{|\nabla|}. \tag{2.8}$$

Moreover, by the Duhamel principle, we rewrite (2.3) as

$$V = e^{-tA}V_0 + \int_0^t e^{-(t-s)A}B(V)(s)\,\mathrm{d}s, \tag{2.9}$$

where the Green matrix can be computed as

$$e^{-tA(\xi)} = \frac{1}{\lambda_+ - \lambda_-}\begin{pmatrix} \lambda_+ e^{\lambda_- t} - \lambda_- e^{\lambda_+ t} & (e^{\lambda_- t} - e^{\lambda_+ t})a(\xi) \\ (e^{\lambda_+ t} - e^{\lambda_- t})a(\xi) & \lambda_+ e^{\lambda_+ t} - \lambda_- e^{\lambda_- t}\end{pmatrix} := \begin{pmatrix} G_1(t,\xi) & -G_2(t,\xi) \\ G_2(t,\xi) & G_3(t,\xi)\end{pmatrix} \tag{2.10}$$

with the eigenvalues of $-A(\xi)$ defined by

$$\lambda_\pm(\xi) := -\varepsilon|\xi|^2 \pm \mathrm{i}\sqrt{|\xi|^2 + |\xi|^{2-2s} - \varepsilon^2|\xi|^4}.$$

We also introduce the following differential operator

$$b(D) := \sqrt{(-\Delta) + (-\Delta)^{1-s} - \varepsilon^2(-\Delta)^2}. \tag{2.11}$$

The mild formulation (2.9) motivates us to study the property of $e^{-tA}$. In subsection 2.1, we obtain some linear estimates for $e^{-tA}$. More precisely, because $e^{-tA}$ shows the dispersive nature and the dissipative nature at the same time, we discuss $e^{-tA}$ in two cases according to the strength of the scale $\varepsilon$ for viscosity. To prove Theorem 1.3, in subsection 2.2, we first apply the nonlinear energy method to derive the estimate in $H^N$ spaces. Then, we are inspired by the normal form method to obtain the desired estimates in $W^{k,p}$ spaces. To conclude the estimate, we prove some estimates in homogeneous negative Sobolev spaces.

2.1. **Linear Decay.** In this subsection, we will exploit the linear time-decay of $e^{-tA}$ in dispersion-dominated domain $\Omega^L := \{\xi : \varepsilon|\xi|^2 \le 2\kappa_0\}$ and dissipation-dominated domain $\Omega^H := \{\xi : \varepsilon|\xi|^2 \ge 2\kappa_0\}$, respectively. We shall carry out the analysis in a general spatial dimension $d \ge 2$ over $\mathbb{R}^d$, although in this paper, we only utilize the case $d = 3$.

2.1.1. *Estimates in $\Omega^L$.* The main results in this part are stated as follows.

**Proposition 2.1.** *There exists $\kappa_0 > 0$ small enough such that uniformly for $\varepsilon \in (0,1]$, and for every $f \in B^d_{1,2}$, we have the estimate*

$$\|e^{\mathrm{i}tb(D)}\chi_{\varepsilon,\kappa_0}(D)f\|_{B^0_{\infty,2}} \lesssim (1+t)^{-\alpha_{d,s}}\|f\|_{B^d_{1,2}}, \tag{2.12}$$

*Moreover,*

$$\|G_j(t,D)\chi_{\varepsilon,\kappa_0}(D)f\|_{B^0_{\infty,2}} \lesssim (1+t)^{-\alpha_{d,s}}\|f\|_{B^d_{1,2}}. \tag{2.13}$$

*Here, the cutoff function $\chi_{\varepsilon,\kappa_0}$ is defined in* (1.16) *and*

$$\alpha_{d,s} := \begin{cases} \dfrac{d}{2}, & \text{when } 0 < s \le 1/2, \\ \dfrac{d}{2} - \dfrac{1}{6}, & \text{when } 1/2 < s < 1. \end{cases}$$

The proof of Proposition 2.1 will follow from the next two lemmas. The first result concerns the case for $s \in (0, 1/2]$ where the estimates are similar to the ones of the linear Klein-Gordon equation.

**Lemma 2.2.** *For $s \in (0, 1/2]$ and $\lambda \in [1, \infty)$, there exists $\kappa_0 > 0$ small enough such that uniformly for $\varepsilon \in (0, 1]$, the following estimates hold*

$$\|e^{\mathrm{i}tb(D)}\chi_{\varepsilon,\kappa_0}(D)\varphi_0(D)f\|_{L^\infty} \lesssim (1+t)^{-\frac{d}{2}}\|f\|_{L^1}, \tag{2.14}$$

$$\|e^{\mathrm{i}tb(D)}\chi_{\varepsilon,\kappa_0}(D)\varphi(\frac{D}{\lambda})f\|_{L^\infty} \lesssim (1+t)^{-\frac{d}{2}}\lambda^d\|f\|_{L^1}, \tag{2.15}$$

*where the cutoff functions $\chi_{\varepsilon,\kappa_0}$, $\varphi_0$, and $\varphi$ are defined as* (1.16), (1.14), *and* (1.18), *respectively.*

*Proof.* First, let us prove (2.14). According to the definition of Fourier multipliers, by Young's convolution inequality, we have

$$\|e^{\mathrm{i}tb(D)}\chi_{\varepsilon,\kappa_0}(D)\varphi_0(D)f\|_{L^\infty} \leq \|\mathcal{F}^{-1}[e^{\mathrm{i}tb(\xi)}\chi_{\varepsilon,\kappa_0}(\xi)\varphi_0(\xi)]\|_{L^\infty}\|f\|_{L^1}.$$

For $0 < t < 1$, we have

$$\|\mathcal{F}^{-1}[e^{\mathrm{i}tb(\xi)}\chi_{\varepsilon,\kappa_0}(\xi)\varphi_0(\xi)]\|_{L^\infty} \lesssim \|\varphi_0\|_{L^1},$$

therefore

$$\|e^{\mathrm{i}tb(D)}\chi_{\varepsilon,\kappa_0}(D)\varphi_0(\xi)f\|_{L^\infty} \leq \|f\|_{L^1}. \tag{2.16}$$

For $t \geq 1$, since $\chi_{\varepsilon,\kappa_0}$ and $\varphi_0$ are chosen to be spherically symmetric, it is known that

$$\mathcal{F}^{-1}[e^{\mathrm{i}tb(\xi)}\chi_{\varepsilon,\kappa_0}(\xi)\varphi_0(\xi)] = \int_0^2 e^{\mathrm{i}tb(r)}\chi_{\varepsilon,\kappa_0}(r)\varphi_0(r)\mathcal{F}(\sigma_{\mathbb{S}^{d-1}})(|x|r)r^{d-1}\,\mathrm{d}r,$$

where

$$\mathcal{F}(\sigma_{\mathbb{S}^{d-1}})(x) := |x|^{-\frac{d-2}{2}}J_{\frac{d-2}{2}}(|x|) = e^{\mathrm{i}|x|}Z(|x|) - e^{-\mathrm{i}|x|}\bar{Z}(|x|)$$

is the Fourier transform of the Lebesgue measure on the sphere $\mathbb{S}^{d-1}$, $J_{\frac{d-2}{2}}$ is the Bessel function, $Z(r)$ satisfies

$$|\partial_r^k Z(r)| \lesssim (1+r)^{-\frac{d-1}{2}-k}, \tag{2.17}$$

for all $k \geq 0$ and all $r > 0$, and $\bar{Z}$ is the complex conjugate of $Z$. Therefore, we write

$$\begin{aligned}\mathcal{F}^{-1}[e^{\mathrm{i}tb(\xi)}\chi_{\varepsilon,\kappa_0}(\xi)\varphi_0(\xi)] &= \sum_{\pm}\int_0^2 e^{\mathrm{i}t\Phi_\pm(r,t,x)}\chi_{\varepsilon,\kappa_0}(r)\varphi_0(r)Z_\pm(|x|r)r^{d-1}\,\mathrm{d}r \\ &:= I_+ + I_-,\end{aligned} \tag{2.18}$$

where we denote

$$\Phi_\pm(r,t,x) := b(r) \pm r\frac{|x|}{t}, \quad Z_+ := Z, \quad Z_- := \bar{Z}. \tag{2.19}$$

For $I_+$, using integration by parts $N_{d,s} := [\frac{d}{1-s}]$ times, we have

$$\begin{aligned}I_+ =& \int_0^2 \frac{\partial_r e^{\mathrm{i}t\Phi_+(r,t,x)}}{\mathrm{i}t\partial_r\Phi_+(r,t,x)}\chi_{\varepsilon,\kappa_0}(r)\varphi_0(r)Z_+(|x|r)r^{d-1}\,\mathrm{d}r \\ =& -\int_0^2 e^{\mathrm{i}t\Phi_+(r,t,x)}\partial_r\left(\frac{1}{\mathrm{i}t\partial_r\Phi_+(r,t,x)}\chi_{\varepsilon,\kappa_0}(r)\varphi_0(r)Z_+(|x|r)r^{d-1}\right)\mathrm{d}r \\ \simeq& (-\frac{1}{\mathrm{i}t})^{N_{d,s}}\sum_{m=0}^{N_{d,s}}\sum_{l_1,\ldots,l_{N_{d,s}}\in\Lambda_m^{N_{d,s}}}\int_0^2 e^{\mathrm{i}t\Phi_+(r,t,x)}\prod_{k=1}^{N_{d,s}}\partial_r^{l_k}\left(\frac{1}{\partial_r\Phi_+(r,t,x)}\right)\partial_r^{N_{d,s}-m}(\chi_{\varepsilon,\kappa_0}(r)\varphi_0(r)Z_+(|x|r)r^{d-1})\,\mathrm{d}r,\end{aligned} \tag{2.20}$$

where $\Lambda_m^{N_{d,s}} := \{l_1, ..., l_{N_{d,s}} \in \mathbb{N} : 0 \le l_1 < \cdots < l_{N_{d,s}} \le N_{d,s}, l_1 + \cdots + l_{N_{d,s}} = N_{d,s}\}$. Since for $r \in (0,2)$, by the definition in (2.19), (2.17), and Lemma A.1, we deduce that

$$\left|\partial_r^n \left(\frac{1}{\partial_r \Phi_+(r,t,x)}\right)\right| \lesssim r^{s-n}, \quad \left|\partial_r^n (\chi_{\varepsilon,\kappa_0}(r)\varphi_0(r) Z_+(|x|r) r^{d-1})\right| \lesssim r^{d-1-n}.$$

Therefore, we conclude

$$|I_+| \lesssim t^{-N_{d,s}}. \tag{2.21}$$

For $I_-$, if $|x|/t < 10^{-2}$, following a similar procedure as (2.20), together with the following estimates deduced by Lemma A.1,

$$\left|\partial_r^\alpha \left(\frac{1}{\partial_r \Phi_-(r,t,x)}\right)\right| \lesssim r^{s-\alpha}, \quad \left|\partial_r^\alpha (\chi_{\varepsilon,\kappa_0}(r)\varphi_0(r) Z_-(|x|r) r^{d-1})\right| \lesssim r^{d-1-\alpha},$$

we obtain

$$|I_-| \lesssim t^{-N_{d,s}}. \tag{2.22}$$

If $|x|/t \ge 10^{-2}$, then by Lemma A.1, $|b''(r)| \ge C^{-1}$ for $r \in [0,2]$, therefore by Van der Corput Lemma (Lemma B.1), we have

$$\begin{aligned} |I_-| &\lesssim t^{-\frac{1}{2}} \sup_{0\le r\le 2} \partial_r \left(\chi_{\varepsilon,\kappa_0}(r)\varphi_0(r) Z_-(|x|r) r^{d-1}\right) \\ &\lesssim t^{-\frac{1}{2}} (1+|x|)^{-\frac{d-1}{2}} \\ &\lesssim t^{-\frac{d}{2}}. \end{aligned} \tag{2.23}$$

Therefore, the estimate (2.14) is proved by combining (2.16), (2.21), (2.22), and (2.23).

It suffices to prove (2.15). According to the definition of Fourier multipliers, by Young's convolution inequality, we have

$$\|e^{\mathrm{i}tb(D)} \chi_{\varepsilon,\kappa_0}(D)\varphi(\frac{D}{\lambda}) f\|_{L^\infty} \le \|\mathcal{F}^{-1}[e^{\mathrm{i}tb(\xi)} \chi_{\varepsilon,\kappa_0}(\xi)\varphi(\frac{\xi}{\lambda})]\|_{L^\infty} \|f\|_{L^1}.$$

For $0 < t < 1$, by the change of variable, we have

$$\mathcal{F}^{-1}[e^{\mathrm{i}tb(\xi)} \chi_{\varepsilon,\kappa_0}(\xi)\varphi(\frac{\xi}{\lambda})] = \int_{\mathbb{R}^d} e^{\mathrm{i}tb(\xi)} e^{\mathrm{i}x\cdot\xi} \chi_{\varepsilon,\kappa_0}(\xi)\varphi(\frac{\xi}{\lambda})\, \mathrm{d}\xi = \lambda^d \int_{\mathbb{R}^d} e^{\mathrm{i}tb(\lambda\xi)} e^{\mathrm{i}\lambda x\cdot\xi} \chi_{\varepsilon,\kappa_0}(\lambda\xi)\varphi(\xi)\, \mathrm{d}\xi \lesssim \lambda^d.$$

therefore

$$\|e^{\mathrm{i}tb(D)} \chi_{\varepsilon,\kappa_0}(D)\varphi(\frac{D}{\lambda}) f\|_{L^\infty} \lesssim \lambda^d \|f\|_{L^1} \tag{2.24}$$

Similarly to (2.18), since $\chi_{\varepsilon,\kappa_0}$ and $\varphi_0$ are chosen to be spherically symmetric, we write

$$\begin{aligned} \mathcal{F}^{-1}[e^{\mathrm{i}tb(\xi)} \chi_{\varepsilon,\kappa_0}(\xi)\varphi(\frac{\xi}{\lambda})] &= \int_0^\infty e^{\mathrm{i}tb(r)} \chi_{\varepsilon,\kappa_0}(r)\varphi(\frac{r}{\lambda}) \mathcal{F}(\sigma_{\mathbb{S}^{d-1}})(|x|r) r^{d-1}\, \mathrm{d}r \\ &= \lambda^d \int_1^4 e^{\mathrm{i}tb(\lambda r)} \chi_{\varepsilon,\kappa_0}(\lambda r)\varphi(r) \mathcal{F}(\sigma_{\mathbb{S}^{d-1}})(\lambda|x|r) r^{d-1}\, \mathrm{d}r \\ &= \sum_{\pm} \lambda^d \int_1^4 e^{\mathrm{i}t\Phi_\pm(\lambda r)} \chi_{\varepsilon,\kappa_0}(\lambda r)\varphi(r) Z_\pm(\lambda|x|r) r^{d-1}\, \mathrm{d}r \\ &:= \lambda^d (J_+ + J_-). \end{aligned}$$

For $J_+$, using integration by parts $N$ times, we have

$$\begin{aligned}
J_+ &= \int_1^4 \frac{\partial_r e^{\mathrm{i}t\Phi_+(\lambda r,t,x)}}{\mathrm{i}\lambda t \partial_r \Phi_+(\lambda r,t,x)} \chi_{\varepsilon,\kappa_0}(\lambda r)\varphi(r) Z_+(\lambda|x|r) r^{d-1}\,\mathrm{d}r \\
&= \int_1^4 e^{\mathrm{i}t\Phi_+(\lambda r,t,x)} \partial_r \left( \frac{1}{\mathrm{i}\lambda t \partial_r \Phi_+(\lambda r,t,x)} \chi_{\varepsilon,\kappa_0}(\lambda r)\varphi(r) Z_+(\lambda|x|r) r^{d-1} \right) \mathrm{d}r \\
&\simeq \left(\frac{1}{\mathrm{i}\lambda t}\right)^N \sum_{m=0}^{N} \sum_{l_1,\dots,l_N\in\Lambda_m^N} \int_1^4 e^{\mathrm{i}\lambda t\Phi_+(\lambda r,t,x)} \prod_{k=1}^{N} \partial_r^{l_k} \left( \frac{1}{\partial_r\Phi_+(\lambda r,t,x)} \right) \partial_r^{N-m} (\chi_{\varepsilon,\kappa_0}(\lambda r)\varphi(r) Z_+(\lambda|x|r) r^{d-1})\,\mathrm{d}r,
\end{aligned} \tag{2.25}$$

where $\Lambda_m^N := \{l_1,...,l_N \in \mathbb{N} : 0 \le l_1 < \cdots < l_N \le N, l_1 + \cdots + l_N = m\}$. Since for $r \in (1,4)$, by the definition (2.19), (2.17), and Lemma A.1, we deduce that

$$\left|\partial_r^\alpha \left( \frac{1}{\partial_r\Phi_+(\lambda r,t,x)} \right)\right| \lesssim \begin{cases} 1, & n = 0, \\ \lambda^{-2s}, & n \ge 1, \end{cases} \quad \left|\partial_r^\alpha (\chi_{\varepsilon,\kappa_0}(\lambda r)\varphi(r) Z_+(\lambda|x|r) r^{d-1})\right| \lesssim \lambda^\alpha.$$

Therefore, we conclude

$$|J_+| \lesssim t^{-N}. \tag{2.26}$$

For $J_-$, if $|x|/t < 10^{-2}$, following a similar procedure as (2.25), together with the following estimates deduced by Lemma A.1

$$\left|\partial_r \left( \frac{1}{\partial_r\Phi_-(\lambda r,t,x)} \right)\right| \lesssim \begin{cases} 1, & \alpha = 0, \\ \lambda^{-2s}, & \alpha \ge 1, \end{cases} \quad |\partial_r^\alpha (\chi_{\varepsilon,\kappa_0}(\lambda r)\varphi(r) Z_-(\lambda|x|r) r^{d-1})| \le \lambda^\alpha,$$

we obtain

$$|J_-| \lesssim t^{-N}. \tag{2.27}$$

If $|x|/t \ge 10^{-2}$, then by the definition (2.19) and Lemma A.1,

$$\left|\partial_r^2 \Phi_-(\lambda r,t,x)\right| \gtrsim \begin{cases} \lambda^{1-2s}, & 0 < s < 1/2, \\ \lambda^{-1}, & s = 1/2, \end{cases}$$

by Van der Corput Lemma (Lemma B.1), for $0 < s < 1/2$,

$$\begin{aligned}
|J_-| &\lesssim (\lambda^{1-2s} t)^{-\frac{1}{2}} \sup_{1\le r\le 4} |\partial_r (\chi_{\varepsilon,\kappa_0}(\lambda r)\varphi(r) Z_-(\lambda|x|r) r^{d-1})| \\
&\lesssim (\lambda^{1-2s} t)^{-\frac{1}{2}} (1+\lambda|x|)^{-\frac{d-1}{2}} \\
&\lesssim \lambda^{-\frac{d}{2}+s} t^{-\frac{d}{2}}.
\end{aligned} \tag{2.28}$$

And for $s = 1/2$,

$$\begin{aligned}
|J_-| &\lesssim (\lambda^{-1} t)^{-\frac{1}{2}} \sup_{1\le r\le 4} |\partial_r (\chi_{\varepsilon,\kappa_0}(\lambda r)\varphi(r) Z_-(\lambda|x|r) r^{d-1})| \\
&\lesssim (\lambda^{-1} t)^{-\frac{1}{2}} (1+\lambda|x|)^{-\frac{d-1}{2}} \\
&\lesssim \lambda^{-\frac{d}{2}+1} t^{-\frac{d}{2}}.
\end{aligned} \tag{2.29}$$

Therefore, the estimate (2.15) is proved by combining (2.24), (2.26), (2.27), (2.28), and (2.29). □

The second result concerns the case for $s \in (1/2, 1)$ where the estimates are much closer to the ones of the linear wave equation than the ones of the linear Klein-Gordon equation.

**Lemma 2.3.** *For $1/2 < s < 1$ and for every $\lambda \geq 1$, there exist $\kappa_0, \delta > 0$ small enough such that uniformly for $\varepsilon \in (0,1]$, the following estimates hold*

$$\|e^{\mathrm{i}tb(D)}\chi_{\varepsilon,\kappa_0}(D)\tilde{\varphi}_{r_0,\delta}(D)f\|_{L^\infty} \lesssim (1+t)^{-\frac{d}{2}+\frac{1}{6}}\|f\|_{L^1}, \tag{2.30}$$

$$\|e^{\mathrm{i}tb(D)}\chi_{\varepsilon,\kappa_0}(D)\varphi_0(D)(1-\tilde{\varphi}_{r_0,\delta}(D))f\|_{L^\infty} \lesssim (1+t)^{-\frac{d}{2}}\|f\|_{L^1}, \tag{2.31}$$

$$\|e^{\mathrm{i}tb(D)}\chi_{\varepsilon,\kappa_0}(D)\varphi(\frac{D}{\lambda})(1-\tilde{\varphi}_{r_0,\delta}(D))f\|_{L^\infty} \lesssim (1+t)^{-\frac{d}{2}}\lambda^d\|f\|_{L^1}, \tag{2.32}$$

*where $r_0$ is the zero of $b''(r)$ determined by Lemma A.1 if it exists, otherwise, we set $r_0 := 0$. The cutoff functions $\chi_{\varepsilon,\kappa_0}$, $\varphi_{\eta,\delta}$, $\varphi_0$, and $\varphi$ are defined as (1.16), (1.17), (1.14), and (1.18), respectively.*

*Proof.* It suffices to consider the case where $b''$ admits a zero. First, we prove (2.30). According to the definition of Fourier multipliers, by Young's convolution inequality, we have

$$\|e^{\mathrm{i}tb(D)}\chi_{\varepsilon,\kappa_0}(D)\tilde{\varphi}_{r_0,\delta}(D)f\|_{L^\infty} \leq \|\mathcal{F}^{-1}[e^{\mathrm{i}tb(\xi)}\chi_{\varepsilon,\kappa_0}(\xi)\tilde{\varphi}_{r_0,\delta}(\xi)]\|_{L^\infty}\|f\|_{L^1}.$$

For $0 < t < 1$, we have

$$\|\mathcal{F}^{-1}[e^{\mathrm{i}tb(\xi)}\chi_{\varepsilon,\kappa_0}(\xi)\tilde{\varphi}_{r_0,\delta}(\xi)]\|_{L^\infty} \lesssim \|\varphi_0\|_{L^1},$$

then

$$\|e^{\mathrm{i}tb(D)}\chi_{\varepsilon,\kappa_0}(D)\tilde{\varphi}_{r_0,\delta}(\xi)f\|_{L^\infty} \leq \|f\|_{L^1}. \tag{2.33}$$

For $t \geq 1$, as in (2.18), since $\chi_{\varepsilon,\kappa_0}$ and $\varphi_0$ are chosen to be spherically symmetric, we write

$$\begin{aligned}
\mathcal{F}^{-1}[e^{\mathrm{i}tb(\xi)}\chi_{\varepsilon,\kappa_0}(\xi)\tilde{\varphi}_{r_0,\delta}(\xi)] &= \sum_{\pm}\int_0^\infty e^{\mathrm{i}t\Phi_\pm(r,t,x)}\chi_{\varepsilon,\kappa_0}(r)\varphi_0(\frac{r-r_0}{\delta})Z_\pm(|x|r)r^{d-1}\,\mathrm{d}r \\
&= \sum_{\pm}\int_{-2\delta}^{2\delta} e^{\mathrm{i}t\Phi_\pm(r_0+r,t,x)}\chi_{\varepsilon,\kappa_0}(r_0+r)\varphi_0(\frac{r}{\delta})Z_\pm(|x|(r_0+r))(r_0+r)^{d-1}\,\mathrm{d}r \\
&:= I_+ + I_-.
\end{aligned}$$

For $I_+$, using integration by parts $N \in \mathbb{N}$ times,

$$\begin{aligned}
I_+ &= \int_{-2\delta}^{2\delta} \frac{\partial_r e^{\mathrm{i}t\Phi_+(r_0+r,t,x)}}{\mathrm{i}t\partial_r\Phi_+(r_0+r,t,x)}\chi_{\varepsilon,\kappa_0}(r_0+r)\varphi_0(\frac{r}{\delta})Z_+(|x|(r_0+r))(r_0+r)^{d-1}\,\mathrm{d}r \\
&= -\int_{-2\delta}^{2\delta} e^{\mathrm{i}t\Phi_+(r_0+r,t,x)}\partial_r\left(\frac{1}{\mathrm{i}t\partial_r\Phi_+(r_0+r,t,x)}\chi_{\varepsilon,\kappa_0}(r_0+r)\varphi_0(\frac{r}{\delta})Z_+(|x|(r_0+r))(r_0+r)^{d-1}\right)\mathrm{d}r \\
&\simeq (-\frac{1}{\mathrm{i}t})^N \sum_{m=0}^N \sum_{l_1+\cdots+l_N=m} \\
&\quad \int_{-2\delta}^{2\delta} e^{\mathrm{i}t\Phi_+(r_0+r)}\prod_{k=1}^N \partial_r^{l_k}\left(\frac{1}{\partial_r\Phi_+(r_0+r,t,x)}\right)\partial_r^{N-m}(\chi_{\varepsilon,\kappa_0}(r_0+r)\varphi_0(\frac{r}{\delta})Z_+(|x|(r_0+r))(r_0+r)^{d-1})\,\mathrm{d}r.
\end{aligned} \tag{2.34}$$

Since for $r \in (0,2)$, by Lemma A.1, we deduce that

$$\left|\partial_r^n\left(\frac{1}{\partial_r\Phi_+(r,t,x)}\right)\right| \lesssim \sum_{m=1}^n \frac{1}{|\partial_r b(r_0+r)|^{m+1}} \sum_{\substack{k_1+\cdots+k_j=m,\\ k_1+\cdots+jk_j=n}} \prod_{i=1}^j |\partial_r^{i+1}b(r_0+r)|^{k_i} \lesssim 1,$$

and

$$\left|\partial_r^\alpha(\chi_{\varepsilon,\kappa_0}(r_0+r)\varphi_0(\frac{r}{\delta})Z_+(|x|(r_0+r))(r_0+r)^{d-1})\right| \lesssim 1,$$

which implies that

$$|I_+| \lesssim t^{-N}. \tag{2.35}$$

Following a similar procedure as (2.34), and using the estimates deduced from Lemma A.1,

$$\left|\partial_r^\alpha\left(\frac{1}{\partial_r\Phi_-(r_0+r)}\right)\right|\lesssim 1,\quad \left|\partial_r^\alpha(\chi_{\varepsilon,\kappa_0}(r_0+r)\varphi_0(\frac{r}{\delta})Z_-(|x|(r_0+r))(r_0+r)^{d-1})\right|\lesssim 1,$$

we obtain

$$|I_-|\lesssim t^{-N}, \tag{2.36}$$

for $\frac{|x|}{t}<\frac{1}{100}$. If $\frac{|x|}{t}\geq\frac{1}{100}$, then by Lemma A.1 and choosing $\delta>0$ sufficiently small, we have $|\partial_r^3\Phi_-(r,t,x)|>1/2|b'''(r_0)|$. Applying the Van der Corput Lemma (Lemma B.1), we obtain

$$\begin{aligned}|I_-|&\lesssim t^{-\frac13}\sup_{|r|\leq 2\delta}\partial_r\left(\chi_{\varepsilon,\kappa_0}(r_0+r)\varphi_0(\frac{r}{\delta})Z_-(|x|(r_0+r))(r_0+r)^{d-1}\right)\\&\lesssim t^{-\frac13}(1+|x|)^{-\frac{d-1}{2}}\\&\lesssim t^{-\frac{d}{2}+\frac16}.\end{aligned}\tag{2.37}$$

Therefore, the estimate (2.30) is proved by combining (2.33), (2.35), (2.36), and (2.37)

Then (2.31) and (2.32) can be proved by following a similar argument as the proof of Lemma 2.2. □

**Lemma 2.4.** *For $j=1,2,3$ and for every $\varepsilon\in(0,1]$, we have the estimate*

$$|\chi_{\varepsilon,\kappa_0}(\xi)G_j(t,\xi)|\lesssim e^{-\frac{\varepsilon|\xi|^2t}{2}},$$

*where the cutoff functions $\chi_{\varepsilon,\kappa_0}$ is defined in* (1.16).

*Proof.* Recall the definition in (2.10), we compute that

$$G_1(t,\xi)=e^{-\varepsilon|\xi|^2t}\left(\cos(b(\xi)t)+\varepsilon\frac{\sin(b(\xi)t)}{b(\xi)}|\xi|^2\right), \tag{2.38}$$

$$G_2(t,\xi)=e^{-\varepsilon|\xi|^2t}\frac{\sin(b(\xi)t)}{b(\xi)}a(\xi), \tag{2.39}$$

$$G_3(t,\xi)=e^{-\varepsilon|\xi|^2t}\left(\cos(b(\xi)t)-\varepsilon\frac{\sin(b(\xi)t)}{b(\xi)}|\xi|^2\right). \tag{2.40}$$

Then we have

$$\begin{aligned}|\chi_{\varepsilon,\kappa_0}(\xi)G_j(t,\xi)|&\lesssim e^{-\varepsilon|\xi|^2t}\left(1+\varepsilon|\xi|^2t\right)\lesssim e^{-\frac{\varepsilon|\xi|^2t}{2}},\quad j=1,3,\\|\chi_{\varepsilon,\kappa_0}(\xi)G_2(t,\xi)|&\lesssim e^{-\varepsilon|\xi|^2t}\sqrt{\frac{|\xi|^{2s}+1}{(1-4\kappa_0\varepsilon)|\xi|^{2s}+1}}\lesssim e^{-\varepsilon|\xi|^2t}.\end{aligned}$$

This concludes the proof. □

We are now in a position to prove Proposition 2.1.

*Proof of Proposition 2.1.* According to the definition of the Besov norm (1.19), (2.12) is obtained by using Lemma 2.2 for $s\in(0,1/2]$ and using Lemma 2.3 for $s\in(1/2,1)$. Recall the definition in (2.10),

$$\begin{aligned}G_1(t,D)&=\frac12 e^{\varepsilon t\Delta}(e^{\mathrm{i}tb(D)}+e^{-\mathrm{i}tb(D)})+\mathrm{i}\frac{\varepsilon\Delta}{b(D)}(e^{\mathrm{i}tb(D)}-e^{-\mathrm{i}tb(D)}),\\G_2(t,D)&=e^{\varepsilon t\Delta}\frac{a(D)}{2b(D)}(e^{\mathrm{i}tb(D)}-e^{-\mathrm{i}tb(D)}),\\G_3(t,D)&=\frac12 e^{\varepsilon t\Delta}(e^{-\mathrm{i}tb(D)}+e^{\mathrm{i}tb(D)})+\mathrm{i}\frac{\varepsilon\Delta}{b(D)}(e^{-\mathrm{i}tb(D)}-e^{\mathrm{i}tb(D)}),\end{aligned}$$

by Lemma A.1, for $r \in \mathrm{supp}\chi_{\varepsilon,\kappa_0}$,

$$\left|\partial_r^n\left(\frac{\varepsilon r^2}{b(r)}\right)\right|, \quad \left|\partial_r^n\left(\frac{a(r)}{b(r)}\right)\right| \lesssim 1.$$

Then (2.13) is obtained by the continuous property of the heat operator $e^{\varepsilon t\Delta}$ on $L^p$ $(1 \le p \le +\infty)$ and following a similar procedure as the proof of Lemma 2.2 and Lemma 2.3. □

2.1.2. *Estimates in $\Omega^H$.* The main result in this part is stated as follows.

**Lemma 2.5.** *There exists $c_0 > 0$ such that, for $j = 1, 2, 3$ and for every $\varepsilon \in (0, 1]$, we have the estimate*

$$|(1 - \chi_{\varepsilon,\kappa_0}(\xi))G_j(t, \xi)| \lesssim e^{-c_0 t},$$

*where the cutoff functions $\chi_{\varepsilon,\kappa_0}$ is defined in* (1.16)*.*

*Proof.* First, we consider the case $\xi : a^2(\xi) \ge \varepsilon^2|\xi|^4$. Recall $G_1$, $G_2$, and $G_3$ in (2.38), (2.39), and (2.40), respectively. Therefore, for $k = 1, 3$, we have

$$|(1 - \chi_{\varepsilon,\kappa_0}(\xi))G_k(t, \xi)| \le e^{-\varepsilon|\xi|^2 t}(1 + \varepsilon|\xi|^2 t)|(1 - \chi_{\varepsilon.\kappa_0}(\xi))| \lesssim e^{-\frac{1}{2}\varepsilon|\xi|^2 t}|(1 - \chi_{\varepsilon,\kappa_0}(\xi))| \lesssim e^{-\frac{1}{2}\kappa_0 t}.$$

For $G_2$, if $b(\xi) \ge a(\xi)/2$, then we have

$$|(1 - \chi_{\varepsilon,\kappa_0}(\xi))G_2(t, \xi)| \le e^{-\varepsilon|\xi|^2 t}|(1 - \chi_{\varepsilon,\kappa_0}(\xi))| \lesssim e^{-\frac{1}{2}\kappa_0 t}.$$

If $b(\xi) \le a(\xi)/2$, which implies $a(\xi) \le \frac{2}{\sqrt{3}}\varepsilon|\xi|^2$, then we have

$$|(1 - \chi_{\varepsilon,\kappa_0}(\xi))G_2(t, \xi)| \le 2\varepsilon|\xi|^2 t e^{-\varepsilon|\xi|^2 t}|1 - \chi_{\varepsilon,\kappa_0}(\xi)| \lesssim e^{-\frac{1}{2}\kappa_0 t}.$$

It remains to consider the case $a^2(\xi) < \varepsilon^2|\xi|^4$. Let us introduce $\widetilde{b}(\xi) := \sqrt{\varepsilon^2|\xi|^4 - a^2(\xi)}$, then recall the definition in (2.10), we compute

$$G_1(t, \xi) = e^{\lambda_-(\xi)t}\left(1 - \frac{1 - e^{-2\widetilde{b}(\xi)t}}{2\widetilde{b}(\xi)}\lambda_-(\xi)\right),$$

$$G_2(t, \xi) = e^{\lambda_-(\xi)t}\frac{1 - e^{-2\widetilde{b}(\xi)t}}{2\widetilde{b}(\xi)}a(\xi),$$

$$G_3(t, \xi) = e^{\lambda_-(\xi)t}\left(1 + \frac{1 - e^{-2\widetilde{b}(\xi)t}}{2\widetilde{b}(\xi)}\lambda_+(\xi)\right),$$

where $\lambda_\pm(\xi) = -\varepsilon|\xi|^2 \mp \widetilde{b}(\xi)$. Moreover, we have

$$\lambda_+(\xi) > -2\varepsilon|\xi|^2, \quad \lambda_-(\xi) = \frac{|\xi|^2 + |\xi|^{2-2s}}{\lambda_+(\xi)} \le -\frac{1}{2\varepsilon}.$$

Therefore, for $G_1$, we have

$$|(1 - \chi_{\varepsilon,\kappa_0})(\xi)G_1(t, \xi)| \le e^{\lambda_-(\xi)t}(1 - \lambda_- t) \lesssim e^{\frac{1}{2}\lambda_-(\xi)t} \lesssim e^{-\frac{1}{4}t}.$$

For $G_2$, if $\widetilde{b}(\xi) \le \varepsilon|\xi|^2/2$, then we have

$$|(1 - \chi_{\varepsilon,\kappa_0}(\xi))G_2(t, \xi)| \lesssim \varepsilon|\xi|^2 t e^{-\frac{1}{2}\varepsilon|\xi|^2 t} \lesssim e^{-\frac{1}{4}\varepsilon|\xi|^2 t} \lesssim e^{-\frac{1}{4}\kappa_0 t}.$$

If $\widetilde{b}(\xi) > \varepsilon|\xi|^2/2$, which implies $a^2(\xi) \le 3\widetilde{b}^2(\xi)$, then we have

$$|(1 - \chi_{\varepsilon,\kappa_0}(\xi))G_2(t, \xi)| \lesssim e^{\lambda_-(\xi)t} \lesssim e^{-\frac{1}{4}t}.$$

For $G_3$, if $\widetilde{b}(\xi) > \varepsilon|\xi|^2/2$, then

$$|(1-\chi_{\varepsilon,\kappa_0}(\xi))G_3(t,\xi)| \le e^{\lambda_- t}(1-\frac{\lambda_+(\xi)}{\widetilde{b}(\xi)}) \lesssim e^{\lambda_-(\xi)t} \lesssim e^{-\frac{1}{2}t}.$$

If $0 \le \widetilde{b}(\xi) \le \varepsilon|\xi|^2/2$, which implies $\lambda_-(\xi) \le -\varepsilon|\xi|^2/2$, then

$$|(1-\chi_{\varepsilon,\kappa_0}(\xi))G_3(t,\xi)| \le e^{\lambda_-(\xi)t}(1-\lambda_+(\xi)t) \lesssim e^{-\frac{\varepsilon}{2}|\xi|^2 t}(1+2\varepsilon|\xi|^2 t) \lesssim e^{-\frac{1}{4}\kappa_0 t}.$$

Combining these estimates, we conclude the proof of Lemma 2.5. □

As a matter of fact, we obtain the uniform boundedness under the $\dot{H}^\sigma$ norm in the whole frequency domain.

**Corollary 2.6.** *For $s \in (0,1)$ and for every $\sigma \in \mathbb{R}$, we have, uniformly in $\varepsilon \in (0,1]$, the estimate*

$$\|G_j(t,D)f\|_{\dot{H}^\sigma} \lesssim \|f\|_{\dot{H}^\sigma}, \quad j=1,2,3. \tag{2.41}$$

*Proof.* By Lemma 2.4 and Lemma 2.5, we have

$$|G_j(t,\xi)| \lesssim 1, \quad j=1,2,3,$$

then by the Plancherel theorem, we deduce that

$$\|G_j(t,D)f\|_{\dot{H}^\sigma} \lesssim \|G_j(t,\xi)\|_{L^\infty}\|f\|_{\dot{H}^\sigma},$$

which proves (2.41). □

We can further derive the linear time-decay rate in $L^p$ space for $2 \le p < +\infty$.

**Corollary 2.7.** *For $2 \le p < +\infty$, we have, uniformly in $\varepsilon \in (0,1]$, the estimate*

$$\|e^{\mathrm{i}tb(D)}\chi_{\varepsilon,\kappa_0}(D)f\|_{L^p} \lesssim (1+t)^{-\beta_{d,s,p}}\|f\|_{W^{(1-\frac{2}{p})d,\frac{p}{p-1}}}. \tag{2.42}$$

*Moreover,*

$$\|G_j(t,D)\chi_{\varepsilon,\kappa_0}(D)f\|_{L^p} \lesssim (1+t)^{-\beta_{d,s,p}}\|f\|_{W^{(1-\frac{2}{p})d,\frac{p}{p-1}}}. \tag{2.43}$$

*Here, the cutoff functions $\chi_{\varepsilon,\kappa_0}$ is defined in* (1.16), *and $\beta_{d,s,p}$ is defined as*

$$\beta_{d,s,p} := \begin{cases} \frac{d}{2}\left(1-\frac{2}{p}\right), & \text{when } 0 < s \le 1/2, \\ \left(\frac{d}{2}-\frac{1}{6}\right)\left(1-\frac{2}{p}\right), & \text{when } 1/2 < s < 1. \end{cases}$$

*Proof.* The estimates (2.42) and (2.43) can be obtained by interpolating in a classical way between (2.41) in Corollary 2.6 and (2.12) in Proposition 2.1, with the embeddings $B^k_{p,2}(\mathbb{R}^d) \hookrightarrow W^{k,p}(\mathbb{R}^d)$, $W^{k,p/(p-1)}(\mathbb{R}^d) \hookrightarrow B^k_{p/(p-1),2}(\mathbb{R}^d)$ with $2 \le p < +\infty$ and $k \ge 0$. □

2.2. **Global Existence Results.** To prove Theorem 1.3, we introduce the following norms:

$$\|f\|_Y := \|f\|_{W^{N_0+3(1-\frac{2}{p}),\frac{p}{p-1}}} + \|f\|_{\dot{H}^{-(1-s)}\cap H^N} + \|f\|_{L^1},$$

$$\|f\|_{X_t} := \sup_{0\le\tau\le t}\left((1+\tau)^{\beta_s}\|f(\tau)\|_{W^{N_0,p}} + (1+\tau)^{\beta_s}\|\chi^H(D)f(\tau)\|_{H^{N-1}} + \|f(\tau)\|_{\dot{H}^{-(1-s)}\cap H^N}\right),$$

where $\chi^H(D) = (1-\chi_{\varepsilon,\kappa_0})(D)$, $N_0 \ge 6, N \ge N_0+14$, $\beta_s$ is defined by (1.9) and $p$ is chosen in (1.8).

Let us denote $U := (\varrho, u, |\nabla|^s\phi)$. For the subsequent analysis, we impose the following a priori assumption: There exists a sufficiently small constant $\delta_1 > 0$ such that

$$\|U\|_{X_T} \le \delta_1. \tag{2.44}$$

We shall use the smallness condition (2.44) as a bootstrap hypothesis. More precisely, for a fixed $T > 0$, we assume that a smooth solution $(\varrho, u, \phi)$ on $[0, T]$ satisfies (2.44); all estimates derived below are uniform with respect to $T$.

2.2.1. *$H^N$ Estimates.* In this part, we will derive the estimate in $H^N$ norm by using the nonlinear energy method. We introduce the new quantity $\mathfrak{c}$

$$\mathfrak{c} := \frac{2}{\gamma - 1}(\rho^{\frac{\gamma-1}{2}} - 1), \tag{2.45}$$

which together with $u$ transforms (1.7) into a symmetric hyperbolic–parabolic system. More precisely, setting $\mathfrak{U} = (\mathfrak{c}, u)^{\mathsf{T}}$, the system takes the form

$$\partial_t \mathfrak{U} + \sum_{i=1}^{3} \mathbf{M}_i(\mathfrak{U})\partial_{x_i}\mathfrak{U} - 2\varepsilon \mathbf{N}\Delta\mathfrak{U} + (0, \nabla\phi)^{\mathsf{T}} = 0, \tag{2.46}$$

where

$$\mathbf{M}_i(\mathfrak{U}) = \begin{pmatrix} u_i & \left(\frac{\gamma-1}{2}\mathfrak{c} + 1\right) e_i \\ \left(\frac{\gamma-1}{2}\mathfrak{c} + 1\right) e_i^{\mathsf{T}} & u_i \mathbb{I}_3 \end{pmatrix}, \quad \mathbf{N} = \begin{pmatrix} 0 & 0 \\ 0 & \mathbb{I}_3 \end{pmatrix}. \tag{2.47}$$

We introduce the energy and dissipation functionals order by order. For the base case $k = 0$, the zero-order energy functional is defined as

$$\mathbb{E}^0(t) := \|(\mathfrak{c}, u, |\nabla|^s\phi)\|_{L^2}^2.$$

For the higher-order derivatives ($k \geq 1$), the energy functional is defined according to the fractional exponent $s$

$$\mathbb{E}^k(t) := \begin{cases} \left\|(\nabla^k\mathfrak{c}, \nabla^k u, \nabla^k|\nabla|^s\phi)\right\|_{L^2}^2 & \text{if } 1/2 \leq s < 1, \\ \left\|\left(\nabla^k\mathfrak{c}, \nabla^k u, \frac{1}{\sqrt{\rho}}\nabla^k|\nabla|^s\phi\right)\right\|_{L^2}^2 & \text{if } 0 < s < 1/2. \end{cases}$$

Based on these, we introduce the cumulative energy functional

$$\mathbb{E}_N^k(t) := \sum_{m=k}^{N} \mathbb{E}^m(t)$$

along with its associated cumulative dissipation rate $\mathbb{D}_N^k(t)$

$$\mathbb{D}_N^k(t) := \|\nabla^k(\nabla\mathfrak{c}, \nabla|\nabla|^s\phi)\|_{H^{N-1-k}}^2 + \|\nabla^k\nabla u\|_{H^{N-k}}^2. \tag{2.48}$$

For simplicity, we drop the superscript 0 when considering the full range from order 0 to $N$, writing the total energy and dissipation rate functionals as $\mathbb{E}_N(t) := \mathbb{E}_N^0(t)$ and $\mathbb{D}_N(t) := \mathbb{D}_N^0(t)$, respectively.

To close the energy estimates in the subsequent analysis, we deduce from (2.45) and (2.44) the following equivalence relations

$$\begin{cases} \|\mathfrak{c}\|_{W^{k,\infty}} \simeq \|\varrho\|_{W^{k,\infty}}, & \text{for any } k \in \mathbb{N} \text{ and } k \leq N-2, \\ \|\nabla^k\mathfrak{c}\|_{L^2} \simeq \|\nabla^k\varrho\|_{L^2}, & \text{for any } k \in \mathbb{N} \text{ and } k \leq N, \\ \||\nabla|^{-\vartheta}\mathfrak{c}\|_{L^2} \simeq \||\nabla|^{-\vartheta}\varrho\|_{L^2}, & \text{for any } 0 < \vartheta < 3/2. \end{cases} \tag{2.49}$$

**Proposition 2.8.** *If* (2.44) *holds, then we have*

$$\frac{1}{2}\frac{d}{dt}\mathbb{E}_N(t) + 2\varepsilon\|\nabla u\|_{H^N}^2 \lesssim (\|(\varrho, \mathfrak{c}, u)\|_{W^{1,\infty}} + \mathbf{1}_{0<s<\frac{1}{2}}\|\nabla\varrho\|_{L^{\frac{3}{s}}})\,\mathbb{E}_N(t), \tag{2.50}$$

$$\frac{1}{2}\frac{d}{dt}\mathbb{E}_N^1(t) + 2\varepsilon\|\nabla^2 u\|_{H^{N-1}}^2 \lesssim (\|(\varrho, \mathfrak{c}, u)\|_{W^{1,\infty}} + \mathbf{1}_{0<s<\frac{1}{2}}\|\nabla\varrho\|_{L^{\frac{3}{s}}})\,\mathbb{E}_N^1(t). \tag{2.51}$$

*Proof.* Applying $\nabla^k$ to (2.46), and taking the $L^2$ inner product with $\nabla^k\mathfrak{U}$, we obtain

$$\frac{1}{2}\frac{d}{dt}\|\nabla^k\mathfrak{U}\|_{L^2}^2 + 2\varepsilon\|\nabla^k\nabla u\|_{L^2}^2 = \sum_{i=1}^{3} I_i, \tag{2.52}$$

where

$$\begin{aligned}
I_1 &:= \frac{1}{2}\int_{\mathbb{R}^3} \nabla^k\mathfrak{U}^{\mathsf{T}}\,\partial_{x_i}M_i(\mathfrak{U})\,\nabla^k\mathfrak{U}\,\mathrm{d}x,\\
I_2 &:= -\mathbf{1}_{k\geq 1}\int_{\mathbb{R}^3} \nabla^k\mathfrak{U}^{\mathsf{T}}[\nabla^k, M_i(\mathfrak{U})]\partial_{x_i}\mathfrak{U}\,\mathrm{d}x,\\
I_3 &:= \int_{\mathbb{R}^3} \nabla^k\phi\,\nabla^k\mathrm{div}\,u\,\mathrm{d}x.
\end{aligned}$$

By Hölder's inequality and Lemma B.4, we readily obtain

$$|I_1| \lesssim \|\nabla\mathfrak{U}\|_{L^\infty}\|\nabla^k\mathfrak{U}\|_{L^2}^2, \tag{2.53}$$

$$|I_2| \lesssim \mathbf{1}_{k\geq 1}\left(\|\mathfrak{U}\|_{\dot{H}^k}\|\nabla\mathfrak{U}\|_{L^\infty} + \|\nabla\mathfrak{U}\|_{L^\infty}\|\nabla\mathfrak{U}\|_{\dot{H}^{k-1}}\right)\|\nabla^k\mathfrak{U}\|_{L^2} \lesssim \|\nabla\mathfrak{U}\|_{L^\infty}\|\nabla^k\mathfrak{U}\|_{L^2}^2. \tag{2.54}$$

It remains to estimate $I_3$, which is treated separately as follows.
**Case 1: $1/2 \leq s < 1$.** Substituting $(1.7)_1$ and $(1.7)_3$ into $I_3$ and integrating by parts, we obtain

$$I_3 = -\frac{1}{2}\frac{d}{dt}\|\nabla^k|\nabla|^s\phi\|_{L^2}^2 + \int_{\mathbb{R}^3}\nabla^k\nabla\phi\cdot\nabla^k(\varrho u)\,\mathrm{d}x. \tag{2.55}$$

Using the Sobolev interpolation and Lemma B.4, the last term of the above equality can be estimated as

$$\begin{aligned}
\left|\int_{\mathbb{R}^3}\nabla^k\nabla\phi\cdot\nabla^k(\varrho u)\,\mathrm{d}x\right| &\lesssim \|\nabla^k\nabla\phi\|_{L^2}\|\nabla^k(\varrho u)\|_{L^2}\\
&\lesssim \|\nabla^k|\nabla|^s\phi\|_{L^2}^{\frac{2s-1}{s}}\|\nabla^k(-\Delta)^s\phi\|_{L^2}^{\frac{1-s}{s}}\|(\varrho,u)\|_{L^\infty}\|\nabla^k(\varrho,u)\|_{L^2}\\
&\lesssim \|(\varrho,u)\|_{L^\infty}\|\nabla^k(\varrho,u,|\nabla|^s\phi)\|_{L^2}^2,
\end{aligned} \tag{2.56}$$

which is valid for $1/2 \leq s < 1$. Plugging (2.53), (2.54), (2.55), and (2.56) into (2.52), using (2.49), and summing over $k \leq N$ and $1 \leq k \leq N$ respectively, yields (2.50) and (2.51) for $1/2 \leq s < 1$.
**Case 2: $0 < s < 1/2$.** In this case, the estimate (2.56) is no longer valid. We use a different argument. For $k = 0$, we have from Lemma B.4 that

$$\left|\int_{\mathbb{R}^3}\nabla\phi\cdot\varrho u\,\mathrm{d}x\right| \lesssim \|(-\Delta)^s\phi\|_{L^2}\||\nabla|^{1-2s}(\varrho u)\|_{L^2} \lesssim \|(\varrho,u)\|_{L^\infty}\|(\varrho,u)\|_{H^1}^2,$$

which yields that

$$\frac{1}{2}\frac{d}{dt}\mathbb{E}^0(t) + 2\varepsilon\|\nabla u\|_{L^2}^2 \lesssim \|(\varrho,\mathfrak{c},u)\|_{W^{1,\infty}}\|(\varrho,\mathfrak{c},u)\|_{H^1}^2. \tag{2.57}$$

For the case $k \geq 1$, we first note

$$\frac{3}{4} \leq \|\rho\|_{L^\infty} \leq \frac{5}{4}, \tag{2.58}$$

due to (2.44) and the inequality $1 - \|\varrho\|_{L^\infty} \leq \|\rho\|_{L^\infty} \leq 1 + \|\varrho\|_{L^\infty}$, which allows us to apply commutator estimates. We thus use integration by parts to get

$$I_3 = \int_{\mathbb{R}^3}\nabla^k|\nabla|^s\phi\nabla^k|\nabla|^{-s}\mathrm{div}\,u\,\mathrm{d}x = I_4 + I_5, \tag{2.59}$$

where

$$I_4 := -\int_{\mathbb{R}^3} \frac{1}{\rho}\left[\nabla^k|\nabla|^{-s},\rho\right]\operatorname{div} u \nabla^k|\nabla|^s\phi\,\mathrm{d}x,$$
$$I_5 := \int_{\mathbb{R}^3} \frac{1}{\rho}\nabla^k|\nabla|^{-s}(\rho \operatorname{div} u)\nabla^k|\nabla|^s\phi\,\mathrm{d}x.$$

Since $k-s>0$, Lemma B.4 is applicable, which allows us to estimate $I_4$ as follows

$$\begin{aligned}|I_4| &\lesssim \left(\|\rho\|_{\dot{H}^{k-s}}\|\nabla u\|_{L^\infty} + \|\nabla\rho\|_{L^{\frac{3}{s}}}\|\nabla u\|_{\dot{W}^{k-s-1,\frac{6}{3-2s}}}\right)\|\nabla^k|\nabla|^s\phi\|_{L^2}\\ &\lesssim (\|\nabla u\|_{L^\infty} + \|\nabla\varrho\|_{L^{\frac{3}{s}}})\|\nabla^k(u,|\nabla|^s\phi)\|_{L^2}^2. \end{aligned} \tag{2.60}$$

Here, we have used the Sobolev embedding $\dot{H}^s(\mathbb{R}^3)\hookrightarrow L^{\frac{6}{3-2s}}(\mathbb{R}^3)$. For the term $I_5$, using $(1.7)_1$ and $(1.7)_3$, we rewrite

$$\begin{aligned} I_5 =& -\int_{\mathbb{R}^3}\frac{1}{\rho}\,\nabla^k|\nabla|^{-s}(\partial_t\rho + u\cdot\nabla\rho)\,\nabla^k|\nabla|^s\phi\;\mathrm{d}x\\ =& -\frac{1}{2}\frac{d}{dt}\int_{\mathbb{R}^3}\frac{1}{\rho}|\nabla^k|\nabla|^s\phi|^2\,\mathrm{d}x - \frac{1}{2}\int_{\mathbb{R}^3}\frac{1}{\rho^2}\partial_t\rho|\nabla^k|\nabla|^s\phi|^2\,\mathrm{d}x\\ &-\int_{\mathbb{R}^3}\frac{1}{\rho}\nabla^k|\nabla|^{-s}(u\cdot\nabla\rho)\nabla^k|\nabla|^s\phi\,\mathrm{d}x\\ =&\sum_{i=1}^{5} I_{5i}, \end{aligned}$$

where

$$I_{51} := -\frac{1}{2}\frac{d}{dt}\int_{\mathbb{R}^3}\frac{1}{\rho}|\nabla^k|\nabla|^s\phi|^2\,\mathrm{d}x, \qquad I_{52} := \frac{1}{2}\int_{\mathbb{R}^3}\frac{1}{\rho}\operatorname{div} u|\nabla^k|\nabla|^s\phi|^2\,\mathrm{d}x,$$
$$I_{53} := \frac{1}{2}\int_{\mathbb{R}^3}\frac{1}{\rho^2}u\cdot\nabla\rho|\nabla^k|\nabla|^s\phi|^2\,\mathrm{d}x, \qquad I_{54} := -\int_{\mathbb{R}^3}\frac{1}{\rho}\left[\nabla^k|\nabla|^{-s},u\right]\cdot\nabla\rho\,\nabla^k|\nabla|^s\phi\,\mathrm{d}x,$$
$$I_{55} := -\frac{1}{2}\int_{\mathbb{R}^3}\frac{1}{\rho}u\cdot\nabla|\nabla^k|\nabla|^s\phi|^2\,\mathrm{d}x.$$

By integration by parts and (2.58), we have

$$|I_{52}+I_{53}+I_{55}| = 2|I_{52}| \lesssim \|\nabla u\|_{L^\infty}\|\nabla^k|\nabla|^s\phi\|_{L^2}^2. \tag{2.61}$$

As for the remaining term $I_4$, following an argument similar to the one for $I_4$, we get

$$|I_{54}| \lesssim (\|\nabla u\|_{L^\infty} + \|\nabla\varrho\|_{L^{\frac{3}{s}}})\|\nabla^k(u,|\nabla|^s\phi)\|_{L^2}^2. \tag{2.62}$$

Substituting (2.53), (2.54), and (2.59)-(2.62) into (2.52) and taking summation over $1\le k\le N$, we combine the result with (2.57) and (2.49) to obtain (2.50) and (2.51) for $0<s<1/2$. □

**Corollary 2.9.** *If* (2.44) *holds, then we have*

$$\|U(t)\|_{H^N}^2 + 2\varepsilon\int_0^t\|\nabla u(\tau)\|_{H^N}^2\,\mathrm{d}\tau \lesssim \|U_0\|_{H^N}^2 + \|U\|_{X_t}^3. \tag{2.63}$$

*Proof.* For $0<s<1/2$, we first derive the following decay estimates for the $W^{1,\frac{3}{s}}$-norm of $\varrho$

$$\|\varrho\|_{W^{1,\frac{3}{s}}} \lesssim \begin{cases} \|\varrho\|_{H^1}^{1-\theta}\|\varrho\|_{W^{1,p}}^{\theta} \lesssim (1+t)^{-\frac{3}{2}(1-\frac{2s}{3})}\|\varrho\|_{X_t}, & \text{if } \frac{3}{s}\le p,\\ \|\varrho\|_{W^{2,p}} \lesssim (1+t)^{-\frac{3}{2}(1-\frac{2}{p})}\|\varrho\|_{X_t}, & \text{if } p<\frac{3}{s}, \end{cases} \tag{2.64}$$

where $\theta = \frac{p(3-2s)}{3(p-2)}$ is the Sobolev interpolation parameter. Then by substituting (2.64) into the Lyapunov-type inequality (2.50) and using the Sobolev embedding $W^{1+\frac{3}{p}+\delta,p}(\mathbb{R}^3) \hookrightarrow W^{1,\infty}(\mathbb{R}^3)$ for any $\delta > 0$ along with (2.58), a standard application of Gronwall's inequality directly yields the desired result (2.63). □

To obtain the time-decay estimates in $H^N$ norm, we recover the estimates for dissipation in perturbations of density and Riesz potential.

**Proposition 2.10.** *If* (2.44) *holds, then for* $0 \le k \le N-1$ *we have*

$$\frac{d}{dt}\int_{\mathbb{R}^3} \nabla^k u \cdot \nabla^k \nabla \varrho \, \mathrm{d}x + \|\nabla^k(\nabla\varrho, \nabla|\nabla|^s\phi)\|_{L^2}^2 \lesssim \|\nabla^k \nabla u\|_{H^1}^2. \tag{2.65}$$

*Proof.* We rewrite the system (1.7) into the system for perturbation $(\varrho, u)$,

$$\begin{cases} \partial_t \varrho + \operatorname{div} u = -\operatorname{div}(\varrho u), \\ \partial_t u + u\cdot\nabla u + \nabla\varrho + \nabla\phi - \varepsilon\mathcal{L}u = -\nabla\left[\frac{1}{\gamma-1}(1+\varrho)^{\gamma-1} - \frac{1}{\gamma-1} - \varrho\right], \\ (-\Delta)^s\phi = \varrho. \end{cases} \tag{2.66}$$

Applying $\nabla^k$ to $(2.66)_2$, taking the $L^2$ inner product with $\nabla^k\nabla\varrho$, and utilizing $(2.66)_1$ and $(2.66)_3$, we obtain

$$\frac{d}{dt}\int_{\mathbb{R}^3} \nabla^k u \cdot \nabla^k \nabla \varrho \, \mathrm{d}x + \|\nabla^k(\nabla\varrho, \nabla|\nabla|^s\phi)\|_{L^2}^2 = \sum_{i=6}^{9} I_i,$$

where

$$I_6 := -\int_{\mathbb{R}^3} \nabla^k u \cdot \nabla^k\nabla(\operatorname{div} u + \operatorname{div}(\varrho u))\,\mathrm{d}x, \qquad I_7 := -\int_{\mathbb{R}^3} \nabla^k(u\cdot\nabla u)\cdot\nabla^k\nabla\varrho\,\mathrm{d}x,$$

$$I_8 := -\int_{\mathbb{R}^3} \nabla^k\nabla\left[\frac{1}{\gamma-1}(1+\varrho)^{\gamma-1} - \frac{1}{\gamma-1} - \varrho\right]\cdot\nabla^k\nabla\varrho\,\mathrm{d}x, \quad I_9 := \varepsilon\int_{\mathbb{R}^3} \nabla^k\mathcal{L}u\cdot\nabla^k\nabla\varrho\,\mathrm{d}x.$$

Next, we estimate the terms $I_6$–$I_9$. For the term $I_6$, integrating by parts, applying Hölder's inequality along with Sobolev embeddings, and utilizing Lemma B.4 yields

$$\begin{aligned} |I_6| &\lesssim \|\nabla^k\nabla u\|_{L^2}^2 + \|\nabla^k \operatorname{div}(\varrho u)\|_{L^2}\|\nabla^k\nabla u\|_{L^2} \\ &\lesssim \|\nabla^k\nabla u\|_{L^2}^2 + \left(\|\nabla^k\nabla\varrho\|_{L^2}\|u\|_{L^\infty} + \|\nabla^k\nabla u\|_{L^2}\|\varrho\|_{L^\infty}\right)\|\nabla^k\nabla u\|_{L^2} \\ &\lesssim \|\nabla^k\nabla u\|_{L^2}^2 + \|(\varrho,u)\|_{H^2}\|\nabla^k(\nabla\varrho,\nabla u)\|_{L^2}^2. \end{aligned}$$

By a similar argument for $I_7$, one can deduce

$$|I_7| \lesssim \|\nabla^k(u\cdot\nabla u)\|_{L^2}\|\nabla^k\nabla\varrho\|_{L^2} \lesssim \|u\|_{H^2}\|\nabla^k(\nabla\varrho,\nabla u)\|_{L^2}^2.$$

To handle $I_8$, we first employ a Taylor expansion to express the nonlinear term in the form

$$\frac{1}{\gamma-1}(1+\varrho)^{\gamma-1} - \frac{1}{\gamma-1} - \varrho = (\gamma-2)\varrho^2\int_0^1 (1+\theta\varrho)^{\gamma-3}(1-\theta)\,\mathrm{d}\theta := \varrho^2\mathcal{N}(\varrho).$$

This representation, together with Hölder's inequality, Lemma B.4, Lemma B.5 and Sobolev embeddings, implies that

$$\begin{aligned} |I_8| &\lesssim \|\nabla^k\nabla(\varrho^2\mathcal{N}(\varrho))\|_{L^2}\|\nabla^k\nabla\varrho\|_{L^2} \\ &\lesssim \|\nabla^k\nabla\varrho\|_{L^2}\left(\|\nabla^k\nabla(\varrho^2)\|_{L^2}\|\mathcal{N}(\varrho)\|_{L^\infty} + \|\nabla^k\nabla\mathcal{N}(\varrho)\|_{L^2}\|\varrho^2\|_{L^\infty}\right) \\ &\lesssim \|\nabla^k\nabla\varrho\|_{L^2}\left(\|\nabla^k\nabla\varrho\|_{L^2}\|\varrho\|_{L^\infty} + \|\nabla^k\nabla\varrho\|_{L^2}\|\varrho\|_{L^\infty}^2\right) \\ &\lesssim \|\varrho\|_{H^2}(1+\|\varrho\|_{H^2})\|\nabla^k\nabla\varrho\|_{L^2}^2. \end{aligned}$$

For $I_9$, an application of the Cauchy-Schwarz inequality gives

$$|I_9| \lesssim \eta\|\nabla^k\nabla\varrho\|_{L^2}^2 + \varepsilon^2\|\nabla^{k+1}\nabla u\|_{L^2}^2.$$

Consequently, the desired estimate (2.65) follows from the above estimates and the *a priori* assumption (2.44), provided that $\eta$ is chosen sufficiently small. □

We next combine the basic energy inequalities in Proposition 2.8 with the compensating estimate in Proposition 2.10. Then the dissipation of $(\nabla\varrho, \nabla|\nabla|^s\phi)$ is recovered. This leads to modified energies equivalent to $\mathbb{E}_N$ and $\mathbb{E}_N^1$, and hence to the following Lyapunov-type estimates.

**Corollary 2.11.** *If* (2.44) *holds, then we have*

$$\frac{d}{dt}\bar{\mathbb{E}}_N(t) + C\varepsilon\kappa_1\mathbb{D}_N(t) \lesssim (\|(\varrho, \mathfrak{c}, u)\|_{W^{1,\infty}} + \mathbf{1}_{0<s<\frac{1}{2}}\|\nabla\varrho\|_{L^{\frac{3}{s}}})\,\mathbb{E}_N(t), \tag{2.67}$$

$$\frac{d}{dt}\bar{\mathbb{E}}_N^1(t) + C\varepsilon\kappa_2\mathbb{D}_N^1(t) \lesssim (\|(\varrho, \mathfrak{c}, u)\|_{W^{1,\infty}} + \mathbf{1}_{0<s<\frac{1}{2}}\|\nabla\varrho\|_{L^{\frac{3}{s}}})\,\mathbb{E}_N^1(t), \tag{2.68}$$

*where $\mathbb{D}_N$ and $\mathbb{D}_N^1$ are defined in* (2.48)*, the modified energy functionals are defined as*

$$\bar{\mathbb{E}}_N(t) := \mathbb{E}_N(t) + C\varepsilon\kappa_1\sum_{k=0}^{N-1}\int_{\mathbb{R}^3}\nabla^k u\cdot\nabla^{k+1}\varrho\,\mathrm{d}x,\quad \bar{\mathbb{E}}_N^1(t) := \mathbb{E}_N^1(t) + C\varepsilon\kappa_2\sum_{k=1}^{N-1}\int_{\mathbb{R}^3}\nabla^k u\cdot\nabla^{k+1}\varrho\,\mathrm{d}x,$$

*with constants $\kappa_1, \kappa_2 > 0$ chosen sufficiently small such that $\bar{\mathbb{E}}_N(t) \simeq \mathbb{E}_N(t)$ and $\bar{\mathbb{E}}_N^1(t) \simeq \mathbb{E}_N^1(t)$.*

*Proof.* First, in view of Proposition 2.10, by summing the corresponding estimates (2.65) over $0 \le k \le N-1$, we obtain the following estimate

$$\frac{d}{dt}\sum_{k=0}^{N-1}\int_{\mathbb{R}^3}\nabla^k u\cdot\nabla^{k+1}\varrho\,\mathrm{d}x + \|(\nabla\varrho, \nabla|\nabla|^s\phi)\|_{H^{N-1}}^2 \lesssim \|\nabla u\|_{H^N}^2. \tag{2.69}$$

By multiplying (2.69) by $\varepsilon\kappa_1$ for a sufficiently small constant $\kappa_1 > 0$, and adding it to (2.50), we establish (2.67) in light of the *a priori* assumption (2.44). Similarly, by summing (2.65) over $1 \le k \le N-1$, multiplying the resulting inequality by $\varepsilon\kappa_2$ with $0 < \kappa_2 \ll 1$, and adding it to (2.51), we deduce (2.68). This concludes the proof. □

2.2.2. $W^{k,p}$ *Estimates.* In this part, we prove the result concerning the estimate in the $W^{k,p}$ norm.

**Proposition 2.12.** *If* (2.44) *holds, then we have*

$$(1+t)^{\beta_s}\|U(t)\|_{W^{N_0,p}} + (1+t)^{\beta_s}\|(1-\chi_{\varepsilon,\kappa_0}(D))U(t)\|_{H^{N-1}} \lesssim \|U_0\|_Y + \|U\|_{X_t}^2 + \|U\|_{X_t}^4,$$

*where $\chi_{\varepsilon,\kappa_0}$ is defined in* (1.16).

*Proof.* We recall (2.1) and $V = (h, q)^{\mathsf{T}}$, then $\|U\|_{X_T} \simeq \|V\|_{X_T}$. Indeed, we have

$$\varrho = \frac{|\nabla|^s}{(1+|\nabla|^{2s})^{1/2}}h, \qquad |\nabla|^s\phi = \frac{1}{(1+|\nabla|^{2s})^{1/2}}h,$$

and

$$h = \frac{|\nabla|^s}{(1+|\nabla|^{2s})^{1/2}}\varrho + \frac{1}{(1+|\nabla|^{2s})^{1/2}}|\nabla|^s\phi.$$

The multipliers appearing above satisfy the Hörmander-Mikhlin conditions in Lemma B.6. Since they commute with $\nabla^k$, it follows that

$$\|\nabla^k h\|_{L^p} \simeq \|\nabla^k(\varrho, |\nabla|^s\phi)\|_{L^p}.$$

Moreover, $q = -\mathcal{R}^*u$, and the irrotationality of $u$ gives $u = -\mathcal{R}q$. Hence, by the $L^p$-boundedness of the Riesz transforms with $p \ge 2$,

$$\|\nabla^k q\|_{L^p} \simeq \|\nabla^k u\|_{L^p}.$$

Combining the preceding two estimates yields

$$\|\nabla^k U\|_{L^p} \simeq \|\nabla^k V\|_{L^p}.$$

Consequently, it suffices to establish the corresponding estimate for $V$. For notational convenience, we decompose it into the low-frequency and high-frequency parts:

$$V^L := \chi^L(D)V := \chi_{\varepsilon,\kappa_0}(D)V, \quad V^H := \chi^H(D)V := (1-\chi_{\varepsilon,\kappa_0})(D)V. \tag{2.70}$$

By using (2.9), Lemma 2.5, the Sobolev embeddings $W^{N_0,p}(\mathbb{R}^3) \hookrightarrow L^\infty(\mathbb{R}^3)$, $H^N(\mathbb{R}^3) \hookrightarrow L^\infty(\mathbb{R}^3)$, Lemma B.7, and Lemma B.10, we have

$$\begin{aligned}
\|V^H(t)\|_{H^{N-1}} &\lesssim e^{-c_0 t}\|V(0)\|_{H^{N-1}} + \int_0^t e^{-c_0(t-\tau)}\|B(V)(\tau)\|_{H^{N-1}}\,\mathrm{d}\tau \\
&\lesssim e^{-c_0 t}\|V(0)\|_{Y} + \int_0^t e^{-c_0(t-\tau)}\|V(\tau)\|_{H^N}\|V(\tau)\|_{W^{N_0,p}}(1+\|V(\tau)\|_{H^N})\,\mathrm{d}\tau \\
&\lesssim e^{-c_0 t}\|V(0)\|_{Y} + (1+t)^{-\beta_s}\left(\|V\|_{X_t}^2 + \|V\|_{X_t}^3\right).
\end{aligned}$$

For the $W^{N_0,p}$ estimate of $V^H(t)$, by the Sobolev embedding $H^{N_0+3(\frac{1}{2}-\frac{1}{p})}(\mathbb{R}^3) \hookrightarrow W^{N_0,p}(\mathbb{R}^3)$, we follow the above steps to obtain

$$\|V^H(t)\|_{W^{N_0,p}} \lesssim e^{-c_0 t}\|V(0)\|_{Y} + (1+t)^{-\beta_s}\left(\|V\|_{X_t}^2 + \|V\|_{X_t}^3\right).$$

We shall now prove the decay estimate for low frequencies. To do this, we diagonalize $A$ as follows

$$\begin{aligned}
A(D) &= \begin{pmatrix} 1 & 1 \\ \frac{-\lambda_-(D)}{a(D)} & \frac{-\lambda_+(D)}{a(D)} \end{pmatrix} \begin{pmatrix} -\lambda_-(D) & 0 \\ 0 & -\lambda_+(D) \end{pmatrix} \begin{pmatrix} \lambda_+(D) & a(D) \\ -\lambda_-(D) & -a(D) \end{pmatrix} \frac{1}{2\mathrm{i}b(D)} \\
&:= Q \begin{pmatrix} -\lambda_-(D) & 0 \\ 0 & -\lambda_+(D) \end{pmatrix} Q^{-1}.
\end{aligned} \tag{2.71}$$

Then by applying $Q^{-1}\chi^L$ to the system (2.3) and setting $R = (R_1, R_2)^\mathsf{T} = Q^{-1}\chi^L(D)V$, we obtain the diagonalized equation

$$\partial_t R + \begin{pmatrix} -\lambda_-(D) & 0 \\ 0 & -\lambda_+(D) \end{pmatrix} R = Q^{-1}\chi^L(D)B(V) \tag{2.72}$$

with initial data $R(0) = Q^{-1}\chi^L(D)V(0)$. We thus obtain from the Duhamel formula that

$$R = \begin{pmatrix} e^{\lambda_-(D)t} & 0 \\ 0 & e^{\lambda_+(D)t} \end{pmatrix} R(0) + \int_0^t \begin{pmatrix} e^{\lambda_-(D)(t-\tau)} & 0 \\ 0 & e^{\lambda_+(D)(t-\tau)} \end{pmatrix} Q^{-1}\chi^L(D)B(V)\,\mathrm{d}\tau = \sum_{i=1}^5 J_i, \tag{2.73}$$

where

$$J_1 := \begin{pmatrix} e^{\lambda_-(D)t} & 0 \\ 0 & e^{\lambda_+(D)t} \end{pmatrix} R(0),$$

$$J_2 := \int_0^t \begin{pmatrix} e^{\lambda_-(D)(t-\tau)} & 0 \\ 0 & e^{\lambda_+(D)(t-\tau)} \end{pmatrix} Q^{-1}\chi^L(D)\mathcal{B}(V^H, V)\,\mathrm{d}\tau,$$

$$J_3 := \int_0^t \begin{pmatrix} e^{\lambda_-(D)(t-\tau)} & 0 \\ 0 & e^{\lambda_+(D)(t-\tau)} \end{pmatrix} Q^{-1}\chi^L(D)\mathcal{B}(V^L, V^H)\,\mathrm{d}\tau,$$

$$J_4 := \int_0^t \begin{pmatrix} e^{\lambda_-(D)(t-\tau)} & 0 \\ 0 & e^{\lambda_+(D)(t-\tau)} \end{pmatrix} Q^{-1}\chi^L(D)\mathcal{B}(V^L, V^L)\,\mathrm{d}\tau,$$

$$J_5 := \int_0^t \begin{pmatrix} e^{\lambda_-(D)(t-\tau)} & 0 \\ 0 & e^{\lambda_+(D)(t-\tau)} \end{pmatrix} Q^{-1}\chi^L(D)\mathcal{Q}(V)\,\mathrm{d}\tau.$$

For the term $J_1$, noting that $R(0) = Q^{-1}\chi^L(D)V(0) = \widetilde{\chi}_{\varepsilon,\kappa_0}(D)Q^{-1}\chi^L(D)V(0)$, thus by Corollary 2.7 and Corollary B.8, we have

$$\|J_1\|_{W^{N_0,p}} \lesssim (1+t)^{-\beta_s}\|Q^{-1}\chi^L(D)V(0)\|_{W^{N_0+3(1-\frac{2}{p}),\frac{p}{p-1}}} \lesssim (1+t)^{-\beta_s}\|V(0)\|_Y.$$

For the term $J_2$, we use the Sobolev embedding $H^{\frac{3}{p}}(\mathbb{R}^3) \hookrightarrow L^{\frac{2p}{p-2}}(\mathbb{R}^3)$, Corollary 2.7, Corollary B.8 and Lemma B.10 to obtain

$$\begin{aligned}
\|J_2\|_{W^{N_0,p}} &\lesssim \int_0^t (1+t-\tau)^{-\beta_s}\|Q^{-1}\chi^L(D)\mathcal{B}(V^H,V)\|_{W^{N_0+3(1-\frac{2}{p}),\frac{p}{p-1}}}\,\mathrm{d}\tau \\
&\lesssim \int_0^t (1+t-\tau)^{-\beta_s}\left(\|V^H\|_{H^{N_0+3(1-\frac{2}{p})+1}}\|V\|_{L^{\frac{2p}{p-2}}} + \|V^H\|_{L^2}\|V\|_{W^{N_0+3(1-\frac{2}{p})+1,\frac{2p}{p-2}}}\right)\mathrm{d}\tau \\
&\lesssim \int_0^t (1+t-\tau)^{-\beta_s}(1+\tau)^{-\beta_s}\|V\|_{X_t}^2\,\mathrm{d}\tau \lesssim (1+t)^{-\beta_s}\|V\|_{X_t}^2.
\end{aligned}$$

The estimate for $J_3$ is similar to the one for $J_2$, so we skip it.

For the term $J_5$, it follows from the Sobolev embeddings $W^{N_0,p}(\mathbb{R}^3) \hookrightarrow L^\infty(\mathbb{R}^3)$, $H^{\frac{3}{p}}(\mathbb{R}^3) \hookrightarrow L^{\frac{2p}{p-2}}(\mathbb{R}^3)$, Corollary 2.7, Corollary B.8, Lemma B.7, and Lemma B.10 that

$$\begin{aligned}
\|J_5\|_{W^{N_0,p}} &\lesssim \int_0^t (1+t-\tau)^{-\beta_s}\|Q^{-1}\chi^L(D)\mathcal{Q}(V)\|_{W^{N_0+3(1-\frac{2}{p}),\frac{p}{p-1}}}\,\mathrm{d}\tau \\
&\lesssim \int_0^t (1+t-\tau)^{-\beta_s}\|n_3(D)V\|_{L^\infty}\left(\|V\|_{H^{N_0+3(1-\frac{2}{p})+1}}\|V\|_{L^{\frac{2p}{p-2}}} + \|V\|_{L^2}\|V\|_{W^{N_0+3(1-\frac{2}{p})+1,\frac{2p}{p-2}}}\right)\mathrm{d}\tau \\
&\lesssim \int_0^t (1+t-\tau)^{-\beta_s}(1+\tau)^{-\beta_s}\|V\|_{X_t}^3\,\mathrm{d}\tau \lesssim (1+t)^{-\beta_s}\|V\|_{X_t}^3.
\end{aligned}$$

It remains to estimate $J_4$, which is the most difficult one. To overcome the difficulty of the quadratic nonlinearity, we need to use the normal form method. By the definition of $V^L$, $Q$ and $V^L = QR$, we have

$$V^L = QR = \begin{pmatrix} 1 & 1 \\ \frac{-\lambda_-(D)}{a(D)} & \frac{-\lambda_+(D)}{a(D)} \end{pmatrix}\begin{pmatrix} R_1 \\ R_2 \end{pmatrix} = \begin{pmatrix} R_1 + R_2 \\ \mathrm{i}\frac{b(D)}{a(D)}(R_1 - R_2) + \frac{-\varepsilon\Delta}{a(D)}(R_1 + R_2) \end{pmatrix},$$

which implies that

$$\mathcal{B}(V^L, V^L) = \mathcal{B}(QR, QR) := \begin{pmatrix} \mathcal{B}_1(R,R) \\ \mathcal{B}_2(R,R) \end{pmatrix}$$

with

$$\mathcal{B}_1(R,R) := -a(D)\mathcal{R}^*\left\{\frac{|\nabla|}{a(D)}(R_1+R_2)\mathcal{R}\left[\mathrm{i}\frac{b(D)}{a(D)}(R_1-R_2)+\frac{-\varepsilon\Delta}{a(D)}(R_1+R_2)\right]\right\},$$

$$\mathcal{B}_2(R,R) := |\nabla|\left\{\frac{1}{2}\left|\mathcal{R}\left[\mathrm{i}\frac{b(D)}{a(D)}(R_1-R_2)+\frac{-\varepsilon\Delta}{a(D)}(R_1+R_2)\right]\right|^2+\frac{\gamma-2}{2}\left|\frac{|\nabla|}{a(D)}(R_1+R_2)\right|^2\right\}.$$

Define

$$\mathcal{A}(R,R) = Q^{-1}\mathcal{B}(QR,QR) = \begin{pmatrix} \frac{1}{2\mathrm{i}b(D)}(\lambda_+(D)\mathcal{B}_1+a(D)\mathcal{B}_2) \\ -\frac{1}{2\mathrm{i}b(D)}(\lambda_-(D)\mathcal{B}_1+a(D)\mathcal{B}_2)\end{pmatrix} := \begin{pmatrix}\mathcal{A}_1(R,R)\\ \mathcal{A}_2(R,R)\end{pmatrix}.$$

According to the definition of $R$, we observe that $R_1=\overline{R_2}$. For notational convenience, we write

$$\mathcal{A}_1(R,R) = \frac{a(D)}{2\mathrm{i}b(D)}\sum_{j,k=1}^{2}\sum n_1(D)(n_2(D)R_j\cdot n_2(D)R_k),$$

where the summation is taken over all possible combinations of $n_1(D)$ and $n_2(D)$ defined as

$$n_1(D) := |\nabla| \quad \text{or} \quad \mathrm{i}b(D)\widetilde{\chi}_{\varepsilon,\kappa_0}(D)\mathcal{R}^* \quad \text{or} \quad \varepsilon\Delta\widetilde{\chi}_{\varepsilon,\kappa_0}(D)\mathcal{R}^*,$$

$$n_2(D) := \frac{|\nabla|}{a(D)} \quad \text{or} \quad \widetilde{\chi}_{\varepsilon,\kappa_0}(D)\mathcal{R}\frac{b(D)}{a(D)} \quad \text{or} \quad \widetilde{\chi}_{\varepsilon,\kappa_0}(D)\mathcal{R}\frac{-\varepsilon\Delta}{a(D)}.$$

The properties of $n_1$ and $n_2$ are discussed in Appendix B. Since $\mathcal{A}_2$ exhibits a structure similar to that of $\mathcal{A}_1$, it suffices to restrict our attention to the first component, namely,

$$J_{41} = \int_0^t e^{\lambda_-(D)(t-\tau)}\chi^L(D)\mathcal{A}_1(R,R)\,\mathrm{d}\tau.$$

Let $\widetilde{\chi}_{\varepsilon,\kappa_0}(\xi)=\widetilde{\chi}\big(\sqrt{\frac{\varepsilon}{\kappa_0}}\xi\big)$ where $\widetilde{\chi}$ is defined in (1.15), we have $\widetilde{\chi}_{\varepsilon,\kappa_0}\chi_{\varepsilon,\kappa_0}=\chi_{\varepsilon,\kappa_0}$. Then, $J_{41}$ can be expressed as a sum of terms of the following form:

$$K_{jk} := e^{-\mathrm{i}tb(D)}\mathcal{F}^{-1}\left(\int_0^t\int_{\mathbb{R}^3} e^{-\varepsilon|\xi|^2(t-\tau)}e^{\mathrm{i}\tau b(\xi)}m(\xi,\eta)\chi_{\varepsilon,\kappa_0}(\xi)\widehat{R}_j(\tau,\xi-\eta)\widehat{R}_k(\tau,\eta)\,\mathrm{d}\eta\,\mathrm{d}\tau\right),$$

where the multiplier $m(\xi,\eta)$ is given by

$$m(\xi,\eta) := \widetilde{\chi}_{\varepsilon,\kappa_0}(\xi-\eta)\widetilde{\chi}_{\varepsilon,\kappa_0}(\eta)\widetilde{\chi}_{\varepsilon,\kappa_0}(\xi)\frac{a(\xi)}{2\mathrm{i}b(\xi)}n_1(\xi)n_2(\xi-\eta)n_2(\eta). \tag{2.74}$$

Let us introduce the new variable $W$:

$$W := (W_1,W_2)^{\mathsf{T}} := \begin{pmatrix} e^{\mathrm{i}tb(D)} & 0 \\ 0 & e^{-\mathrm{i}tb(D)}\end{pmatrix}R.$$

Then, it follows from (2.72) that $W$ satisfies the following equation:

$$\partial_t W = \begin{pmatrix} e^{\mathrm{i}tb(D)} & 0 \\ 0 & e^{-\mathrm{i}tb(D)}\end{pmatrix}\left[\varepsilon\Delta R + Q^{-1}\chi^L(D)B(V)\right]. \tag{2.75}$$

By the definition of $W$, we have $W_j(t) := e^{\mathrm{i}(-1)^{j+1}tb(D)}R_j(t)$, thus we get

$$K_{jk} = e^{-\mathrm{i}tb(D)}\mathcal{F}^{-1}\left(\int_0^t\int_{\mathbb{R}^3} e^{-\varepsilon|\xi|^2(t-\tau)}e^{-\mathrm{i}\tau\Phi_{jk}}m(\xi,\eta)\chi_{\varepsilon,\kappa_0}(\xi)\widehat{W}_j(\tau,\xi-\eta)\widehat{W}_k(\tau,\eta)\,\mathrm{d}\eta\,\mathrm{d}\tau\right)$$

with the phase function $\Phi_{jk}$ defined as

$$-\Phi_{jk} := b(\xi)+(-1)^j b(\xi-\eta)+(-1)^k b(\eta). \tag{2.76}$$

Integrating by parts in $\tau$ and making use of (2.75), we deduce that

$$\begin{aligned} K_{jk} &= \mathrm{i}e^{-\mathrm{i}tb(D)}\mathcal{F}^{-1}\left(\int_0^t\int_{\mathbb{R}^3} e^{-\varepsilon|\xi|^2(t-\tau)}\frac{\partial_\tau e^{-\mathrm{i}\tau\Phi_{jk}}}{\Phi_{jk}} m(\xi,\eta)\chi_{\varepsilon,\kappa_0}(\xi)\widehat{W}_j(\tau,\xi-\eta)\widehat{W}_k(\tau,\eta)\,\mathrm{d}\eta\,\mathrm{d}\tau\right) \\ &= \sum_{i=1}^{7} L_{ijk}, \end{aligned} \tag{2.77}$$

where

$$\begin{aligned} L_{1jk} &:= \mathrm{i}\chi_{\varepsilon,\kappa_0}(D)T_{\frac{m}{\Phi_{jk}}}\left(R_j(t),R_k(t)\right), \\ L_{2jk} &:= -\mathrm{i}e^{\varepsilon t\Delta}e^{-\mathrm{i}tb(D)}\chi_{\varepsilon,\kappa_0}(D)T_{\frac{m}{\Phi_{jk}}}\left(R_j(0),R_k(0)\right), \\ L_{3jk} &:= \mathrm{i}\int_0^t e^{\varepsilon(t-\tau)\Delta}e^{\mathrm{i}(t-\tau)b(D)}(\varepsilon\Delta)\chi_{\varepsilon,\kappa_0}(D)T_{\frac{m}{\Phi_{jk}}}\left(R_j(\tau),R_k(\tau)\right)\,\mathrm{d}\tau, \\ L_{4jk} &:= -\mathrm{i}\int_0^t e^{\varepsilon(t-\tau)\Delta}e^{\mathrm{i}(t-\tau)b(D)}\chi_{\varepsilon,\kappa_0}(D)T_{\frac{m}{\Phi_{jk}}}\left(\varepsilon\Delta R_j(\tau),R_k(\tau)\right)\,\mathrm{d}\tau, \\ L_{5jk} &:= -\mathrm{i}\int_0^t e^{\varepsilon(t-\tau)\Delta}e^{\mathrm{i}(t-\tau)b(D)}\chi_{\varepsilon,\kappa_0}(D)T_{\frac{m}{\Phi_{jk}}}\left(R_j(\tau),\varepsilon\Delta R_k(\tau)\right)\,\mathrm{d}\tau, \\ L_{6jk} &:= -\mathrm{i}\int_0^t e^{\varepsilon(t-\tau)\Delta}e^{\mathrm{i}(t-\tau)b(D)}\chi_{\varepsilon,\kappa_0}(D)T_{\frac{m}{\Phi_{jk}}}\left(\widetilde{B}_j(\tau),R_k(\tau)\right)\,\mathrm{d}\tau, \\ L_{7jk} &:= -\mathrm{i}\int_0^t e^{\varepsilon(t-\tau)\Delta}e^{\mathrm{i}(t-\tau)b(D)}\chi_{\varepsilon,\kappa_0}(D)T_{\frac{m}{\Phi_{jk}}}\left(R_j(\tau),\widetilde{B}_k(\tau)\right)\,\mathrm{d}\tau, \end{aligned}$$

with $\widetilde{B} := Q^{-1}\chi^L(D)B(V)$ and $T_{\frac{m}{\Phi_{jk}}}$ denotes the bilinear operator defined by (1.13). In the following, we omit the subscript $jk$ for brevity, and we estimate $L_1$–$L_7$ term by term.

To estimate $L_1$, we introduce following quantities:

$$\begin{aligned} \mathfrak{M}_{jk} :=& \frac{|\xi|^{1-s}|\eta|^{1-s}|\xi-\eta|^{1-s}}{\langle\xi-\eta\rangle^{2\lambda}\langle\eta\rangle^{2\lambda}\Phi_{jk}}\widetilde{\chi}_{\varepsilon,\kappa_0}(\xi-\eta)\widetilde{\chi}_{\varepsilon,\kappa_0}(\eta)\widetilde{\chi}_{\varepsilon,\kappa_0}(\xi)\frac{a(\xi)}{2\mathrm{i}b(\xi)}, \\ \mathfrak{R} :=& (\mathfrak{R}_1,\mathfrak{R}_2)^{\mathsf{T}} \quad \text{with} \quad \widehat{\mathfrak{R}}_i(\xi) := \frac{n_2(\xi)\langle\xi\rangle^{2\lambda}}{|\xi|^{1-s}}\widehat{R}_i(\xi), \quad i=1,2, \\ \mathfrak{B} :=& (\mathfrak{B}_1,\mathfrak{B}_2)^{\mathsf{T}} \quad \text{with} \quad \widehat{\mathfrak{B}}_i(\xi) := \frac{n_2(\xi)\langle\xi\rangle^{2\lambda}}{|\xi|^{1-s}}\widehat{\widetilde{B}}_i(\xi), \quad i=1,2, \end{aligned}$$

where $\lambda := 7/4+\ell s - s/2+\delta = 7/4+(1-\delta)s+\delta$. We split the analysis into two cases: $|\xi|\le 1$ and $|\xi|\ge 1$. For the case $|\xi|\le 1$, utilizing the Plancherel theorem, Corollary B.8, the Sobolev embeddings $\dot{H}^{3(\frac{1}{2}-\frac{1}{p})}(\mathbb{R}^3)\hookrightarrow L^p(\mathbb{R}^3)$ and $H^{\frac{3}{p}-\delta}(\mathbb{R}^3)\hookrightarrow L^{\frac{1}{\frac{1}{2}-\frac{1}{p}+\frac{\delta}{3}}}(\mathbb{R}^3)$, and Lemma B.13, we obtain

$$\begin{aligned} \|L_1\|_{W^{N_0,p}} &\lesssim \left\|\langle\xi\rangle^{N_0+s}|\xi|^{3(\frac{1}{2}-\frac{1}{p})}\int_{\mathbb{R}^3}\mathfrak{M}_{jk}\widehat{\mathfrak{R}}_j(t,\xi-\eta)\widehat{\mathfrak{R}}_k(t,\eta)\,\mathrm{d}\eta\right\|_{L^2_\xi} \\ &\lesssim \left\||\xi|^{1-s}\int_{\mathbb{R}^3}\mathfrak{M}_{jk}\widehat{\mathfrak{R}}_j(t,\xi-\eta)\widehat{\mathfrak{R}}_k(t,\eta)\,\mathrm{d}\eta\right\|_{L^2_\xi} \\ &\lesssim \||\nabla|^{1-s}\mathfrak{R}\|_{L^p}\|\mathfrak{R}\|_{L^{\frac{1}{\frac{1}{2}-\frac{1}{p}+\frac{\delta}{3}}}} \\ &\lesssim \|R\|_{W^{2\lambda,p}}\||\nabla|^{-(1-s)}R\|_{H^{2\lambda+\frac{3}{p}-\delta}} \\ &\lesssim (1+t)^{-\beta_s}\|V\|_{X_t}^2, \end{aligned}$$

where we notice that $2\lambda \le N_0$, $3(1/2-1/p) \ge 1-s$ and $|\xi|^{1-s} \le |\xi-\eta|^{1-s} + |\eta|^{1-s}$. Turning to the high-frequency regime where $|\xi| \ge 1$, applying Lemma B.13, we deduce that

$$\begin{aligned}
\|L_1\|_{W^{N_0,p}} &\lesssim \left\| |\xi|^{N_0+s+3(\frac{1}{2}-\frac{1}{p})} \int_{\mathbb{R}^3} \mathfrak{M}_{jk} \widehat{\mathfrak{R}}_j(t,\xi-\eta) \widehat{\mathfrak{R}}_k(t,\eta)\, \mathrm{d}\eta \right\|_{L^2_\xi} \\
&\lesssim \left\| |\nabla|^{N_0+s+3(\frac{1}{2}-\frac{1}{p})} \mathfrak{R} \right\|_{L^{l_1}} \|\mathfrak{R}\|_{L^{l_2}},
\end{aligned}$$

where we choose $\frac{1}{l_1} = \sigma\left(\frac{1}{2} - \frac{1}{p} + \frac{1-s}{3}\right) + \frac{2\delta}{3}$, $\frac{1}{l_2} = \frac{1}{2} - \sigma\left(\frac{1}{2} - \frac{1}{p} + \frac{1-s}{3}\right) - \frac{\delta}{3}$, and $\sigma := 3/5 - \delta$ such that $2 \le l_1, l_2 \le 3/\delta$. To estimate the first factor, by Lemma B.7 and the Sobolev interpolation, we obtain

$$\begin{aligned}
\left\| |\nabla|^{N_0+s+3(\frac{1}{2}-\frac{1}{p})} \mathfrak{R} \right\|_{L^{l_1}} &\lesssim \left\| |\nabla|^{N_0+2s-1+3(\frac{1}{2}-\frac{1}{p})} \langle\nabla\rangle^{2\lambda} R \right\|_{L^{l_1}} \\
&\lesssim \|\langle\nabla\rangle^{2\lambda} R\|_{W^{N_0-2\lambda,p}}^{1-\sigma} \|\langle\nabla\rangle^{2\lambda} R\|_{H^{\widetilde{N}_1-2\lambda}}^{\sigma} \\
&\lesssim \|R\|_{W^{N_0,p}}^{1-\sigma} \|R\|_{H^N}^{\sigma},
\end{aligned}$$

where $\widetilde{N}_1 := N_0 + s - 1 + (2\lambda + 2s - 2\delta + 1/2)/\sigma \le N$. To bound the second factor, we apply the Bernstein's inequality (Lemma B.2) and the Sobolev interpolation to the low-frequency part to obtain

$$\begin{aligned}
\|\varphi_0(D)\mathfrak{R}\|_{L^{l_2}} &\lesssim \sum_{j=-\infty}^{0} 2^{-j(1-s)} \|\varphi(\frac{D}{2^j}) R\|_{L^{l_2}} \\
&\lesssim \sum_{j=-\infty}^{0} 2^{-j(1-s)} 2^{3j(\frac{1}{l_3}-\frac{1}{l_2})} \|\varphi(\frac{D}{2^j}) R\|_{L^{l_3}} \\
&\lesssim \sum_{j=-\infty}^{0} 2^{j\delta} \left( 2^{-j(1-s)} \|\varphi(\frac{D}{2^j}) R\|_{L^2} \right)^{1-\sigma} \|\varphi(\frac{D}{2^j}) R\|_{L^p}^{\sigma} \\
&\lesssim \||\nabla|^{-(1-s)} R\|_{L^2}^{1-\sigma} \|R\|_{L^p}^{\sigma},
\end{aligned}$$

where we choose $1/l_3 = (1-\sigma)/2 + \sigma/p$. Moreover, for the high-frequency part, by the Sobolev interpolation, we have

$$\begin{aligned}
\|\sum_{j=1}^{\infty} \varphi(\frac{D}{2^j}) \mathfrak{R}\|_{L^{l_2}} &\lesssim \sum_{j=1}^{\infty} 2^{-j(1-s)} \|\varphi(\frac{D}{2^j}) \langle\nabla\rangle^{2\lambda} R\|_{L^{l_2}} \\
&\lesssim \|\langle\nabla\rangle^{2\lambda} R\|_{L^{l_2}} \\
&\lesssim \|\langle\nabla\rangle^{2\lambda} R\|_{H^{\widetilde{N}_2-2\lambda}}^{1-\sigma} \|\langle\nabla\rangle^{2\lambda} R\|_{L^p}^{\sigma} \\
&\lesssim \|R\|_{H^N}^{1-\sigma} \|R\|_{W^{N_0,p}}^{\sigma},
\end{aligned}$$

where $\widetilde{N}_2 := (1-s)\sigma/(1-\sigma) + \delta/(1-\sigma) + 2\lambda \le N$. Collecting these estimates, we deduce from the definition of $R$ and Corollary B.8 that

$$\|L_1\|_{W^{N_0,p}} \lesssim \|V\|_{X_t} \|V\|_{W^{N_0,p}} \lesssim (1+t)^{-\beta_s} \|V\|_{X_t}^2.$$

To estimate $L_2$, applying Corollary 2.7, Corollary B.8, Lemma B.15, and the Sobolev embedding $H^{\frac{3}{p}-\delta}(\mathbb{R}^3) \hookrightarrow L^{\frac{1}{\frac{1}{2}-\frac{1}{p}+\frac{\delta}{3}}}(\mathbb{R}^3)$, we obtain

$$\begin{aligned}
\|L_2\|_{W^{N_0,p}} &\lesssim (1+t)^{-\beta_s} \|T_{\frac{m}{\Phi_{jk}}}(R_j(0), R_k(0))\|_{W^{N_0+3(1-\frac{2}{p}),\frac{p}{p-1}}} \\
&\lesssim (1+t)^{-\beta_s} \||\nabla|^{-(1-s)} R(0)\|_{H^{N_0+3(1-\frac{2}{p})+s+2\lambda}} \||\nabla|^{-(1-s)} R(0)\|_{W^{2\lambda,\frac{1}{\frac{1}{2}-\frac{1}{p}+\frac{\delta}{3}}}} \\
&\lesssim (1+t)^{-\beta_s} \||\nabla|^{-(1-s)} R(0)\|_{H^{N_0+3(1-\frac{2}{p})+s+2\lambda}} \||\nabla|^{-(1-s)} R(0)\|_{H^{2\lambda+\frac{3}{p}-\delta}} \\
&\lesssim (1+t)^{-\beta_s} \|V(0)\|_Y^2.
\end{aligned}$$

Returning to the estimates for $L_3$, $L_4$, and $L_5$, we utilize the identity:

$$\varepsilon\Delta T_{\frac{m}{\Phi_{jk}}}(R_j, R_k) = T_{\frac{m}{\Phi_{jk}}}(\varepsilon\Delta R_j, R_k) + T_{\frac{m}{\Phi_{jk}}}(R_j, \varepsilon\Delta R_k) + 2T_{\frac{m}{\Phi_{jk}}}\left(\varepsilon^{\frac{1}{2}}\nabla R_j, \varepsilon^{\frac{1}{2}}\nabla R_k\right).$$

This identity allows us to combine these terms as follows:

$$L_3 + L_4 + L_5 = 2\mathrm{i} \int_0^t e^{\varepsilon(t-\tau)\Delta} e^{\mathrm{i}(t-\tau)b(D)} \chi_{\varepsilon,\kappa_0}(D) T_{\frac{m}{\Phi_{jk}}}\left(\varepsilon^{\frac{1}{2}}\nabla R_j(\tau), \varepsilon^{\frac{1}{2}}\nabla R_k(\tau)\right) \mathrm{d}\tau := M.$$

Therefore, using Corollary 2.7 and Lemma B.15, we obtain

$$\begin{aligned}
\|M\|_{W^{N_0,p}} &\lesssim \int_0^t (1+t-\tau)^{-\beta_s} \left\| T_{\frac{m}{\Phi_{jk}}}\left(\varepsilon^{\frac{1}{2}}\nabla R_j(\tau), \varepsilon^{\frac{1}{2}}\nabla R_k(\tau)\right) \right\|_{W^{N_0+3(1-\frac{2}{p}),\frac{p}{p-1}}} \mathrm{d}\tau \\
&\lesssim \int_0^t (1+t-\tau)^{-\beta_s} \left\| \frac{\varepsilon^{\frac{1}{2}}\nabla}{|\nabla|^{1-s}} R(\tau) \right\|_{H^{N_0+3(1-\frac{2}{p})+s+2\lambda}} \left\| \frac{\varepsilon^{\frac{1}{2}}\nabla}{|\nabla|^{1-s}} R(\tau) \right\|_{W^{2\lambda,\frac{1}{\frac{1}{2}-\frac{1}{p}+\frac{\delta}{3}}}} \mathrm{d}\tau.
\end{aligned}$$

For the terms in the integral, by (2.73), Corollary 2.6, Lemma B.9, Lemma B.11, Corollary 2.7 and Sobolev embeddings, we deduce for $k = N_0 + 3(1-2/p) + s + 2\lambda$ that

$$\begin{aligned}
\left\| \frac{\varepsilon^{\frac{1}{2}}\nabla}{|\nabla|^{1-s}} R(t) \right\|_{H^k} &\lesssim \left\| e^{\varepsilon t\Delta} e^{-itb(D)} \frac{\varepsilon^{\frac{1}{2}}\nabla}{|\nabla|^{1-s}} R(0) \right\|_{H^k} + \int_0^t \left\| e^{\varepsilon(t-\tau)\Delta} e^{-i(t-\tau)b(D)} \frac{\varepsilon^{\frac{1}{2}}\nabla}{|\nabla|^{1-s}} \widetilde{B}(V)(\tau) \right\|_{H^k} \mathrm{d}\tau \\
&\lesssim (1+t)^{-\frac{1}{2}} \left\| |\nabla|^{-(1-s)} R(0) \right\|_{H^k} + \int_0^t (1+t-\tau)^{-\frac{1}{2}} \left\| |\nabla|^{-(1-s)} \widetilde{B}(V)(\tau) \right\|_{H^k} \mathrm{d}\tau,
\end{aligned} \tag{2.78}$$

and

$$\begin{aligned}
\left\| \frac{\varepsilon^{\frac{1}{2}}\nabla}{|\nabla|^{1-s}} R(t) \right\|_{W^{2\lambda,\frac{1}{\frac{1}{2}-\frac{1}{p}+\frac{\delta}{3}}}} &\lesssim (1+t)^{-\widetilde{\beta}_s-\frac{1}{2}} \||\nabla|^{-(1-s)} R(0)\|_{W^{2\lambda+\frac{6}{p}-2\delta,\frac{1}{\frac{1}{p}+\frac{1}{2}-\frac{\delta}{3}}}} \\
&\quad + \int_0^t (1+t-\tau)^{-\widetilde{\beta}_s-\frac{1}{2}} \||\nabla|^{-(1-s)} \widetilde{B}(V)(\tau)\|_{W^{2\lambda+\frac{6}{p}-2\delta,\frac{1}{\frac{1}{p}+\frac{1}{2}-\frac{\delta}{3}}}} \mathrm{d}\tau,
\end{aligned} \tag{2.79}$$

where

$$\widetilde{\beta}_s = \begin{cases} 3\left(\frac{1}{p} - \frac{\delta}{3}\right) & \text{when } 0 < s \le 1/2, \\ \frac{8}{3}\left(\frac{1}{p} - \frac{\delta}{3}\right) & \text{when } 1/2 < s < 1. \end{cases}$$

Furthermore, by the Sobolev embedding $L^{\frac{1}{\frac{1-s}{3}+\frac{1}{p}+\frac{1}{2}-\frac{\delta}{3}}}(\mathbb{R}^3) \hookrightarrow \dot{W}^{-(1-s),\frac{1}{\frac{1}{p}+\frac{1}{2}-\frac{\delta}{3}}}(\mathbb{R}^3)$ and the Sobolev interpolation, we have

$$\||\nabla|^{-(1-s)}R(0)\|_{W^{2\lambda+\frac{6}{p}-2\delta,\frac{1}{\frac{1}{p}+\frac{1}{2}-\frac{\delta}{3}}}} \lesssim \|R(0)\|_{W^{2\lambda+\frac{6}{p}-2\delta,\frac{1}{\frac{1-s}{3}+\frac{1}{p}+\frac{1}{2}-\frac{\delta}{3}}}} \lesssim \|R(0)\|_{L^1}^{\theta}\|R(0)\|_{H^N}^{1-\theta} \lesssim \|V(0)\|_Y, \tag{2.80}$$

with the choice of $\theta \in (0,1)$ satisfying

$$\frac{1-s}{3}+\frac{1}{p}+\frac{1}{2}-\frac{\delta}{3}-\frac{1}{3}\left(2\lambda+\frac{6}{p}-2\delta\right) = \theta + (1-\theta)\left(\frac{1}{2}-\frac{N}{3}\right),$$

and by Lemma B.11 along with Sobolev embeddings $W^{N_0,p}(\mathbb{R}^3) \hookrightarrow L^{\frac{1}{\frac{1}{p}-\frac{\delta}{3}}}(\mathbb{R}^3)$, $H^N(\mathbb{R}^3) \hookrightarrow L^\infty(\mathbb{R}^3)$, and Lemma B.7, we have

$$\begin{aligned}
\||\nabla|^{-(1-s)}\widetilde{B}(V)(t)\|_{W^{2\lambda+\frac{6}{p}-2\delta,\frac{1}{\frac{1}{p}+\frac{1}{2}-\frac{\delta}{3}}}} &\lesssim (1+\|n_3(D)V\|_{L^\infty})\|V\|_{H^{2\lambda+\frac{6}{p}-2\delta+s}}\|V\|_{L^{\frac{1}{\frac{1}{p}-\frac{\delta}{3}}}} \\
&\lesssim (1+\|V\|_{H^{\frac{3}{2}+\delta}})\|V\|_{H^{2\lambda+\frac{6}{p}-2\delta+s}}\|V\|_{W^{\delta,p}} \\
&\lesssim (1+t)^{-\beta_s}(\|V\|_{X_t}^2+\|V\|_{X_t}^3).
\end{aligned} \tag{2.81}$$

Thus, substituting (2.78), (2.79), (2.80), and (2.81) into the estimate for $M$, we arrive at

$$\begin{aligned}
\|M\|_{W^{N_0,p}} &\lesssim \int_0^t (1+t-\tau)^{-\beta_s}(1+\tau)^{-\widetilde{\beta}_s-1}\,\mathrm{d}\tau\left(\|V(0)\|_Y^2+\|V\|_{X_t}^4+\|V\|_{X_t}^6\right) \\
&\lesssim (1+t)^{-\beta_s}\left(\|V(0)\|_Y^2+\|V\|_{X_t}^4+\|V\|_{X_t}^6\right),
\end{aligned}$$

where $p$ is chosen as in (1.8) to ensure $\widetilde{\beta}_s+1 \geq \beta_s$.

To estimate $L_6$, applying Corollary 2.7, Corollary B.8, Lemma B.15, and the Sobolev embedding $H^{\frac{3}{p}-\delta}(\mathbb{R}^3) \hookrightarrow L^{\frac{1}{\frac{1}{2}-\frac{1}{p}+\frac{\delta}{3}}}(\mathbb{R}^3)$, we have

$$\begin{aligned}
\|L_6\|_{W^{N_0,p}} &\lesssim \int_0^t (1+t-\tau)^{-\beta_s}\left\|T_{\frac{m}{\Phi_{jk}}}\left(\widetilde{B}_j(\tau),R_k(\tau)\right)\right\|_{W^{N_0+3(1-\frac{2}{p}),\frac{p}{p-1}}}\,\mathrm{d}\tau \\
&\lesssim \int_0^t (1+t-\tau)^{-\beta_s}\Big(\||\nabla|^{-(1-s)}\widetilde{B}(\tau)\|_{H^{N_0+3(1-\frac{2}{p})+s+2\lambda}}\||\nabla|^{-(1-s)}R(\tau)\|_{H^{2\lambda+\frac{3}{p}-\delta}} \\
&\qquad + \||\nabla|^{-(1-s)}\widetilde{B}(\tau)\|_{H^{2\lambda+\frac{3}{p}-\delta}}\||\nabla|^{-(1-s)}R(\tau)\|_{H^{N_0+3(1-\frac{2}{p})+s+2\lambda}}\Big)\,\mathrm{d}\tau.
\end{aligned}$$

By Lemma B.11 along with Sobolev embeddings $W^{N_0,p}(\mathbb{R}^3) \hookrightarrow L^\infty(\mathbb{R}^3)$, $H^N(\mathbb{R}^3) \hookrightarrow L^\infty(\mathbb{R}^3)$, and Lemma B.7, we find that for $k = N_0+3(1-2/p)+s+2\lambda$ or $k = 2\lambda+3/p-\delta$,

$$\begin{aligned}
\||\nabla|^{-(1-s)}\widetilde{B}(\tau)\|_{H^k} &\lesssim (1+\|n_3(D)V\|_{L^\infty})\,\|n_3(D)V\|_{L^\infty}\,\|V\|_{H^{k+s}} \\
&\lesssim (1+\tau)^{-\beta_s}\left(\|V\|_{X_t}^2+\|V\|_{X_t}^3\right),
\end{aligned} \tag{2.82}$$

where we notice that $N \geq N_0+3(1-2/p)+2s+2\lambda$. Applying (2.81) and (2.82), we obtain the estimate for $L_6$

$$\|L_6\|_{W^{N_0,p}} \lesssim (1+t)^{-\beta_s}\left(\|V\|_{X_t}^3+\|V\|_{X_t}^4\right).$$

The analysis of $L_7$ follows a strictly analogous procedure to that of $L_6$ and is therefore omitted. Combining the estimates for $L_1$–$L_7$, and using the *a priori* assumption (2.44), we obtain

$$\|J_4\|_{W^{N_0,p}} \lesssim (1+t)^{-\beta_s}\left(\|V(0)\|_Y+\|V\|_{X_t}^2+\|V\|_{X_t}^4\right).$$

Together with the estimate for $J_1$, $J_2$, $J_3$, and $J_5$, this gives

$$\|V^L(t)\|_{W^{N_0,p}} \lesssim (1+t)^{-\beta_s}\left(\|V(0)\|_Y + \|V\|_{X_t}^2 + \|V\|_{X_t}^4\right).$$

Combining the $W^{N_0,p}$ and $H^{N-1}$ estimates for $V^H$ with the $W^{N_0,p}$ estimate for $V^L$, we obtain the desired estimate for $V$. This completes the proof of Proposition 2.12. □

2.2.3. *Estimates in negative Sobolev spaces.* First, we complete the last piece of uniform estimates in $X_t$ by proving the estimate in $\dot{H}^{-(1-s)}$.

**Proposition 2.13.** *If* (2.44) *holds, then we have*

$$\|U(t)\|_{\dot{H}^{-(1-s)}} \lesssim \|U_0\|_{\dot{H}^{-(1-s)}} + \|U\|_{X_t}^2 + \|U\|_{X_t}^3. \tag{2.83}$$

*Proof.* Given the equivalence of the norm $\|U\|_{\dot{H}^{-(1-s)}} \simeq \|V\|_{\dot{H}^{-(1-s)}}$, it is sufficient to establish the bound for the variable $V$. By integrating (2.9) with respect to time and applying Corollary 2.6, we deduce

$$\|V(t)\|_{\dot{H}^{-(1-s)}} \lesssim \|V(0)\|_{\dot{H}^{-(1-s)}} + \int_0^t \|B(V)(\tau)\|_{\dot{H}^{-(1-s)}}\,\mathrm{d}\tau.$$

To handle the nonlinear integral term, an application of Lemma B.11 yields

$$\|B(V)\|_{\dot{H}^{-(1-s)}} \lesssim (1+\|n_3(D)V\|_{L^\infty})\|n_3(D)V\|_{L^\infty}\|V\|_{H^s}.$$

Coupling this estimate with the Sobolev embedding $W^{\frac{3}{p}+\delta}(\mathbb{R}^3) \hookrightarrow L^\infty(\mathbb{R}^3)$ (for any $\delta > 0$) ensures the temporal integrability. Specifically, we have

$$\int_0^t \|B(V)(\tau)\|_{\dot{H}^{-(1-s)}}\,\mathrm{d}\tau \lesssim \int_0^t (1+\tau)^{-\beta_s}\,\mathrm{d}\tau\left(\|V\|_{X_t}^2 + \|V\|_{X_t}^3\right) \lesssim \|V\|_{X_t}^2 + \|V\|_{X_t}^3.$$

Synthesizing these estimates yields the desired inequality (2.83), thereby completing the proof. □

To recover the $L^2$ time-decay estimate, we further prove the estimate in $\dot{H}^{-\vartheta}$ with $0 < \vartheta < 3/2$.

**Proposition 2.14.** *If* (2.44) *holds, then we have*

$$\frac{d}{dt}\sqrt{\mathbb{E}^{-\vartheta}(t)} \lesssim \|(\nabla\mathfrak{c},\nabla u)\|_{L^2}\|(\mathfrak{c},u)\|_{L^{\frac{3}{\vartheta}}} + \begin{cases}\|\varrho\|_{L^2}\|u\|_{L^{\frac{3}{\vartheta+s-1}}}, & \text{for } 1<\vartheta+s<5/2,\\ \|(\varrho,u)\|_{H^1}\|(\varrho,u)\|_{L^\infty}, & \text{for } \vartheta+s\le 1,\end{cases} \tag{2.84}$$

*where the negative-order Sobolev energy functional is defined as*

$$\mathbb{E}^{-\vartheta}(t) := \||\nabla|^{-\vartheta}(\mathfrak{c},u,|\nabla|^s\phi)\|_{L^2}^2.$$

*Proof.* Applying $|\nabla|^{-\vartheta}$ to (2.46), taking the $L^2$ inner product with $|\nabla|^{-\vartheta}\mathfrak{U}$, and performing integration by parts yields

$$\frac{1}{2}\frac{d}{dt}\||\nabla|^{-\vartheta}\mathfrak{U}\|_{L^2}^2 + 2\varepsilon\||\nabla|^{-\vartheta}\nabla u\|_{L^2}^2 = \sum_{i=1}^5 \mathbf{I}_i,$$

where

$$\mathbf{I}_1 := -(|\nabla|^{-\vartheta}(u\cdot\nabla\mathfrak{c}),|\nabla|^{-\vartheta}\mathfrak{c}),\ \ \mathbf{I}_2 := -(|\nabla|^{-\vartheta}(u\cdot\nabla u),|\nabla|^{-\vartheta}u),\ \ \mathbf{I}_3 := -\frac{\gamma-1}{2}(|\nabla|^{-\vartheta}(\mathfrak{c}\,\mathrm{div}\,u),|\nabla|^{-\vartheta}\mathfrak{c}),$$

$$\mathbf{I}_4 := -\frac{\gamma-1}{2}(|\nabla|^{-\vartheta}(\mathfrak{c}\nabla\mathfrak{c}),|\nabla|^{-\vartheta}u),\ \ \mathbf{I}_5 := -(|\nabla|^{-\vartheta}\nabla\phi,|\nabla|^{-\vartheta}u).$$

For the terms $\mathbf{I}_1$–$\mathbf{I}_4$, we use Hölder's inequality and the Hardy-Littlewood-Sobolev inequality (Lemma B.3) to get

$$\begin{aligned}\sum_{i=1}^{4}|\mathbf{I}_i| &\lesssim \||\nabla|^{-\vartheta}(\mathfrak{c},u)\|_{L^2}\||\nabla|^{-\vartheta}\left((\mathfrak{c},u)\nabla(\mathfrak{c},u)\right)\|_{L^2}\\ &\lesssim \||\nabla|^{-\vartheta}(\mathfrak{c},u)\|_{L^2}\|(\mathfrak{c},u)\nabla(\mathfrak{c},u)\|_{L^{\frac{1}{\frac12+\frac{\vartheta}{3}}}}\\ &\lesssim \||\nabla|^{-\vartheta}(\mathfrak{c},u)\|_{L^2}\|(\nabla\mathfrak{c},\nabla u)\|_{L^2}\|(\mathfrak{c},u)\|_{L^{\frac{3}{\vartheta}}}.\end{aligned}$$

Regarding $\mathbf{I}_5$, we apply integration by parts and use $(1.7)_1$ and $(1.7)_3$ to deduce

$$\mathbf{I}_5=-\frac{1}{2}\frac{d}{dt}\||\nabla|^{-\vartheta}|\nabla|^{s}\phi\|_{L^2}^2\underbrace{-(|\nabla|^{-\vartheta}\phi,|\nabla|^{-\vartheta}\mathrm{div}\,(\varrho u))}_{\mathbf{I}_6}.$$

The term $\mathbf{I}_6$ is estimated in two regimes according to the index $\vartheta+s$:

**Case 1:** $1<\vartheta+s<5/2$. By Hölder's inequality and the Hardy-Littlewood-Sobolev inequality (Lemma B.3), we have

$$|\mathbf{I}_6|\lesssim \||\nabla|^{-\vartheta}|\nabla|^{s}\phi\|_{L^2}\|\varrho u\|_{L^{\frac{1}{\frac12+\frac{\vartheta+s-1}{3}}}}\lesssim \||\nabla|^{-\vartheta}|\nabla|^{s}\phi\|_{L^2}\|\varrho\|_{L^2}\|u\|_{L^{\frac{3}{\vartheta+s-1}}}.$$

**Case 2:** $\vartheta+s\le 1$. Applying Hölder's inequality and Lemma B.4, one has

$$|\mathbf{I}_6|\lesssim \||\nabla|^{-\vartheta}|\nabla|^{s}\phi\|_{L^2}\||\nabla|^{-\vartheta+1-s}(\varrho u)\|_{L^2}\lesssim \||\nabla|^{-\vartheta}|\nabla|^{s}\phi\|_{L^2}\|(\varrho,u)\|_{H^1}\|(\varrho,u)\|_{L^\infty}.$$

Collecting the above estimates, omitting the non-negative dissipation term on the left-hand side and dividing both sides by $\sqrt{\mathbb{E}^{-\vartheta}(t)}$, the desired differential inequality (2.84) is proved. □

2.2.4. *Proof of Theorem 1.3.* In this part, we prove Theorem 1.3.

First, following a standard iteration argument, one can show that system (1.7) admits a unique solution in $C([0,T_\varepsilon),H^3)$ for some $T_\varepsilon>0$ and that if the initial data are in $H^l$, $l\ge 3$, then this additional regularity also propagates on $[0,T_\varepsilon)$.

To extend this local solution globally in time, we define

$$T_*:=\sup\{T\in[0,T_\varepsilon):\|U\|_{X_T}\le\delta_1\}.$$

Then choosing $\delta_0$ and $\delta_1$ sufficiently small, by Corollary 2.9, Proposition 2.12, and Proposition 2.13, we have

$$\|U\|_{X_T}\le C\left(\|U_0\|_Y+\|U\|_{X_T}^{\frac32}+\|U\|_{X_T}^4\right)\le C(\delta_0+\delta_1^{\frac32})\le\frac{\delta_1}{2},$$

for all $T<T_*$, A standard bootstrap argument implies that $T_*=T_\varepsilon=\infty$. This ensures the global existence of the solution $U(t)$ as well as the uniform-in-time bound $\|U\|_{X_{T_*}}\le\delta_1$.

Having established the global-in-time existence and uniform bounds, we now derive the $\varepsilon$-dependent $L^2$ decay rates using the negative Sobolev estimates from Proposition 2.14.

**Step 1.** By the Sobolev interpolation, for $k=0,1$, we have

$$\|\nabla^k U\|_{L^2}\le\||\nabla|^{-\vartheta}U\|_{L^2}^{\frac{1}{k+1+\vartheta}}\|\nabla\nabla^k U\|_{L^2}^{\frac{k+\vartheta}{k+1+\vartheta}},$$

which implies that for $k=0,1$,

$$\begin{aligned}\mathbb{E}_N^k(t)&\lesssim\left\{\mathbb{E}_N^k(t)+\mathbb{E}^{-\vartheta}(t)\right\}^{\frac{1}{k+1+\vartheta}}\left\{\mathbb{D}_N^k(t)\right\}^{\frac{k+\vartheta}{k+1+\vartheta}}\\ &\lesssim\left\{\sup_{t\ge0}\mathbb{E}_N^k(t)+\sup_{t\ge0}\mathbb{E}^{-\vartheta}(t)\right\}^{\frac{1}{k+1+\vartheta}}\left\{\mathbb{D}_N^k(t)\right\}^{\frac{k+\vartheta}{k+1+\vartheta}}.\end{aligned}$$

Combining this with the energy estimates in Corollary 2.11 and using the equivalence $\mathbb{E}_N^k \simeq \bar{\mathbb{E}}_N^k$ $(k=0,1)$, we obtain that for $k=0,1$,

$$\begin{aligned}&\frac{d}{dt}\bar{\mathbb{E}}_N^k(t) + C\varepsilon\left\{\sup_{t\ge 0}\mathbb{E}_N^k(t) + \sup_{t\ge 0}\mathbb{E}^{-\vartheta}(t)\right\}^{-\frac{1}{k+\vartheta}}\left\{\bar{\mathbb{E}}_N^k(t)\right\}^{1+\frac{1}{k+\vartheta}} \\ &\lesssim \left(\|(\varrho,\mathfrak{c},u)\|_{W^{1,\infty}} + \mathbf{1}_{0<s<\frac12}\|\nabla\varrho\|_{L^{\frac{3}{s}}}\right)\bar{\mathbb{E}}_N^k(t).\end{aligned} \tag{2.85}$$

Denote

$$I(t) = \exp\left\{-\int_0^t \left(\|(\varrho,\mathfrak{c},u)(\tau)\|_{W^{1,\infty}} + \mathbf{1}_{0<s<\frac12}\|\nabla\varrho(\tau)\|_{L^{\frac{3}{s}}}\right)\mathrm{d}\tau\right\}.$$

From (2.64) and (2.44), we have $I(t)\in[e^{-C\delta_1},1]$. Solving the differential inequality (2.85) directly yields that for $k=0,1$,

$$\mathbb{E}_N^k(t) \lesssim \left[\left\{\mathbb{E}_N^k(0)\right\}^{-\frac{1}{k+\vartheta}} + \varepsilon t\left\{\sup_{t\ge 0}\mathbb{E}_N^k(t) + \sup_{t\ge 0}\mathbb{E}^{-\vartheta}(t)\right\}^{-\frac{1}{k+\vartheta}}\right]^{-(k+\vartheta)}. \tag{2.86}$$

**Step 2.** We derive the decay rates by dividing the proof into two cases: $0<s\le 1/2$ and $1/2<s<1$.
**Case 1:** $0<s\le 1/2$. We establish the decay rates sequentially according to the range of $\vartheta$.

First, for $0<\vartheta<\frac12$, we have $\vartheta+s<1$. A direct calculation gives

$$\|(\nabla\mathfrak{c},\nabla u)\|_{L^2}\|(\mathfrak{c},u)\|_{L^{\frac{3}{\vartheta}}} + \|(\varrho,u)\|_{H^1}\|(\varrho,u)\|_{L^\infty} \lesssim (1+t)^{-(\frac32-\vartheta)} + (1+t)^{-\frac32(1-\frac2p)}.$$

Combining this with (2.84) yields $\mathbb{E}^{-\vartheta}(t)\lesssim 1$. Substituting this bound into (2.86) and using (2.44), we obtain

$$\mathbb{E}_N^k(t) \lesssim (1+\varepsilon t)^{-(k+\vartheta)} \quad \text{for } \ k=0,1. \tag{2.87}$$

Next, for $\frac12\le\vartheta\le 1$, utilizing (2.87) and the embedding $\dot{H}^{-\vartheta}\cap L^2\subset \dot{H}^{-\vartheta'}$ for any fixed $\vartheta'\in(0,\frac12)$, we have

$$\begin{aligned}&\|(\nabla\mathfrak{c},\nabla u)\|_{L^2}\|(\mathfrak{c},u)\|_{L^{\frac{3}{\vartheta}}} + \|(\varrho,u)\|_{H^1}\|(\varrho,u)\|_{L^\infty} + \|\varrho\|_{L^2}\|u\|_{L^{\frac{3}{\vartheta+s-1}}} \\ &\lesssim (1+\varepsilon t)^{-\frac{1+\vartheta'}{2}}(1+t)^{-(\frac32-\vartheta)} + (1+t)^{-\frac32\left(1-\frac2p\right)} + (1+\varepsilon t)^{-\frac{\vartheta'}{2}}(1+t)^{-(\frac52-\vartheta-s)}.\end{aligned}$$

The temporal integration is dominated by the first term with decay exponent $\beta=\frac32-\vartheta$. Combining this with (2.84) and applying Lemma B.16, we obtain

$$\mathbb{E}^{-\vartheta}(t) \lesssim \begin{cases}\left[\ln(1+\varepsilon^{-1})\right]^2, & \text{if } \vartheta=\frac12,\\ \varepsilon^{-(2\vartheta-1)}, & \text{if } \frac12<\vartheta\le 1.\end{cases} \tag{2.88}$$

Inserting (2.88) into (2.86), we can unify the two scenarios into a single, compact decay estimate

$$\mathbb{E}_N^k(t) \lesssim \left(1+\varepsilon^{\Gamma_1(k,\vartheta)}t\right)^{-(k+\vartheta)} \quad \text{for } \ k=0,1, \tag{2.89}$$

where the effective exponent $\Gamma_1\ge 1$ is defined by

$$\Gamma_1(k,\vartheta) := \begin{cases}1^+, & \text{if } \vartheta=\frac12,\\ 1+\frac{2\vartheta-1}{k+\vartheta}, & \text{if } \frac12<\vartheta\le 1.\end{cases}$$

Finally, we extend the analysis to the regime $1<\vartheta<\frac32$. Using the embedding $\dot{H}^{-\vartheta}\cap L^2\subset\dot{H}^{-1}$ along with the condition $\vartheta+s>1$, we can use the estimate (2.89) for $\vartheta'=1$. This yields

$$\begin{aligned}&\|(\nabla\mathfrak{c},\nabla u)\|_{L^2}\|(\mathfrak{c},u)\|_{L^{\frac{3}{\vartheta}}} + \|\varrho\|_{L^2}\|u\|_{L^{\frac{3}{\vartheta+s-1}}} \\ &\lesssim (1+\varepsilon^{\frac32}t)^{-1}(1+t)^{-(\frac32-\vartheta)} + (1+\varepsilon^2 t)^{-\frac12}(1+t)^{-(\frac52-\vartheta-s)}.\end{aligned}$$

The temporal integration is dominated by the first term with decay exponent $\beta = \frac{3}{2} - \vartheta$. Combining this with (2.84) and applying Lemma B.16, we obtain $\mathbb{E}^{-\vartheta}(t) \lesssim \varepsilon^{-\frac{3}{2}(2\vartheta-1)}$. Inserting this into (2.86) gives the exact decay rate

$$\mathbb{E}_N^k(t) \lesssim \left(1 + \varepsilon^{1+\frac{3(2\vartheta-1)}{2(k+\vartheta)}} t\right)^{-(k+\vartheta)} \quad \text{for } \ k = 0, 1. \tag{2.90}$$

Summarizing (2.89) and (2.90), the overall decay rates for $0 < s \le 1/2$ and $k = 0, 1$ can be unified as

$$\mathbb{E}_N^k(t) \lesssim \left(1 + \varepsilon^{\Upsilon_1} t\right)^{-(k+\vartheta)},$$

where the effective exponent $\Upsilon_1 \ge 1$ is defined by

$$\Upsilon_1 := \begin{cases} 1, & \text{if } 0 < \vartheta < \frac{1}{2}, \\ 1^+, & \text{if } \vartheta = \frac{1}{2}, \\ 1 + \frac{2\vartheta-1}{\vartheta}, & \text{if } \frac{1}{2} < \vartheta \le 1, \\ 1 + \frac{3(2\vartheta-1)}{2\vartheta}, & \text{if } 1 < \vartheta < \frac{3}{2}. \end{cases} \tag{2.91}$$

Consequently, we obtain the following enhanced decay estimate for $\varrho$ as follows

$$\|\varrho\|_{L^2} \lesssim \||\nabla|^s \phi\|_{L^2}^{1-s} \|\nabla |\nabla|^s \phi\|_{L^2}^{s} \lesssim \left(1 + \varepsilon^{\Upsilon_1} t\right)^{-\left(\frac{s}{2} + \frac{\vartheta}{2}\right)}.$$

**Case 2:** $1/2 < s < 1$**.** Due to the slower algebraic decay rate, the critical integrability thresholds are changed. We establish the decay rates for $\vartheta \in (0, \frac{3}{2})$ step by step.

First, for $0 < \vartheta < \frac{3}{8}$, we have

$$\|(\nabla \mathfrak{c}, \nabla u)\|_{L^2} \|(\mathfrak{c}, u)\|_{L^{\frac{3}{\vartheta}}} + \|(\varrho, u)\|_{H^1} \|(\varrho, u)\|_{L^\infty} \lesssim (1+t)^{-\frac{4}{3}\left(1 - \frac{2\vartheta}{3}\right)} + (1+t)^{-\frac{4}{3}\left(1 - \frac{2}{p}\right)}.$$

In view of (2.84), it follows that $\mathbb{E}^{-\vartheta}(t) \lesssim 1$. Substituting this bound into (2.86) and using (2.44), we obtain

$$\mathbb{E}_N^k(t) \lesssim (1 + \varepsilon t)^{-(k+\vartheta)} \quad \text{for } \ k = 0, 1, \tag{2.92}$$

and subsequently,

$$\|\varrho\|_{L^2} \lesssim (1 + \varepsilon t)^{-\left(\frac{s}{2} + \frac{\vartheta}{2}\right)}.$$

Next, for $\frac{3}{8} \le \vartheta \le 1$, utilizing (2.92) and the embedding $\dot{H}^{-\vartheta} \cap L^2 \subset \dot{H}^{-\frac{1}{4}}$, we have

$$\begin{aligned} &\|(\nabla \mathfrak{c}, \nabla u)\|_{L^2} \|(\mathfrak{c}, u)\|_{L^{\frac{3}{\vartheta}}} + \|(\varrho, u)\|_{H^1} \|(\varrho, u)\|_{L^\infty} + \|\varrho\|_{L^2} \|u\|_{L^{\frac{3}{\vartheta+s-1}}} \\ &\lesssim (1 + \varepsilon t)^{-\frac{5}{8}} (1+t)^{-\frac{4}{3}\left(1 - \frac{2\vartheta}{3}\right)} + (1+t)^{-\frac{4}{3}\left(1 - \frac{2}{p}\right)} + (1 + \varepsilon t)^{-\frac{1}{8} - \frac{s}{2}} (1+t)^{-\left(\frac{20}{9} - \frac{8}{9}(\vartheta + s)\right)}. \end{aligned}$$

The temporal integration is dominated by the first term with decay exponent $\beta = \frac{4}{3} - \frac{8\vartheta}{9}$. Combining this with (2.84) and applying Lemma B.16, we obtain

$$\mathbb{E}^{-\vartheta}(t) \lesssim \begin{cases} \left[\ln(1 + \varepsilon^{-1})\right]^2, & \text{if } \vartheta = \frac{3}{8}, \\ \varepsilon^{-\frac{16\vartheta - 6}{9}}, & \text{if } \frac{3}{8} < \vartheta \le 1. \end{cases} \tag{2.93}$$

Inserting (2.93) into (2.86) leads to

$$\mathbb{E}_N^k(t) \lesssim \left(1 + \varepsilon^{\Gamma_2(k,\vartheta)} t\right)^{-(k+\vartheta)} \quad \text{for } \ k = 0, 1, \quad \text{and } \ \|\varrho\|_{L^2} \lesssim \left(1 + \varepsilon^{\Gamma_2(0,\vartheta)} t\right)^{-\left(\frac{s}{2} + \frac{\vartheta}{2}\right)}, \tag{2.94}$$

where the effective exponent $\Gamma_2 \ge 1$ is defined by

$$\Gamma_2(k, \vartheta) := \begin{cases} 1^+, & \text{if } \vartheta = \frac{3}{8}, \\ 1 + \frac{16\vartheta - 6}{9(k+\vartheta)}, & \text{if } \frac{3}{8} < \vartheta \le 1. \end{cases}$$

Finally, we extend the analysis to the regime $1 < \vartheta < \frac{3}{2}$. Using the embedding $\dot{H}^{-\vartheta} \cap L^2 \subset \dot{H}^{-1}$ along with the condition $\vartheta + s > 1$, we can use the estimate (2.94) for $\vartheta = 1$. This yields

$$\begin{aligned}&\|(\nabla \mathfrak{c}, \nabla u)\|_{L^2} \|(\mathfrak{c}, u)\|_{L^{\frac{3}{\vartheta}}} + \|\varrho\|_{L^2} \|u\|_{L^{\frac{3}{\vartheta+s-1}}} \\ &\lesssim (1 + \varepsilon^{\frac{14}{9}} t)^{-1} (1+t)^{-\frac{4}{3}\left(1 - \frac{2\vartheta}{3}\right)} + (1 + \varepsilon^{\frac{19}{9}} t)^{-\frac{1+s}{2}} (1+t)^{-\left(\frac{20}{9} - \frac{8}{9}(\vartheta+s)\right)}.\end{aligned}$$

The temporal integration is governed by the decay exponents $\beta_1 = \frac{4}{3} - \frac{8\vartheta}{9}$ and $\beta_2 = \frac{20}{9} - \frac{8}{9}(\vartheta + s)$. Combining this with (2.84) and applying Lemma B.16, we obtain

$$\mathbb{E}^{-\vartheta}(t) \lesssim \varepsilon^{-\frac{2}{81} \max\left\{14(8\vartheta - 3), 19(8\vartheta + 8s - 11)\right\}}.$$

Substituting this bound into (2.86) gives

$$\mathbb{E}_N^k(t) \lesssim \left(1 + \varepsilon^{\Gamma_3(k,\vartheta,s)} t\right)^{-(k+\vartheta)} \quad \text{for } \ k = 0, 1, \quad \text{and} \quad \|\varrho\|_{L^2} \lesssim \left(1 + \varepsilon^{\Gamma_3(0,\vartheta,s)} t\right)^{-\left(\frac{s}{2} + \frac{\vartheta}{2}\right)}, \tag{2.95}$$

where the effective exponent $\Gamma_3 \geq 1$ is defined by

$$\Gamma_3(k, \vartheta, s) := 1 + \frac{2 \max\left\{14(8\vartheta - 3), 19(8\vartheta + 8s - 11)\right\}}{81(k + \vartheta)}.$$

Summarizing (2.92), (2.94) and (2.95), the overall decay rates for $1/2 < s < 1$ and $k = 0, 1$ can be unified as

$$\mathbb{E}_N^k(t) \lesssim \left(1 + \varepsilon^{\Upsilon_2} t\right)^{-(k+\vartheta)}, \quad \text{and} \quad \|\varrho\|_{L^2} \lesssim \left(1 + \varepsilon^{\Upsilon_2} t\right)^{-\left(\frac{s}{2} + \frac{\vartheta}{2}\right)},$$

where the effective exponent $\Upsilon_2 \geq 1$ is defined by

$$\Upsilon_2 := \begin{cases} 1, & \text{if } 0 < \vartheta < \frac{3}{8}, \\ 1^+, & \text{if } \vartheta = \frac{3}{8}, \\ 1 + \frac{16\vartheta - 6}{9\vartheta}, & \text{if } \frac{3}{8} < \vartheta \leq 1, \\ 1 + \frac{2 \max\left\{14(8\vartheta - 3), 19(8\vartheta + 8s - 11)\right\}}{81\vartheta}, & \text{if } 1 < \vartheta < \frac{3}{2}. \end{cases} \tag{2.96}$$

Thus, we conclude the proof of Theorem 1.3.

## 3. Correction System

In this section, we study the correction system (1.10). In order to derive the basic $L^2$ energy estimates, we introduce the new quantity $\mathfrak{n}$

$$\mathfrak{n} := \frac{2}{\gamma - 1} \left[ (\rho + n)^{\frac{\gamma-1}{2}} - \rho^{\frac{\gamma-1}{2}} \right].$$

Then $\mathfrak{V} := (\mathfrak{n}, v)$ satisfies a symmetric hyperbolic–parabolic system

$$\partial_t \mathfrak{V} + \sum_{i=1}^{3} \mathbb{M}_i(\mathfrak{U}, \mathfrak{V}) \partial_{x_i} \mathfrak{V} - \varepsilon \mathbf{N} \mathcal{L} \mathfrak{V} + (0, \nabla \psi)^\intercal = S, \tag{3.1}$$

where $\mathbf{N}$ is defined in (2.47) and

$$\mathbb{M}_i(\mathfrak{U}, \mathfrak{V}) := \begin{pmatrix} u_i + v_i & \frac{\gamma-1}{2}(\mathfrak{n} + \mathfrak{c} + \frac{2}{\gamma-1}) e_i \\ \frac{\gamma-1}{2}(\mathfrak{n} + \mathfrak{c} + \frac{2}{\gamma-1}) e_i^\intercal & (u_i + v_i)\mathbb{I}_3 \end{pmatrix},$$

$$S := \begin{pmatrix} -v \cdot \nabla \mathfrak{c} - \frac{\gamma-1}{2} \mathfrak{n} \operatorname{div} u \\ -v \cdot \nabla u - \frac{\gamma-1}{2} \mathfrak{n} \nabla \mathfrak{c} + \varepsilon \left(\frac{1}{\rho+n} - 1\right) \mathcal{L}(v + u) \end{pmatrix}.$$

We introduce the energy and dissipation functionals order by order. For the base case $k = 0$, the zero-order energy functional is defined as

$$\mathcal{E}^0(t) := \|(\mathfrak{n}, v, |\nabla|^s\psi)\|_{L^2}^2.$$

For the higher-order derivatives ($k \geq 1$), the energy functional is defined according to the fractional exponent $s$

$$\mathcal{E}^k(t) := \begin{cases} \left\|(\nabla^k\mathfrak{n}, \nabla^k v, \nabla^k|\nabla|^s\psi)\right\|_{L^2}^2 & \text{if } 1/2 \leq s < 1, \\ \left\|\left(\nabla^k\mathfrak{n}, \nabla^k v, \dfrac{1}{\sqrt{\rho+n}}\nabla^k|\nabla|^s\psi\right)\right\|_{L^2}^2 & \text{if } 0 < s < 1/2. \end{cases}$$

Based on these, for an integer $m \geq 1$, we introduce the total energy functional

$$\mathcal{E}_m(t) := \sum_{k=0}^{m} \mathcal{E}^k(t),$$

and the corresponding total dissipation rate functional

$$\mathcal{D}_m(t) := \|(\nabla v, \operatorname{div} v)\|_{H^m}^2 + \|(\nabla\mathfrak{n}, \nabla|\nabla|^s\psi)\|_{H^{m-1}}^2.$$

For the subsequent analysis, we impose the following *a priori* assumption: there exist a constant $0 < \delta_1 \ll 1$ and an integer $N_1$ with $3 \leq N_1 \leq N_0 - 3$ such that

$$\sup_{0\leq t\leq T} \mathcal{E}_{N_1}(t) \leq \varepsilon^2\delta_1^2. \tag{3.2}$$

Under the *a priori* assumptions (2.44) and (3.2), it is easy to see the following equivalence relation

$$\begin{cases} \|\mathfrak{n}\|_{L^\infty} \simeq \|n\|_{L^\infty}, & \\ \|\mathfrak{n}\|_{H^k} \simeq \|n\|_{H^k}, & \text{for any } k \in \mathbb{N} \text{ and } 2 \leq k \leq N_1, \\ \|\nabla\mathfrak{n}\|_{H^k} \simeq \|\nabla n\|_{H^k}, & \text{for any } k \in \mathbb{N} \text{ and } 1 \leq k \leq N_1 - 1, \\ \||\nabla|^{-\vartheta}\mathfrak{n}\|_{L^2} \simeq \||\nabla|^{-\vartheta}n\|_{L^2}, & \text{for any } 0 < \vartheta \leq 1/2. \end{cases} \tag{3.3}$$

We shall use the smallness condition (3.2) as a bootstrap hypothesis. More precisely, for a fixed $T > 0$, we assume that a smooth solution $(\mathfrak{n}, v, \psi)$ on $[0, T]$ satisfies (2.44); all estimates derived below are uniform with respect to $T$.

**Proposition 3.1.** *If* (3.2) *holds, then we have*

$$\frac{1}{2}\frac{d}{dt}\mathcal{E}^0(t) + \varepsilon\|\nabla v\|_{L^2}^2 \lesssim \|(\varrho, u)\|_{W^{2,\infty}}\left\{\mathcal{E}_{N_1}(t) + \varepsilon^2\mathbb{E}_N(t)\right\} + \sqrt{\mathcal{E}_{N_1}(t)}\mathcal{D}_{N_1}(t).$$

*Proof.* Taking the $L^2$ inner product of $\mathfrak{V}$ with the system (3.1), we obtain

$$\frac{1}{2}\frac{d}{dt}\|(\mathfrak{n}, v)\|_{L^2}^2 + \varepsilon(\|\nabla v\|_{L^2}^2 + \|\operatorname{div} v\|_{L^2}^2) = \sum_{i=1}^{6}\mathfrak{I}_i,$$

where

$$\mathfrak{I}_1 := -\int_{\mathbb{R}^3} [(u+v)\cdot\nabla\mathfrak{n}]\,\mathfrak{n} + [(u+v)\cdot\nabla v]\cdot v\,\mathrm{d}x,$$

$$\mathfrak{I}_2 := -\frac{\gamma-1}{2}\int_{\mathbb{R}^3}\left[\left(\mathfrak{n} + \mathfrak{c} + \frac{2}{\gamma-1}\right)\operatorname{div} v\right]\mathfrak{n} + \left[\left(\mathfrak{n} + \mathfrak{c} + \frac{2}{\gamma-1}\right)\nabla\mathfrak{n}\right]\cdot v\,\mathrm{d}x,$$

$$\mathfrak{I}_3 := -\frac{\gamma-1}{2}\int_{\mathbb{R}^3}(\mathfrak{n}\operatorname{div} u)\mathfrak{n} + (\mathfrak{n}\nabla\mathfrak{c})\cdot v\,\mathrm{d}x, \qquad \mathfrak{I}_4 := -\int_{\mathbb{R}^3}(v\cdot\nabla\mathfrak{c})\mathfrak{n} + (v\cdot\nabla u)\cdot v\,\mathrm{d}x,$$

$$\mathfrak{I}_5 := \varepsilon\int_{\mathbb{R}^3}\Big(\frac{1}{\rho+n} - 1\Big)(\mathcal{L}v + \mathcal{L}u)\cdot v\,\mathrm{d}x, \qquad \mathfrak{I}_6 := -\int_{\mathbb{R}^3}\nabla\psi\cdot v\,\mathrm{d}x.$$

In view of integration by parts and Hölder's inequality, we have

$$
\begin{aligned}
|\mathfrak{I}_1| &\lesssim \|\nabla u\|_{L^\infty}\|(\mathfrak{n},v)\|_{L^2}^2 + \|\nabla v\|_{L^2}\|(\mathfrak{n},v)\|_{L^6}\|(\mathfrak{n},v)\|_{L^3},\\
|\mathfrak{I}_2| &\lesssim \|\mathfrak{c}\|_{L^\infty}\|(\nabla\mathfrak{n},\nabla v)\|_{L^2}\|(\mathfrak{n},v)\|_{L^2} + \|(\nabla\mathfrak{n},\nabla v)\|_{L^2}\|(\mathfrak{n},v)\|_{L^6}\|(\mathfrak{n},v)\|_{L^3},\\
|\mathfrak{I}_3| + |\mathfrak{I}_4| &\lesssim \|(\nabla\mathfrak{c},\nabla u)\|_{L^\infty}\|(\mathfrak{n},v)\|_{L^2}^2,
\end{aligned}
$$

which, combined with (2.49), implies that

$$
\sum_{i=1}^4 |\mathfrak{I}_i| \lesssim \|(\varrho,u)\|_{W^{1,\infty}}\mathcal{E}_{N_1}(t) + \sqrt{\mathcal{E}_{N_1}(t)}\mathcal{D}_{N_1}(t).
$$

For $\mathfrak{I}_5$, making use of the facts

$$
\frac{1}{\rho+n} - 1 = -\frac{\varrho+n}{\rho+n}, \quad \text{and} \quad \frac{1}{2} \le |\rho+n| \le \frac{3}{2}, \tag{3.4}
$$

we have from Hölder's inequality and Sobolev embeddings that

$$
\begin{aligned}
|\mathfrak{I}_5| &\lesssim \varepsilon\|\varrho\|_{L^\infty}(\|\nabla^2 v\|_{L^2} + \|\nabla^2 u\|_{L^2})\|v\|_{L^2} + \varepsilon\|n\|_{L^3}\|\nabla^2 v\|_{L^2}\|v\|_{L^6} + \varepsilon\|n\|_{L^2}\|\nabla^2 u\|_{L^\infty}\|v\|_{L^2}\\
&\lesssim \|(\varrho,u)\|_{W^{2,\infty}}\|(n,v)\|_{H^2}^2 + \varepsilon^2\|\varrho\|_{L^\infty}\|\nabla^2 u\|_{L^2}^2 + \|n\|_{H^1}\|\nabla v\|_{H^1}^2\\
&\lesssim \|(\varrho,u)\|_{W^{2,\infty}}\mathcal{E}_{N_1}(t) + \varepsilon^2\|\varrho\|_{L^\infty}\mathbb{E}_N(t) + \sqrt{\mathcal{E}_{N_1}(t)}\mathcal{D}_{N_1}(t).
\end{aligned}
$$

For $\mathfrak{I}_6$, using integration by parts and substituting $(1.10)_1$ and $(1.10)_3$ into $\mathfrak{I}_6$, we obtain

$$
\mathfrak{I}_6 = \int_{\mathbb{R}^3}\psi \operatorname{div} v\,\mathrm{d}x = -\frac{1}{2}\frac{d}{dt}\||\nabla|^s\psi\|_{L^2}^2 + \int_{\mathbb{R}^3}(\varrho v + nu + nv)\cdot\nabla\psi\,\mathrm{d}x.
$$

Using Hölder's inequality and Lemma B.4, we have

$$
\left|\int_{\mathbb{R}^3}(\varrho v + nu)\cdot\nabla\psi\,\mathrm{d}x\right| \lesssim \||\nabla|^s\psi\|_{L^2}\|(\varrho v, nu)\|_{H^1} \lesssim \|(\varrho,u)\|_{W^{1,\infty}}\mathcal{E}_{N_1}(t).
$$

As for the remaining part, with the help of Hölder's inequality, Sobolev embeddings and (3.3), we have that for $0 < s < 1/2$,

$$
\begin{aligned}
\left|\int_{\mathbb{R}^3} nv\cdot\nabla\psi\,\mathrm{d}x\right| &\lesssim \|n\|_{L^6}\|v\|_{L^3}\|\nabla\psi\|_{L^2}\\
&\lesssim \|\nabla n\|_{L^2}\|v\|_{L^2}^{\frac{1}{2}}\|\nabla v\|_{L^2}^{\frac{1}{2}}\||\nabla|^s\psi\|_{L^2}^{s}\|\nabla|\nabla|^s\psi\|_{L^2}^{1-s}\\
&\lesssim \sqrt{\mathcal{E}_{N_1}(t)}\mathcal{D}_{N_1}(t),
\end{aligned}
$$

and for $1/2 \le s < 1$,

$$
\begin{aligned}
\left|\int_{\mathbb{R}^3} nv\cdot\nabla\psi\,\mathrm{d}x\right| &= \int_{\mathbb{R}^3}(-\Delta)^s\psi v\cdot\nabla\psi\,\mathrm{d}x \lesssim \|(-\Delta)^s\psi\|_{L^2}\|v\|_{L^6}\|\nabla\psi\|_{L^3}\\
&\lesssim \||\nabla|^s\psi\|_{L^2}^{1-s}\|\nabla|\nabla|^s\psi\|_{L^2}^{s}\||\nabla|^s\psi\|_{L^2}^{\frac{2s-1}{2}}\|\nabla|\nabla|^s\psi\|_{L^2}^{\frac{3-2s}{2}}\|\nabla v\|_{L^2}\\
&\lesssim \sqrt{\mathcal{E}_{N_1}(t)}\mathcal{D}_{N_1}(t).
\end{aligned}
$$

In the above estimates, we have used the fact

$$
\|\nabla|\nabla|^s\psi\|_{L^2}^2 \lesssim \min\{\mathcal{E}_{N_1}(t), \mathcal{D}_{N_1}(t)\}.
$$

It is worth noting that the Sobolev interpolation is employed to ensure sufficient dissipation. Specifically, for $0 < s < 1/2$, the resulting dissipation index is given by $1 + 1/2 + 1 - s = 5/2 - s$. Given that this index exceeds 2, we can strategically allocate the "quadratic part" as the dissipation rate functional, while the remaining part contributes to the energy functional. Collecting all the estimates of $\mathfrak{I}_1$–$\mathfrak{I}_6$, we complete the proof of the result. □

**Proposition 3.2.** *If* (3.2) *holds, then for* $1 \le k \le N_1$*, we have*

$$\begin{aligned}
&\frac{1}{2}\frac{d}{dt}\mathcal{E}^k(t) + \varepsilon\|\nabla^k\nabla v\|_{L^2}^2 \\
&\lesssim (\|(\varrho, u, |\nabla|^s\phi)\|_{W^{N_1+2,\infty}} + \|\varrho\|_{W^{N_1-1,6}}^2)\left\{\mathcal{E}_{N_1}(t) + \varepsilon^2\mathbb{E}_N(t)\right\} \\
&\quad + \left\{\varepsilon\eta + \varepsilon\sqrt{\mathbb{E}_N(t)} + \sqrt{\mathcal{E}_{N_1}(t)}\right\}\mathcal{D}_{N_1}(t).
\end{aligned}$$

*Proof.* Applying $\nabla^k$ to (3.1), then taking the $L^2$ inner product with $\nabla^k\mathfrak{V}$, we have

$$\frac{1}{2}\frac{d}{dt}\|(\nabla^k\mathfrak{n}, \nabla^k v)\|_{L^2}^2 + \varepsilon\|\nabla^k\nabla v\|_{L^2}^2 + \varepsilon\|\nabla^k \mathrm{div}\, v\|_{L^2}^2 = \sum_{i=1}^{9}\mathfrak{J}_i,$$

where

$$\begin{aligned}
\mathfrak{J}_1 &:= -\int_{\mathbb{R}^3}\left[(u+v)\cdot\nabla\nabla^k\mathfrak{n}\right]\nabla^k\mathfrak{n} + \left[(u+v)\cdot\nabla\nabla^k v\right]\cdot\nabla^k v\,\mathrm{d}x,\\
\mathfrak{J}_2 &:= -\int_{\mathbb{R}^3}\left[\nabla^k, (u+v)\cdot\nabla\right]\mathfrak{n}\,\nabla^k\mathfrak{n} + \left[\nabla^k, (u+v)\cdot\nabla\right]v\cdot\nabla^k v\,\mathrm{d}x,\\
\mathfrak{J}_3 &:= -\frac{\gamma-1}{2}\int_{\mathbb{R}^3}\left(\mathfrak{n}\mathrm{div}\,\nabla^k v\right)\nabla^k\mathfrak{n} + \left(\mathfrak{n}\nabla\nabla^k\mathfrak{n}\right)\cdot\nabla^k v\,\mathrm{d}x,\\
\mathfrak{J}_4 &:= -\frac{\gamma-1}{2}\int_{\mathbb{R}^3}\left[\nabla^k, \mathfrak{n}\right]\mathrm{div}\, v\nabla^k\mathfrak{n} + \left[\nabla^k, \mathfrak{n}\cdot\nabla\right]\mathfrak{n}\cdot\nabla^k v\,\mathrm{d}x,\\
\mathfrak{J}_5 &:= -\frac{\gamma-1}{2}\int_{\mathbb{R}^3}\nabla^k\left[\left(\mathfrak{c} + \frac{2}{\gamma-1}\right)\mathrm{div}\, v\right]\nabla^k\mathfrak{n} + \nabla^k\left[\left(\mathfrak{c} + \frac{2}{\gamma-1}\right)\nabla\mathfrak{n}\right]\cdot\nabla^k v\,\mathrm{d}x,\\
\mathfrak{J}_6 &:= -\frac{\gamma-1}{2}\int_{\mathbb{R}^3}\nabla^k(\mathfrak{n}\mathrm{div}\, u)\nabla^k\mathfrak{n} + \nabla^k(\mathfrak{n}\nabla\mathfrak{c})\cdot\nabla^k v\,\mathrm{d}x,\\
\mathfrak{J}_7 &:= -\int_{\mathbb{R}^3}\nabla^k(v\cdot\nabla\mathfrak{c})\nabla^k\mathfrak{n} + \nabla^k(v\cdot\nabla u)\cdot\nabla^k v\,\mathrm{d}x,\\
\mathfrak{J}_8 &:= \varepsilon\int_{\mathbb{R}^3}\nabla^k\left[\left(\frac{1}{\rho+n} - 1\right)(\mathcal{L}v + \mathcal{L}u)\right]\cdot\nabla^k v\,\mathrm{d}x,\\
\mathfrak{J}_9 &:= -\int_{\mathbb{R}^3}\nabla^k\nabla\psi\cdot\nabla^k v\,\mathrm{d}x.
\end{aligned}$$

For $\mathfrak{J}_1$, by integration by parts and Hölder's inequality, we have

$$|\mathfrak{J}_1| \lesssim \|(\nabla u, \nabla v)\|_{L^\infty}\|(\nabla^k\mathfrak{n}, \nabla^k v)\|_{L^2}^2.$$

For $\mathfrak{J}_2$, a direct calculation using the Leibniz rule, Hölder's inequality and Sobolev embeddings gives

$$\begin{aligned}
|\mathfrak{J}_2| &= \left|\sum_{k_1=1}^{k} C_k^{k_1}\int_{\mathbb{R}^3}\nabla^{k_1}(u+v)\cdot\left(\nabla\nabla^{k-k_1}\mathfrak{n}\nabla^k\mathfrak{n} + \nabla\nabla^{k-k_1}v\cdot\nabla^k v\right)\,\mathrm{d}x\right|\\
&\lesssim \|\nabla u\|_{W^{k-1,\infty}}\|(\mathfrak{n}, v)\|_{H^k}^2 + \|v\|_{H^{\max\{k,3\}}}\|(\nabla\mathfrak{n}, \nabla v)\|_{H^{k-1}}^2.
\end{aligned}$$

Similar to $\mathfrak{J}_1$, the term $\mathfrak{J}_3$ is bounded by

$$|\mathfrak{J}_3| \lesssim \|\nabla\mathfrak{n}\|_{L^\infty}\|(\nabla^k\mathfrak{n}, \nabla^k v)\|_{L^2}^2.$$

As for $\mathfrak{J}_4$, similar to $\mathfrak{J}_2$, we have

$$
\begin{aligned}
|\mathfrak{J}_4| &= \frac{\gamma-1}{2}\left|\sum_{k_1=1}^{k} C_k^{k_1}\int_{\mathbb{R}^3}\nabla^{k_1}\mathfrak{n}\left(\nabla^{k-k_1}\mathrm{div}v\nabla^k\mathfrak{n}+\nabla\nabla^{k-k_1}\mathfrak{n}\cdot\nabla^k v\right)\mathrm{d}x\right| \\
&\lesssim \|\mathfrak{n}\|_{H^{\max\{k,3\}}}\|(\nabla\mathfrak{n},\nabla v)\|_{H^{k-1}}^2.
\end{aligned}
$$

By a similar argument, we can estimate $\mathfrak{J}_5$ as follows:

$$
\begin{aligned}
|\mathfrak{J}_5| &= \frac{\gamma-1}{2}\left|\sum_{k_1=0}^{k} C_k^{k_1}\int_{\mathbb{R}^3}\nabla^{k_1}\mathfrak{c}\left(\nabla^{k-k_1}\mathrm{div}v\nabla^k\mathfrak{n}+\nabla\nabla^{k-k_1}\mathfrak{n}\cdot\nabla^k v\right)\mathrm{d}x\right| \\
&\lesssim \|\nabla\mathfrak{c}\|_{W^{k-1,\infty}}\|(\mathfrak{n},v)\|_{H^k}^2.
\end{aligned}
$$

Similarly, for the terms $\mathfrak{J}_6$ and $\mathfrak{J}_7$, it is easy to obtain

$$
|\mathfrak{J}_6|+|\mathfrak{J}_7|\lesssim \|(\nabla\mathfrak{c},\nabla u)\|_{W^{k,\infty}}\|(\mathfrak{n},v)\|_{H^k}^2.
$$

For $\mathfrak{J}_8$, we apply integration by parts to deduce

$$
\mathfrak{J}_8=-\varepsilon\int_{\mathbb{R}^3}\nabla^{k-1}\left[\Big(\frac{1}{\rho+n}-1\Big)\mathcal{L}v\right]\cdot\nabla^k\nabla v\,\mathrm{d}x-\varepsilon\int_{\mathbb{R}^3}\nabla^{k}\left[\Big(\frac{1}{\rho+n}-1\Big)\mathcal{L}u\right]\cdot\nabla^k v\,\mathrm{d}x,
$$

which is then decomposed into $\mathfrak{J}_{81}$ and $\mathfrak{J}_{82}$. By using Hölder's inequality, Sobolev embeddings, Lemma B.4, Lemma B.5 and (3.3), we have

$$
\begin{aligned}
|\mathfrak{J}_{81}| &\lesssim \varepsilon\left\|\nabla^{k-1}\left[\Big(\frac{1}{\rho+n}-1\Big)\mathcal{L}v\right]\right\|_{L^2}\|\nabla^k\nabla v\|_{L^2} \\
&\lesssim \varepsilon\|\nabla v\|_{H^k}\left(\|\varrho+n\|_{L^\infty}\|\nabla^2 v\|_{H^{k-1}}+\|\varrho+n\|_{W^{k-1,6}}\|\nabla^2 v\|_{L^3}\right) \\
&\lesssim \varepsilon\|(\varrho,\mathfrak{n})\|_{L^\infty}\|\nabla v\|_{H^k}^2+(\varepsilon\eta+\|\nabla v\|_{H^2})(\|\nabla\mathfrak{n}\|_{H^{k-1}}^2+\|\nabla v\|_{H^k}^2)+\varepsilon\|\nabla v\|_{H^2}^2\|\varrho\|_{W^{k-1,6}}^2,
\end{aligned}
$$

and

$$
\begin{aligned}
|\mathfrak{J}_{82}| &\lesssim \varepsilon\left\|\nabla^{k}\left[\Big(\frac{1}{\rho+n}-1\Big)\mathcal{L}u\right]\right\|_{L^2}\|\nabla^k v\|_{L^2} \\
&\lesssim \varepsilon\|\varrho+n\|_{H^k}\|\nabla^2 u\|_{W^{k,\infty}}\|v\|_{H^k} \\
&\lesssim \|\nabla^2 u\|_{W^{k,\infty}}\|(n,v)\|_{H^k}^2+\varepsilon^2\|\nabla^2 u\|_{W^{k,\infty}}\|\varrho\|_{H^k}^2.
\end{aligned}
$$

Then we deduce

$$
|\mathfrak{J}_8|\lesssim \left\{\|\nabla^2 u\|_{W^{k,\infty}}+\|\varrho\|_{W^{k-1,6}}^2\right\}\left\{\mathcal{E}_{N_1}(t)+\varepsilon^2\mathbb{E}_N(t)\right\}+\left(\varepsilon\eta+\varepsilon\sqrt{\mathbb{E}_N(t)}+\sqrt{\mathcal{E}_{N_1}(t)}\right)\mathcal{D}_{N_1}(t).
$$

It thus remains to estimate $\mathfrak{J}_9$, which is treated separately according to whether $0<s<1/2$ or $1/2\le s<1$.

**Case 1:** $1/2\le s<1$. Substituting $(1.10)_1$ and $(1.10)_3$ into $\mathfrak{J}_9$ and integrating by parts, we have

$$
\mathfrak{J}_9=-\frac{1}{2}\frac{d}{dt}\|\nabla^k|\nabla|^s\psi\|_{L^2}^2+\int_{\mathbb{R}^3}\nabla^k\nabla\psi\cdot\nabla^k(\varrho v+nu+nv)\,\mathrm{d}x.
$$

Using Sobolev embedding, the last term of the above equality can be estimated as

$$
\begin{aligned}
&\left|\int_{\mathbb{R}^3}\nabla^k\nabla\psi\cdot\nabla^k(\varrho v+nu+nv)\,\mathrm{d}x\right| \\
&\lesssim \|\nabla^k\nabla\psi\|_{L^2}\left(\|\varrho\|_{W^{k,\infty}}\|v\|_{H^k}+\|u\|_{W^{k,\infty}}\|n\|_{H^k}+\|(\nabla n,\nabla v)\|_{H^{\max\{k-1,1\}}}^2\right) \\
&\lesssim \|(\varrho,u)\|_{W^{k,\infty}}\mathcal{E}_{N_1}(t)+\sqrt{\mathcal{E}_{N_1}(t)}\mathcal{D}_{N_1}(t),
\end{aligned}
$$

where we use the Sobolev interpolation to obtain

$$\|\nabla^k \nabla \psi\|_{L^2} \lesssim \|\nabla^k |\nabla|^s \psi\|_{L^2}^{\frac{2s-1}{s}} \|\nabla^k (-\Delta)^s \psi\|_{L^2}^{\frac{1-s}{s}} \lesssim \sqrt{\mathcal{E}_{N_1}(t)}. \tag{3.5}$$

**Case 2:** $0 < s < 1/2$**.** In this case, the estimate (3.5) is no longer valid and requires a different argument. We thus use integration by parts to get

$$\mathfrak{J}_9 = \int_{\mathbb{R}^3} \nabla^k |\nabla|^s \psi \nabla^k |\nabla|^{-s} \mathrm{div}\, v \, \mathrm{d}x = \mathfrak{J}_{91} + \mathfrak{J}_{92} + \mathfrak{J}_{93},$$

where

$$\begin{aligned}
\mathfrak{J}_{91} &:= -\int_{\mathbb{R}^3} \frac{1}{\rho + n} \left[\nabla^k |\nabla|^{-s}, \rho\right] \mathrm{div}\, v \nabla^k |\nabla|^s \psi \, \mathrm{d}x, \\
\mathfrak{J}_{92} &:= -\int_{\mathbb{R}^3} \frac{1}{\rho + n} \left[\nabla^k |\nabla|^{-s}, n\right] \mathrm{div}\, v \nabla^k |\nabla|^s \psi \, \mathrm{d}x, \\
\mathfrak{J}_{93} &:= \int_{\mathbb{R}^3} \frac{1}{\rho + n} \nabla^k |\nabla|^{-s} \left[(\rho + n) \mathrm{div}\, v\right] \nabla^k |\nabla|^s \psi \, \mathrm{d}x.
\end{aligned}$$

Since $k - s > 0$, Lemma B.4 is applicable, which allows us to estimate $\mathfrak{J}_{91}$ and $\mathfrak{J}_{92}$ as follows

$$\begin{aligned}
|\mathfrak{J}_{91}| &\lesssim \left(\|\nabla^k |\nabla|^{-s} \varrho\|_{L^\infty} \|\nabla v\|_{L^2} + \|\nabla \varrho\|_{L^\infty} \|\nabla^k |\nabla|^{-s} v\|_{L^2}\right) \|\nabla^k |\nabla|^s \psi\|_{L^2} \\
&\lesssim \|(\varrho, |\nabla|^s \phi)\|_{W^{k,\infty}} \mathcal{E}_{N_1}(t), \\
|\mathfrak{J}_{92}| &\lesssim \left(\|\nabla^k |\nabla|^{-s} n\|_{L^2} \|\nabla v\|_{L^\infty} + \|\nabla n\|_{L^\infty} \|\nabla^k |\nabla|^{-s} v\|_{L^2}\right) \|\nabla^k |\nabla|^s \psi\|_{L^2} \\
&\lesssim \sqrt{\mathcal{E}_{N_1}(t)} \mathcal{D}_{N_1}(t).
\end{aligned}$$

For $\mathfrak{J}_{93}$, using $(1.10)_1$ and $(1.10)_3$, we rewrite

$$\begin{aligned}
\mathfrak{J}_{93} &= -\int_{\mathbb{R}^3} \frac{1}{\rho + n} \nabla^k |\nabla|^{-s} \left[\partial_t n + v \cdot \nabla(\rho + n) + \mathrm{div}\,(nu)\right] \nabla^k |\nabla|^s \psi \, \mathrm{d}x \\
&= -\frac{1}{2} \frac{d}{dt} \int_{\mathbb{R}^3} \frac{1}{\rho + n} |\nabla^k |\nabla|^s \psi|^2 \, \mathrm{d}x + \mathfrak{J}_{93}^1 + \mathfrak{J}_{93}^2 + \mathfrak{J}_{93}^3,
\end{aligned}$$

where

$$\begin{aligned}
\mathfrak{J}_{93}^1 &:= -\frac{1}{2} \int_{\mathbb{R}^3} \frac{1}{(\rho + n)^2} \partial_t (\rho + n) |\nabla^k |\nabla|^s \psi|^2 \, \mathrm{d}x, \\
\mathfrak{J}_{93}^2 &:= -\int_{\mathbb{R}^3} \frac{1}{\rho + n} \nabla^k |\nabla|^{-s} \left[v \cdot \nabla(\rho + n)\right] \nabla^k |\nabla|^s \psi \, \mathrm{d}x, \\
\mathfrak{J}_{93}^3 &:= -\int_{\mathbb{R}^3} \frac{1}{\rho + n} \nabla^k |\nabla|^{-s} \mathrm{div}\,(nu) \nabla^k |\nabla|^s \psi \, \mathrm{d}x.
\end{aligned}$$

Next, we estimate them term by term. For $\mathfrak{J}_{93}^1$, it follows from (2.44), (3.2) and

$$\partial_t (\rho + n) + \mathrm{div}\left[(\rho + n)(u + v)\right] = 0,$$

that

$$\begin{aligned}
|\mathfrak{J}_{93}^1| &= \left|\frac{1}{2} \int_{\mathbb{R}^3} \frac{1}{\rho + n} \mathrm{div}\,(u + v) |\nabla^k |\nabla|^s \psi|^2 + \frac{1}{(\rho + n)^2} (u + v) \cdot \nabla(\rho + n) |\nabla^k |\nabla|^s \psi|^2 \, \mathrm{d}x\right| \\
&\lesssim \|(\nabla u, \nabla v)\|_{L^\infty} \|\nabla^k |\nabla|^s \psi\|_{L^2}^2 + \|(u, v)\|_{L^\infty} \|(\nabla \varrho, \nabla n)\|_{L^\infty} \|\nabla^k |\nabla|^s \psi\|_{L^2}^2 \\
&\lesssim \|(\varrho, u)\|_{W^{1,\infty}} \mathcal{E}_{N_1}(t) + \sqrt{\mathcal{E}_{N_1}(t)} \mathcal{D}_{N_1}(t).
\end{aligned}$$

For $\mathfrak{J}_{93}^2$, we first decompose it into $\mathfrak{J}_{93}^{21}$ and $\mathfrak{J}_{93}^{22}$, and then provide their respective estimates as follows

$$
\begin{aligned}
|\mathfrak{J}_{93}^{21}| &= \left| \int_{\mathbb{R}^3} \frac{1}{\rho+n} v \cdot \nabla\nabla^k |\nabla|^{-s}(\rho+n) \nabla^k |\nabla|^s \psi \,\mathrm{d}x \right| \\
&= \left| \int_{\mathbb{R}^3} \frac{1}{\rho+n} v \cdot \nabla\nabla^k |\nabla|^{s}(\phi+\psi) \nabla^k |\nabla|^s \psi \,\mathrm{d}x \right| \\
&\lesssim \|v\|_{L^2} \|\nabla^{k+1}|\nabla|^s \phi\|_{L^\infty} \|\nabla^k |\nabla|^s \psi\|_{L^2} + \left(\|(\nabla\varrho, \nabla n)\|_{L^\infty} \|v\|_{L^2} + \|\nabla v\|_{L^\infty}\right) \|\nabla^k |\nabla|^s \psi\|_{L^2}^2 \\
&\lesssim \|(\varrho, |\nabla|^s \phi)\|_{W^{k+1,\infty}} \mathcal{E}_{N_1}(t) + \sqrt{\mathcal{E}_{N_1}(t)} \mathcal{D}_{N_1}(t), \\
|\mathfrak{J}_{93}^{22}| &= \left| \int_{\mathbb{R}^3} \frac{1}{\rho+n} \left[\nabla^k |\nabla|^{-s}, v\right] \cdot \nabla(\rho+n) \nabla^k |\nabla|^s \psi \,\mathrm{d}x \right| \\
&\lesssim \|\nabla^k |\nabla|^{-s} v\|_{L^2} \|(\nabla\varrho, \nabla n)\|_{L^\infty} \|\nabla^k |\nabla|^s \psi\|_{L^2} \\
&\quad + \left( \|\nabla v\|_{L^2} \|\nabla^k |\nabla|^{-s} \varrho\|_{L^\infty} + \|\nabla v\|_{L^\infty} \|\nabla^k |\nabla|^{-s} n\|_{L^2} \right) \|\nabla^k |\nabla|^s \psi\|_{L^2} \\
&\lesssim \|(\varrho, |\nabla|^s \phi)\|_{W^{k,\infty}} \mathcal{E}_{N_1}(t) + \sqrt{\mathcal{E}_{N_1}(t)} \mathcal{D}_{N_1}(t),
\end{aligned}
$$

where we have used integration by parts, $(1.7)_3$, $(1.10)_3$ and (3.4). Thus we have

$$
|\mathfrak{J}_{93}^2| \lesssim \|(\varrho, |\nabla|^s \phi)\|_{W^{k+1,\infty}} \mathcal{E}_{N_1}(t) + \sqrt{\mathcal{E}_{N_1}(t)} \mathcal{D}_{N_1}(t).
$$

As for $\mathfrak{J}_{93}^3$, we have

$$
\begin{aligned}
\mathfrak{J}_{93}^3 = &- \int_{\mathbb{R}^3} \frac{1}{\rho+n} \nabla^k |\nabla|^{-s} (n \operatorname{div} u) \nabla^k |\nabla|^s \psi \,\mathrm{d}x \\
&- \int_{\mathbb{R}^3} \frac{1}{\rho+n} u \cdot \nabla\nabla^k |\nabla|^{-s} n \nabla^k |\nabla|^s \psi \,\mathrm{d}x \\
&- \int_{\mathbb{R}^3} \frac{1}{\rho+n} [\nabla^k |\nabla|^{-s}, u] \cdot \nabla n \nabla^k |\nabla|^s \psi \,\mathrm{d}x.
\end{aligned}
$$

Similarly to the estimate of $\mathfrak{J}_{93}^2$, we obtain

$$
|\mathfrak{J}_{93}^3| \lesssim \|u\|_{W^{k+1,\infty}} \mathcal{E}_{N_1}(t) + \sqrt{\mathcal{E}_{N_1}(t)} \mathcal{D}_{N_1}(t).
$$

Collecting all the estimates of $\mathfrak{J}_1$–$\mathfrak{J}_9$, and summing over $1 \le k \le N_1$, we thus conclude our proof. □

Next, we derive the dissipation estimate for $n$ and $|\nabla|^s \psi$. We first decompose the pressure term as follows

$$
\frac{1}{\gamma-1} \nabla \left[(\rho+n)^{\gamma-1} - \rho^{\gamma-1}\right] = \nabla n + \nabla \mathfrak{H}_1 + \nabla \mathfrak{H}_2
$$

with

$$
\mathfrak{H}_1 := \frac{1}{\gamma-1} \left[(\rho+n)^{\gamma-1} - \rho^{\gamma-1} - (\gamma-1)\rho^{\gamma-2} n\right], \quad \mathfrak{H}_2 := \left(\rho^{\gamma-2} - 1\right) n,
$$

which allow us to rewrite $(1.10)_2$ as

$$
\partial_t v + u \cdot \nabla v + v \cdot (\nabla u + \nabla v) + \nabla n + \nabla \psi = \varepsilon \frac{1}{\rho+n} \mathcal{L} v + \varepsilon \left( \frac{1}{\rho+n} - 1 \right) \mathcal{L} u - \nabla \mathfrak{H}_1 - \nabla \mathfrak{H}_2. \quad (3.6)
$$

**Proposition 3.3.** *If* (3.2) *holds, then we have*

$$
\begin{aligned}
&\frac{d}{dt}\sum_{k=0}^{N_1-1}\int_{\mathbb{R}^3}\nabla^k v\cdot\nabla^k\nabla n\,\mathrm{d}x+\|(\nabla n,\nabla|\nabla|^s\psi)\|_{H^{N_1-1}}^2\\
&\lesssim\|\nabla v\|_{H^{N_1}}^2+\left\{\varepsilon\sqrt{\mathbb{E}_N(t)}+\sqrt{\mathcal{E}_{N_1}(t)}\right\}\mathcal{D}_{N_1}(t)\\
&\quad+\left\{\|(\varrho,u)\|_{W^{N_1+1,\infty}}+\|\varrho\|_{W^{N_1-1,6}}^2\right\}\left\{\mathcal{E}_{N_1}(t)+\varepsilon^2\mathbb{E}_N(t)\right\}.
\end{aligned}\tag{3.7}
$$

*Proof.* Applying $\nabla^k$ with $k\le N_1-1$ to (3.6), taking the $L^2$ inner product with $\nabla^k\nabla n$ and using $(1.10)_1$ and $(1.10)_3$ yields

$$
\frac{d}{dt}\int_{\mathbb{R}^3}\nabla^k v\cdot\nabla^k\nabla n\,\mathrm{d}x+\|\nabla^k(\nabla n,\nabla|\nabla|^s\psi)\|_{L^2}^2=\sum_{i=1}^{7}\mathfrak{L}_i,
$$

where

$$
\begin{aligned}
&\mathfrak{L}_1:=-\int_{\mathbb{R}^3}\nabla^k v\cdot\nabla^k\nabla\mathrm{div}\,(\rho v+nu)\,\mathrm{d}x, && \mathfrak{L}_2:=-\int_{\mathbb{R}^3}\nabla^k(u\cdot\nabla v+v\cdot\nabla u)\cdot\nabla^k\nabla n\,\mathrm{d}x,\\
&\mathfrak{L}_3:=\varepsilon\int_{\mathbb{R}^3}\nabla^k\Big(\frac{1}{\rho+n}\mathcal{L}v\Big)\cdot\nabla^k\nabla n\,\mathrm{d}x, && \mathfrak{L}_4:=\varepsilon\int_{\mathbb{R}^3}\nabla^k\left[\Big(\frac{1}{\rho+n}-1\Big)\mathcal{L}u\right]\cdot\nabla^k\nabla n\,\mathrm{d}x,\\
&\mathfrak{L}_5:=-\int_{\mathbb{R}^3}\nabla^k(v\cdot\nabla v)\cdot\nabla^k\nabla n+\nabla^k v\cdot\nabla^k\nabla\mathrm{div}\,(nv)\,\mathrm{d}x,\\
&\mathfrak{L}_6:=-\int_{\mathbb{R}^3}\nabla^k\nabla\mathfrak{H}_1\cdot\nabla^k\nabla n\,\mathrm{d}x, && \mathfrak{L}_7:=-\int_{\mathbb{R}^3}\nabla^k\nabla\mathfrak{H}_2\cdot\nabla^k\nabla n\,\mathrm{d}x.
\end{aligned}
$$

By integration by parts, Hölder's inequality and Cauchy-Schwarz inequality, we have

$$
\begin{aligned}
|\mathfrak{L}_1|&=-\int_{\mathbb{R}^3}\nabla^k\mathrm{div}\,v\cdot\nabla^k\mathrm{div}\,(v+\varrho v+nu)\,\mathrm{d}x\lesssim\|\nabla^k\nabla v\|_{L^2}^2+\|(\varrho,u)\|_{W^{k+1,\infty}}\|(n,v)\|_{H^{k+1}}^2,\\
|\mathfrak{L}_2|&\lesssim\|u\|_{W^{k+1,\infty}}\|(n,v)\|_{H^{k+1}}^2,\\
|\mathfrak{L}_3|&\lesssim\varepsilon\left\|\nabla^k\left[\Big(\frac{1}{\rho+n}-1\Big)\mathcal{L}v\right]\right\|_{L^2}\|\nabla^k\nabla n\|_{L^2}+\varepsilon\|\nabla^{k+2}v\|_{L^2}\|\nabla^k\nabla n\|_{L^2}\\
&\lesssim\varepsilon\|\nabla n\|_{H^k}\left(\|\varrho+n\|_{L^\infty}\|\nabla^2 v\|_{H^k}+\|\varrho+n\|_{W^{k,6}}\|\nabla^2 v\|_{L^3}\right)+\varepsilon\|\nabla^{k+2}v\|_{L^2}\|\nabla^k\nabla n\|_{L^2}\\
&\lesssim\varepsilon\|(\varrho,n)\|_{L^\infty}\|(\nabla n,\nabla^2 v)\|_{H^k}^2+(\eta+\|\nabla v\|_{H^2})\|\nabla n\|_{H^k}^2+\varepsilon^2\|\nabla v\|_{H^2}^2\|\varrho\|_{W^{k,6}}^2+\|\nabla^{k+2}v\|_{L^2}^2,\\
|\mathfrak{L}_4|&\lesssim\varepsilon\left\|\nabla^k\left[\Big(\frac{1}{\rho+n}-1\Big)\mathcal{L}u\right]\right\|_{L^2}\|\nabla^k\nabla n\|_{L^2}\\
&\lesssim\varepsilon\|\varrho+n\|_{H^k}\|\nabla^2 u\|_{W^{k,\infty}}\|\nabla n\|_{H^k}\lesssim\|\nabla^2 u\|_{W^{k,\infty}}\|n\|_{H^{k+1}}^2+\varepsilon^2\|\nabla^2 u\|_{W^{k,\infty}}\|\varrho\|_{H^k}^2,\\
|\mathfrak{L}_5|&\lesssim\|(\nabla n,\nabla v)\|_{H^{\max\{k,1\}}}^3.
\end{aligned}
$$

Thus we derive from $k+1\le N_1$ that

$$
\begin{aligned}
\sum_{i=1}^{5}|\mathfrak{L}_i|\lesssim{}&\eta\|\nabla n\|_{H^{N_1-1}}^2+\|\nabla v\|_{H^{N_1}}^2+\left\{\varepsilon\sqrt{\mathbb{E}_N(t)}+\sqrt{\mathcal{E}_{N_1}(t)}\right\}\mathcal{D}_{N_1}(t)\\
&+\left\{\|(\varrho,u)\|_{W^{N_1+1,\infty}}+\|\varrho\|_{W^{N_1-1,6}}^2\right\}\left\{\mathcal{E}_{N_1}(t)+\varepsilon^2\mathbb{E}_N(t)\right\}.
\end{aligned}
$$

For the term $\mathfrak{L}_6$, we note that $\mathfrak{H}_1$ can be expressed via a Taylor expansion in integral form

$$
\mathfrak{H}_1=(\gamma-2)\int_0^1(1-\theta)(\rho+\theta n)^{\gamma-3}\,\mathrm{d}\theta\,n^2:=\Psi(\rho,n)n^2.
$$

This representation, combining with Lemma B.4 and (3.2) leads to

$$\begin{aligned}
|\mathfrak{L}_6| &\lesssim \|\nabla^k\nabla(\Psi(\rho,n)n^2)\|_{L^2}\|\nabla^k\nabla n\|_{L^2}\\
&\lesssim (\|\nabla^{k+1}\Psi(\rho,n)\|_{L^2}\|n^2\|_{L^\infty}+\|\nabla^k\nabla(n^2)\|_{L^2}\|\Psi(\rho,n)\|_{L^\infty})\|\nabla^k\nabla n\|_{L^2}\\
&\lesssim (\|(\varrho,n)\|_{H^{k+1}}\|n\|_{L^\infty}^2+\|n\|_{L^\infty}\|\nabla^k\nabla n\|_{L^2})\|\nabla^k\nabla n\|_{L^2}\\
&\lesssim \sqrt{\mathcal{E}_{N_1}(t)}\mathcal{D}_{N_1}(t).
\end{aligned}$$

Similarly, for the term $\mathfrak{L}_7$, we have

$$|\mathfrak{L}_7| \lesssim \|\varrho\|_{W^{k+1,\infty}}\mathcal{E}_{N_1}(t).$$

Collecting all the estimates of $\mathfrak{L}_1$–$\mathfrak{L}_7$, we thus conclude our proof. □

**Corollary 3.4.** *If* (3.2) *holds, then we have*

$$\begin{aligned}
&\frac{d}{dt}\bar{\mathcal{E}}_{N_1}(t)+C\varepsilon\kappa_3\mathcal{D}_{N_1}(t)\\
&\lesssim (\|(\varrho,u,|\nabla|^s\phi)\|_{W^{N_1+2,\infty}}+\|\varrho\|_{W^{N_1-1,6}}^2)\left\{\mathcal{E}_{N_1}(t)+\varepsilon^2\mathbb{E}_N(t)\right\}\\
&\quad+\left\{\varepsilon\sqrt{\mathbb{E}_N(t)}+\sqrt{\mathcal{E}_{N_1}(t)}\right\}\mathcal{D}_{N_1}(t).
\end{aligned}\tag{3.8}$$

*Here, the modified energy functional is defined as*

$$\bar{\mathcal{E}}_{N_1}(t):=\mathcal{E}_{N_1}(t)+C\varepsilon\kappa_3\sum_{k=0}^{N_1-1}\int_{\mathbb{R}^3}\nabla^k v\cdot\nabla^{k+1}n\,\mathrm{d}x,$$

*with the constant $\kappa_3>0$ chosen sufficiently small such that $\bar{\mathcal{E}}_{N_1}(t)\simeq\mathcal{E}_{N_1}(t)$.*

*Proof.* Summing the estimates in Proposition 3.2 over $1\le k\le N_1$ and combining the result with Proposition 3.1, we deduce that

$$\begin{aligned}
&\frac{1}{2}\frac{d}{dt}\mathcal{E}_{N_1}(t)+\varepsilon\|\nabla v\|_{H^{N_1}}^2\\
&\lesssim (\|(\varrho,u,|\nabla|^s\phi)\|_{W^{N_1+2,\infty}}+\|\varrho\|_{W^{N_1-1,6}}^2)\left\{\mathcal{E}_{N_1}(t)+\varepsilon^2\mathbb{E}_N(t)\right\}\\
&\quad+\left\{\varepsilon\eta+\varepsilon\sqrt{\mathbb{E}_N(t)}+\sqrt{\mathcal{E}_{N_1}(t)}\right\}\mathcal{D}_{N_1}(t).
\end{aligned}\tag{3.9}$$

By Proposition 3.3, we multiply (3.7) by $C\varepsilon\kappa_3$ with a sufficiently small constant $0<\kappa_3\ll 1$, and adding the result to (3.9), we arrive at (3.8) provided $\eta>0$ is chosen sufficiently small. This concludes the proof. □

**Proposition 3.5.** *If* (3.2) *holds, then for $0<\vartheta\le 1/2$ we have*

$$\begin{aligned}
\frac{d}{dt}\sqrt{\mathcal{E}^{-\vartheta}(t)} &\lesssim \|(\nabla\mathfrak{n},\nabla v,\nabla|\nabla|^s\psi)\|_{H^1}^2+\|(\varrho,u)\|_{W^{2,\frac{3}{\vartheta}}}\|(\mathfrak{n},v)\|_{H^2}+\varepsilon\|\nabla^2u\|_{L^{\frac{3}{\vartheta}}}\|\varrho\|_{L^2}\\
&\quad+\begin{cases}\|(n,v)\|_{L^2}\|(\varrho,u)\|_{L^{\frac{3}{\vartheta+s-1}}}, & \text{for } 1<\vartheta+s<3/2,\\ \|(n,v)\|_{H^1}\|(\varrho,u)\|_{W^{1,\infty}}, & \text{for } \vartheta+s\le 1,\end{cases}
\end{aligned}\tag{3.10}$$

*where the negative-order Sobolev energy functional is defined as*

$$\mathcal{E}^{-\vartheta}(t):=\||\nabla|^{-\vartheta}(\mathfrak{n},v,|\nabla|^s\psi)\|_{L^2}^2.$$

*Proof.* Applying $|\nabla|^{-\vartheta}$ to (3.1), taking the $L^2$ inner product with $|\nabla|^{-\vartheta}\mathfrak{V}$, and integrating by parts, we obtain

$$\frac{1}{2}\frac{d}{dt}\||\nabla|^{-\vartheta}\mathfrak{V}\|_{L^2}^2+\varepsilon\||\nabla|^{-\vartheta}\nabla v\|_{L^2}^2+\varepsilon\||\nabla|^{-\vartheta}\mathrm{div}\,v\|_{L^2}^2=\sum_{i=1}^4\mathbf{J}_i,$$

where

$$\mathbf{J}_1 := -\int_{\mathbb{R}^3} |\nabla|^{-\vartheta}\Big[(u+v)\cdot\nabla\mathfrak{n} + \frac{\gamma-1}{2}(\mathfrak{n}+\mathfrak{c})\mathrm{div}\, v + \frac{\gamma-1}{2}\mathfrak{n}\mathrm{div}\, u + v\cdot\nabla\mathfrak{c}\Big]|\nabla|^{-\vartheta}\mathfrak{n}\,\mathrm{d}x,$$
$$\mathbf{J}_2 := -\int_{\mathbb{R}^3} |\nabla|^{-\vartheta}\Big[(u+v)\cdot\nabla v + \frac{\gamma-1}{2}(\mathfrak{n}+\mathfrak{c})\nabla\mathfrak{n} + \frac{\gamma-1}{2}\mathfrak{n}\nabla\mathfrak{c} + v\cdot\nabla u\Big]\cdot|\nabla|^{-\vartheta}v\,\mathrm{d}x,$$
$$\mathbf{J}_3 := -\varepsilon\int_{\mathbb{R}^3} |\nabla|^{-\vartheta}\Big[\Big(\frac{1}{\rho+n}-1\Big)(\mathcal{L}v+\mathcal{L}u)\Big]\cdot|\nabla|^{-\vartheta}v\,\mathrm{d}x,$$
$$\mathbf{J}_4 := -\int_{\mathbb{R}^3} |\nabla|^{-\vartheta}\nabla\psi\cdot|\nabla|^{-\vartheta}v\,\mathrm{d}x.$$

For the terms $\mathbf{J}_1$ and $\mathbf{J}_2$, employing Hölder's inequality and the Hardy-Littlewood-Sobolev inequality (Lemma B.3) yields

$$|\mathbf{J}_1| + |\mathbf{J}_2| \lesssim \||\nabla|^{-\vartheta}(\mathfrak{n},v)\|_{L^2}\big(\|(\nabla\mathfrak{n},\nabla v)\|_{L^2}\|(\mathfrak{n},v)\|_{L^{\frac{3}{\vartheta}}} + \|(\mathfrak{c},u)\|_{W^{1,\frac{3}{\vartheta}}}\|(\mathfrak{n},v)\|_{H^1}\big).$$

Similarly, the term $\mathbf{J}_3$ is estimated as

$$|\mathbf{J}_3| \lesssim \varepsilon\||\nabla|^{-\vartheta}v\|_{L^2}\big(\|(\varrho,u)\|_{W^{2,\frac{3}{\vartheta}}}\|(\mathfrak{n},v)\|_{H^2} + \|\nabla^2 v\|_{L^2}\|\mathfrak{n}\|_{L^{\frac{3}{\vartheta}}} + \|\nabla^2 u\|_{L^{\frac{3}{\vartheta}}}\|\varrho\|_{L^2}\big).$$

To handle $\mathbf{J}_4$, we perform integration by parts and utilize $(1.10)_1$ and $(1.10)_3$ to deduce

$$\mathbf{J}_4 = -\frac{1}{2}\frac{d}{dt}\||\nabla|^{-\vartheta}|\nabla|^s\psi\|_{L^2}^2 \underbrace{-\int_{\mathbb{R}^3}|\nabla|^{-\vartheta}\psi\,|\nabla|^{-\vartheta}\mathrm{div}\,(\varrho v + nv + nu)\,\mathrm{d}x}_{\mathbf{J}_5}.$$

The remaining term $\mathbf{J}_5$ requires separate treatment based on the index $\vartheta+s$.

**Case 1:** $1<\vartheta+s<3/2$**.** By the Sobolev embedding $\dot{H}^{\frac{3}{2}-(\vartheta+s)}(\mathbb{R}^3)\hookrightarrow L^{\frac{3}{\vartheta+s}}(\mathbb{R}^3)$, we have

$$\|nv\|_{L^{\frac{1}{\frac{1}{2}+\frac{\vartheta+s-1}{3}}}} \lesssim \|v\|_{L^6}\|n\|_{L^{\frac{3}{\vartheta+s}}} \lesssim \|\nabla v\|_{L^2}\|n\|_{\dot{H}^{\frac{3}{2}-(\vartheta+s)}} \lesssim \|\nabla v\|_{L^2}\||\nabla|^s\psi\|_{\dot{H}^{\frac{3}{2}-\vartheta}}.$$

Consequently, Hölder's and Hardy-Littlewood-Sobolev inequality (Lemma B.3) imply

$$\begin{aligned}|\mathbf{J}_5| &\lesssim \||\nabla|^{-\vartheta}|\nabla|^s\psi\|_{L^2}\|(\varrho v, nv, nu)\|_{L^{\frac{1}{\frac{1}{2}+\frac{\vartheta+s-1}{3}}}}\\ &\lesssim \||\nabla|^{-\vartheta}|\nabla|^s\psi\|_{L^2}\big(\|(v,n)\|_{L^2}\|(\varrho,u)\|_{L^{\frac{3}{\vartheta+s-1}}} + \|\nabla v\|_{L^2}\||\nabla|^s\psi\|_{\dot{H}^{\frac{3}{2}-\vartheta}}\big).\end{aligned}$$

**Case 2:** $\vartheta+s\le 1$**.** Fractional Sobolev embeddings guarantee that

$$\||\nabla|^{-\vartheta+1-s}(\varrho v + nu)\|_{L^2} \lesssim \|(\varrho,u)\|_{W^{1,\infty}}\|(n,v)\|_{H^1},$$

and furthermore,

$$\begin{aligned}\||\nabla|^{-\vartheta+1-s}(nv)\|_{L^2} &\lesssim \||\nabla|^{-\vartheta+1-s}n\|_{L^3}\|v\|_{L^6} + \||\nabla|^{-\vartheta+1-s}v\|_{L^{\frac{6}{3-2(\vartheta+s)}}}\|n\|_{L^{\frac{3}{\vartheta+s}}}\\ &\lesssim \||\nabla|^s\psi\|_{\dot{H}^{\frac{3}{2}-\vartheta}}\|\nabla v\|_{L^2}.\end{aligned}$$

Using Hölder's inequality alongside Lemma B.4, we bound $\mathbf{J}_5$ by

$$|\mathbf{J}_5| \lesssim \||\nabla|^{-\vartheta}|\nabla|^s\psi\|_{L^2}\big(\|(\varrho,u)\|_{W^{1,\infty}}\|(n,v)\|_{H^1} + \||\nabla|^s\psi\|_{\dot{H}^{\frac{3}{2}-\vartheta}}\|\nabla v\|_{L^2}\big).$$

Combining these bounds, using the fact $\|f\|_{L^{\frac{3}{\vartheta}}} \lesssim \|\nabla f\|_{H^1}$ for $0<\vartheta\le\frac{1}{2}$, and dividing both sides by $\sqrt{\mathcal{E}^{-\vartheta}(t)}$, we establish the desired estimate (3.10). □

*Proof of Theorem 1.7.* By virtue of the *a priori* assumptions (2.44) and (3.2), and the estimate

$$\|\varrho\|_{W^{N_1-1,6}} \lesssim \|\varrho\|_{W^{N_1-1,p}}^{\frac{2p}{3(p-2)}}\|\varrho\|_{H^{N_1-1}}^{\frac{p-6}{3(p-2)}} \lesssim (1+t)^{-\frac{8}{9}}\delta_1,$$

it follows from (3.8) in Corollary 3.4 that

$$\frac{d}{dt}\bar{\mathcal{E}}_{N_1}(t) + C\varepsilon\kappa_3\mathcal{D}_{N_1}(t) \lesssim (1+t)^{-\beta_s}\delta_1\mathcal{E}_{N_1}(t) + (1+t)^{-\beta_s}\varepsilon^2\delta_1^3. \tag{3.11}$$

Applying Gronwall's inequality to (3.11) and using $\bar{\mathcal{E}}_{N_1}(t) \simeq \mathcal{E}_{N_1}(t)$ yields

$$\mathcal{E}_{N_1}(t) + C\varepsilon\kappa_3\int_0^t \mathcal{D}_{N_1}(\tau)\,\mathrm{d}\tau \lesssim \mathcal{E}_{N_1}(0) + \varepsilon^2\delta_1^3. \tag{3.12}$$

The global existence of $(n, |\nabla|^s\psi, v)$ in $C([0,+\infty), H^{N_1})$ then follows from (3.3) and a standard bootstrap argument.

Having established the global-in-time existence and uniform bounds, we now derive the $L^2$ decay rates using the negative Sobolev estimates from Proposition 3.5.

**Step 1.** Define

$$\widetilde{\mathcal{E}}(t) := \varepsilon^2\bar{\mathbb{E}}_N(t) + \bar{\mathcal{E}}_{N_1}(t), \quad \text{and} \quad \widetilde{\mathcal{D}}(t) := \varepsilon^2\mathbb{D}_N(t) + \mathcal{D}_{N_1}(t).$$

It follows from (3.8), (2.67), (2.44) and (3.2) that

$$\frac{d}{dt}\widetilde{\mathcal{E}}(t) + \varepsilon\widetilde{\mathcal{D}}(t) \lesssim (\|(\varrho, u, |\nabla|^s\phi)\|_{W^{N_1+2,\infty}} + \|\varrho\|^2_{W^{N_1-1,6}} + \mathbf{1}_{0<s<\frac{1}{2}}\|\nabla\varrho\|_{L^{\frac{3}{s}}})\,\widetilde{\mathcal{E}}(t).$$

Similar to the estimate of (2.85), we have

$$\widetilde{\mathcal{E}}(t) \lesssim \left[\left\{\widetilde{\mathcal{E}}(0)\right\}^{-\frac{1}{\vartheta}} + \varepsilon t\left\{\sup_{t\geq 0}\widetilde{\mathcal{E}}(t) + \varepsilon^2\sup_{t\geq 0}\mathbb{E}^{-\vartheta}(t) + \sup_{t\geq 0}\mathcal{E}^{-\vartheta}(t)\right\}^{-\frac{1}{\vartheta}}\right]^{-\vartheta}. \tag{3.13}$$

**Step 2.** We derive the decay rates by dividing the proof into two cases: $0 < s \leq 1/2$ and $1/2 < s < 1$.

**Case 1:** $0 < s \leq 1/2$**.** In light of Proposition 3.5, the energy estimate (3.12) and the decay estimate $\|(\varrho, u)\|_{W^{2,\frac{3}{\vartheta}}} \lesssim (1+t)^{-\frac{3}{2}(1-\frac{2\vartheta}{3})}$ for $0 < \vartheta < \frac{1}{2}$, we deduce that

$$\mathcal{E}^{-\vartheta}(t) \lesssim \mathcal{E}^{-\vartheta}(0) + \varepsilon^2\delta_1^4 \lesssim \varepsilon^2. \tag{3.14}$$

Together with (3.13), (2.44) and (3.2), this yields

$$\widetilde{\mathcal{E}}(t) \lesssim \varepsilon^2(1+\varepsilon t)^{-\vartheta}.$$

Next, for $\vartheta = \frac{1}{2}$, utilizing (3.14) and the embedding $\dot{H}^{-\vartheta} \cap L^2 \subset \dot{H}^{-\vartheta'}$ for any fixed $\vartheta' \in (0, \frac{1}{2})$, we have

$$\|(\varrho, u)\|_{W^{2,\frac{3}{\vartheta}}}\|(\mathfrak{n}, v)\|_{H^2} + \varepsilon\|\nabla^2 u\|_{L^{\frac{3}{\vartheta}}}\|\varrho\|_{L^2} \lesssim \varepsilon(1+\varepsilon t)^{-\frac{\vartheta'}{2}}(1+t)^{-1}.$$

Combining this with (3.10) and applying Lemma B.16, we obtain $\mathcal{E}^{-\vartheta}(t) \lesssim \varepsilon^2\left[\ln(1+\varepsilon^{-1})\right]^2$. Inserting this into (3.13), we can derive

$$\widetilde{\mathcal{E}}(t) \lesssim \varepsilon^2(1+\varepsilon^{1^+}t)^{-\vartheta}.$$

**Case 2:** $1/2 < s < 1$**.** The slower algebraic decay shifts the critical integrability thresholds. We first consider $0 < \vartheta < \frac{3}{8}$. Following the identical argument for $0 < \vartheta < \frac{1}{2}$ in **Case 1**, we deduce

$$\widetilde{\mathcal{E}}(t) \lesssim \varepsilon^2(1+\varepsilon t)^{-\vartheta}. \tag{3.15}$$

Next, for $\frac{3}{8} \leq \vartheta \leq \frac{1}{2}$, utilizing (3.15) and the embedding $\dot{H}^{-\vartheta} \cap L^2 \subset \dot{H}^{-\frac{1}{4}}$, we have

$$\begin{aligned}&\|(\varrho, u)\|_{W^{2,\frac{3}{\vartheta}}}\|(\mathfrak{n}, v)\|_{H^2} + \varepsilon\|\nabla^2 u\|_{L^{\frac{3}{\vartheta}}}\|\varrho\|_{L^2} + \|(n, v)\|_{L^2}\|(\varrho, u)\|_{L^{\frac{3}{\vartheta+s-1}}} + \|(n, v)\|_{H^1}\|(\varrho, u)\|_{W^{1,\infty}}\\ &\lesssim \varepsilon(1+\varepsilon t)^{-\frac{1}{8}}(1+t)^{-\frac{4}{3}\left(1-\frac{2\vartheta}{3}\right)} + \varepsilon(1+\varepsilon t)^{-\frac{1}{8}}(1+t)^{-\left(\frac{20}{9}-\frac{8}{9}(\vartheta+s)\right)} + \varepsilon(1+t)^{-\frac{4}{3}(1-\frac{2}{p})}.\end{aligned}$$

The temporal integration is dominated by the first term with decay exponent $\beta = \frac{4}{3} - \frac{8\vartheta}{9}$. Combining this with (3.10) and applying Lemma B.16, we obtain

$$\mathcal{E}^{-\vartheta}(t) \lesssim \begin{cases} \varepsilon^2\big[\ln(1+\varepsilon^{-1})\big]^2, & \text{if } \vartheta = \frac{3}{8}, \\ \varepsilon^{2-\frac{16\vartheta-6}{9}}, & \text{if } \frac{3}{8} < \vartheta \le \frac{1}{2}. \end{cases}$$

Inserting this into (3.13) leads to

$$\widetilde{\mathcal{E}}(t) \lesssim \varepsilon^2\big(1+\varepsilon^{\Gamma_4(\vartheta)}t\big)^{-\vartheta},$$

where the effective exponent $\Gamma_4 \ge 1$ is defined by

$$\Gamma_4(\vartheta) := \begin{cases} 1^+, & \text{if } \vartheta = \frac{3}{8}, \\ 1 + \frac{16\vartheta-6}{9\vartheta}, & \text{if } \frac{3}{8} < \vartheta \le \frac{1}{2}. \end{cases}$$

Summarizing the two cases, the decay estimates for the combined energy functional can be unified as

$$\widetilde{\mathcal{E}}(t) \lesssim \varepsilon^2\left(1+\varepsilon^{\Upsilon_3}t\right)^{-\vartheta} \lesssim \min\left\{\varepsilon^2, (1+t)^{-\vartheta}\right\},$$

where the comprehensive effective exponent $\Upsilon_3 \ge 1$ is defined by

$$\Upsilon_3 := \begin{cases} 1, & \text{if } 0 < s \le 1/2 \text{ with } 0 < \vartheta < \frac{1}{2}, \text{ or } 1/2 < s < 1 \text{ with } 0 < \vartheta < \frac{3}{8}, \\ 1^+, & \text{if } 0 < s \le 1/2 \text{ with } \vartheta = \frac{1}{2}, \text{ or } 1/2 < s < 1 \text{ with } \vartheta = \frac{3}{8}, \\ 1 + \frac{16\vartheta-6}{9\vartheta}, & \text{if } 1/2 < s < 1 \text{ with } \frac{3}{8} < \vartheta \le \frac{1}{2}. \end{cases}$$

Thus, we conclude the proof of Theorem 1.7. □

## 4. Inviscid Limit to the Compressible Euler-Riesz System

In this section, we investigate the inviscid limit from the compressible Navier-Stokes-Riesz system to the compressible Euler-Riesz system as $\varepsilon \to 0$. By estimating the difference between the two solutions, we not only rigorously justify the limit but also provide an explicit convergence rate for all time $t \ge 0$.

We recall that the solution in Theorem 1.9 is constructed through

$$(\rho^\varepsilon, u^\varepsilon) = (\rho, u) + (n, v),$$

where $(\rho, u)$ and $(n, v)$ are obtained by Theorem 1.3 and Theorem 1.7, respectively.

First, we have, following directly from Theorem 1.7,

**Corollary 4.1.** *Under the assumption of Theorem 1.9, for* $3 \le N_1 \le N_0 - 3$*, we have*

$$\sup_{t\ge0} \|(\rho^\varepsilon - \rho, u^\varepsilon - u, |\nabla|^s(\phi^\varepsilon - \phi))(t)\|^2_{H^{N_1}} \le C\varepsilon^2,$$

*where* $C > 0$ *is a constant independent of* $t$ *and* $\varepsilon$*.*

To prove Theorem 1.11, it suffices to prove the following result.

**Corollary 4.2.** *Under the assumption of Theorem 1.3, let* $\varepsilon > 0$ *be sufficiently small and* $(\widetilde{\rho}, \widetilde{u})$ *be the global solution to the compressible Euler-Riesz system in Lemma C.1 with the initial data* $(\widetilde{\rho}_0, \widetilde{u}_0) := (\rho_0^\varepsilon, \mathcal{P}^\perp u_0^\varepsilon)$*, then for any integer* $2 \le N_0' \le N_0 - 2$ *and* $N_2 \le N_0' - 3/2(1-2/p)$*, we have*

$$\sup_{0\le t<+\infty} \|(\rho - \widetilde{\rho}, u - \widetilde{u}, |\nabla|^s(\phi - \widetilde{\phi}))(t)\|_{W^{N_2,p}} \le C\varepsilon^{\lambda_s}.$$

*where* $C$ *is a constant independent of* $t$ *and* $\varepsilon$ *and* $\lambda_s$ *is defined by* (1.11)*.*

Analogous to the transformation (2.45), we introduce

$$\widetilde{\mathfrak{c}} := \frac{2}{\gamma - 1}\left(\widetilde{\rho}^{\frac{\gamma-1}{2}} - 1\right).$$

Consequently, $\widetilde{\mathfrak{U}} := (\widetilde{\mathfrak{c}}, \widetilde{u})$ satisfies

$$\partial_t \widetilde{\mathfrak{U}} + \sum_{i=1}^{3} \mathbf{M}_i(\widetilde{\mathfrak{U}})\partial_{x_i}\widetilde{\mathfrak{U}} + (0, \nabla\widetilde{\phi})^{\mathsf{T}} = 0. \tag{4.1}$$

We define

$$\mathfrak{W} := (\mathfrak{s}, w) := (\mathfrak{c} - \widetilde{\mathfrak{c}}, u - \widetilde{u}), \quad \varpi := \phi - \widetilde{\phi}.$$

It then follows from (4.1) and (2.46) that

$$\partial_t \mathfrak{W} + \sum_{i=1}^{3} \mathbf{M}_i(\mathfrak{U})\partial_{x_i}\mathfrak{W} - 2\varepsilon \mathbf{N}\Delta\mathfrak{W} + (0, \nabla\varpi)^{\mathsf{T}} = \mathfrak{S}, \tag{4.2}$$

where $\mathbf{M}_i$ and $\mathbf{N}$ are defined in (2.47), and

$$\mathfrak{S} := \begin{pmatrix} -w\cdot\nabla\widetilde{\mathfrak{c}} - \frac{\gamma-1}{2}\mathfrak{s}\operatorname{div}\widetilde{u} \\ -w\cdot\nabla\widetilde{u} - \frac{\gamma-1}{2}\mathfrak{s}\nabla\widetilde{\mathfrak{c}} + 2\varepsilon\Delta\widetilde{u} \end{pmatrix}.$$

We then introduce the error energy functionals. For $k = 0$, the zero-order functional is given by

$$\mathscr{E}^0(t) := \|(\mathfrak{s}, w, |\nabla|^s\varpi)\|_{L^2}^2.$$

For $k \geq 1$, the higher-order energy functionals depend on the fractional exponent $s$

$$\mathscr{E}^k(t) := \begin{cases} \left\|(\nabla^k\mathfrak{s}, \nabla^k w, \nabla^k|\nabla|^s\varpi)\right\|_{L^2}^2 & \text{if } 1/2 \leq s < 1, \\ \left\|\left(\nabla^k\mathfrak{s}, \nabla^k w, \frac{1}{\sqrt{\widetilde{\rho}}}\nabla^k|\nabla|^s\varpi\right)\right\|_{L^2}^2 & \text{if } 0 < s < 1/2. \end{cases}$$

Consequently, for an integer $m \geq 0$, the total error energy functional is defined as

$$\mathscr{E}_m(t) := \sum_{k=0}^{m} \mathscr{E}^k(t).$$

For the subsequent analysis, we impose the following *a priori* assumption: there exist a constant $0 < \delta_2 \ll 1$ and an integer $N_0'$ with $2 \leq N_0' \leq N_0 - 2$ such that

$$\sup_{0\leq t\leq T_\varepsilon} \mathscr{E}_{N_0'}(t) \leq \delta_2^2, \tag{4.3}$$

which implies that

$$\|\mathfrak{s}\|_{H^{N_0'}} \simeq \|\varrho - \widetilde{\varrho}\|_{H^{N_0'}}. \tag{4.4}$$

We shall use the smallness condition (4.3) as a bootstrap hypothesis. More precisely, for a fixed $T_\varepsilon > 0$, we assume that a smooth solution $(\mathfrak{s}, w, \varpi)$ on $[0, T_\varepsilon]$ satisfies (2.44). The estimates derived below are uniform with respect to $T_\varepsilon$. In addition, $T_\varepsilon$ will be chosen specifically later.

*Proof of Corollary 4.2.* Applying $\nabla^k$ to (4.2), and taking the $L^2$ inner products with $\nabla^k\mathfrak{W}$, we obtain

$$\frac{1}{2}\frac{d}{dt}\left(\|\nabla^k\mathfrak{s}\|_{L^2}^2 + \|\nabla^k w\|_{L^2}^2\right) + 2\varepsilon\|\nabla^{k+1}w\|_{L^2}^2 = \sum_{i=1}^{7} Q_i,$$

where the terms $Q_1$–$Q_7$ are defined as follows

$$
\begin{aligned}
Q_1 &:= -\int_{\mathbb{R}^3} u\cdot\nabla\nabla^k\mathfrak{s}\,\nabla^k\mathfrak{s}\,\mathrm{d}x - \int_{\mathbb{R}^3} u\cdot\nabla\nabla^k w\cdot\nabla^k w\,\mathrm{d}x,\\
Q_2 &:= -\int_{\mathbb{R}^3}\Big[\nabla^k, u\cdot\nabla\Big]\,\mathfrak{s}\,\nabla^k\mathfrak{s}\,\mathrm{d}x - \int_{\mathbb{R}^3}\Big[\nabla^k, u\cdot\nabla\Big]\,w\cdot\nabla^k w\,\mathrm{d}x,\\
Q_3 &:= -\int_{\mathbb{R}^3}\Big(\frac{\gamma-1}{2}\mathfrak{c}+1\Big)\,\mathrm{div}\,\nabla^k w\,\nabla^k\mathfrak{s} + \Big(\frac{\gamma-1}{2}\mathfrak{c}+1\Big)\,\nabla\nabla^k\mathfrak{s}\cdot\nabla^k w\,\mathrm{d}x,\\
Q_4 &:= -\frac{\gamma-1}{2}\int_{\mathbb{R}^3}\Big[\nabla^k,\mathfrak{c}\Big]\,\mathrm{div}\,w\,\nabla^k\mathfrak{s} + \Big[\nabla^k,\mathfrak{c}\Big]\,\nabla\mathfrak{s}\cdot\nabla^k w\,\mathrm{d}x,\\
Q_5 &:= -\int_{\mathbb{R}^3}\nabla^k\left(w\cdot\nabla\widetilde{\mathfrak{c}}\right)\nabla^k\mathfrak{s} + \nabla^k\left(w\cdot\nabla\widetilde{u}\right)\cdot\nabla^k w\,\mathrm{d}x\\
&\quad - \frac{\gamma-1}{2}\int_{\mathbb{R}^3}\nabla^k\left(\mathfrak{s}\mathrm{div}\,\widetilde{u}\right)\nabla^k\mathfrak{s} + \nabla^k\left(\mathfrak{s}\nabla\widetilde{\mathfrak{c}}\right)\cdot\nabla^k w\,\mathrm{d}x,\\
Q_6 &:= 2\varepsilon\int_{\mathbb{R}^3}\nabla^k\Delta\widetilde{u}\cdot\nabla^k w\,\mathrm{d}x, \qquad Q_7 := -\int_{\mathbb{R}^3}\nabla^k\nabla\varpi\cdot\nabla^k w\,\mathrm{d}x.
\end{aligned}
$$

By integration by parts, Hölder's inequality, and Lemma B.4, the first five terms can be bounded as

$$
\begin{aligned}
|Q_1|+|Q_3| &\lesssim \|(\nabla u,\nabla\mathfrak{c})\|_{L^\infty}\|(\nabla^k\mathfrak{s},\nabla^k w)\|_{L^2}^2,\\
|Q_2|+|Q_4| &\lesssim \Big(\|(\nabla^k u,\nabla^k\mathfrak{c})\|_{L^\infty}\|(\nabla\mathfrak{s},\nabla w)\|_{L^2} + \|(\nabla u,\nabla\mathfrak{c})\|_{L^\infty}\|(\nabla^k\mathfrak{s},\nabla^k w)\|_{L^2}\Big)\|(\nabla^k\mathfrak{s},\nabla^k w)\|_{L^2},\\
|Q_5| &\lesssim \|(\nabla\widetilde{\mathfrak{c}},\nabla\widetilde{u})\|_{W^{k,\infty}}\|(\mathfrak{s},w)\|_{H^k}^2.
\end{aligned}
$$

Since $k\le N_0'\le N_0-2$, we essentially have $k+2\le N_0<N$. Therefore, for $Q_6$, applying Hölder's inequality and Lemma C.1 yields

$$
|Q_6|\lesssim\varepsilon\|\nabla^{k+2}\widetilde{u}\|_{L^2}\|\nabla^k w\|_{L^2}\lesssim\varepsilon\|\widetilde{u}\|_{H^N}\sqrt{\mathscr{E}_{N_0'}(t)}\lesssim\varepsilon\delta_0\sqrt{\mathscr{E}_{N_0'}(t)}.
$$

To estimate $Q_7$, we notice that the continuity equations for $\rho$ and $\widetilde{\rho}$ yield

$$
\partial_t(\varrho-\widetilde{\varrho})+\mathrm{div}\,w+\mathrm{div}\,(\widetilde{\varrho}w)+\mathrm{div}\,\big((\varrho-\widetilde{\varrho})u\big)=0, \tag{4.5}
$$

where we denote $\widetilde{\varrho}:=\widetilde{\rho}-1$. Then, we treat two cases similar as before based on the range of $s$.

**Case 1:** $1/2\le s<1$. By integration by parts, (4.5) and $(-\Delta)^s\varpi=\varrho-\widetilde{\varrho}$, we obtain

$$
Q_7=-\frac{1}{2}\frac{d}{dt}\|\nabla^k|\nabla|^s\varpi\|_{L^2}^2+\int_{\mathbb{R}^3}\nabla^k\nabla\varpi\cdot\nabla^k\big(\widetilde{\varrho}w+(\varrho-\widetilde{\varrho})u\big)\,dx. \tag{4.6}
$$

For the second term on the right-hand side, utilizing the Sobolev interpolation and Lemma B.4, we deduce that

$$
\begin{aligned}
&\left|\int_{\mathbb{R}^3}\nabla^k\nabla\varpi\cdot\nabla^k\left(\widetilde{\varrho}w+(\varrho-\widetilde{\varrho})u\right)\,\mathrm{d}x\right|\\
&\lesssim\|\nabla^k|\nabla|^s\varpi\|_{L^2}^{\frac{2s-1}{s}}\|\nabla^k(-\Delta)^s\varpi\|_{L^2}^{\frac{1-s}{s}}\left\|\nabla^k\left(\widetilde{\varrho}w+(\varrho-\widetilde{\varrho})u\right)\right\|_{L^2}\\
&\lesssim\|\nabla^k|\nabla|^s\varpi\|_{L^2}^{\frac{2s-1}{s}}\|\nabla^k(\varrho-\widetilde{\varrho})\|_{L^2}^{\frac{1-s}{s}}\|(\widetilde{\varrho},u)\|_{W^{k,\infty}}\|(w,\varrho-\widetilde{\varrho})\|_{H^k}\\
&\lesssim\|(\widetilde{\varrho},u)\|_{W^{k,\infty}}\|(\varrho-\widetilde{\varrho},|\nabla|^s\varpi,w)\|_{H^k}^2.
\end{aligned}
$$

**Case 2:** $0 < s < 1/2$. For $k = 0$, we employ (4.6) to obtain

$$\left|\int_{\mathbb{R}^3} \nabla\varpi \cdot (\widetilde{\varrho} w + (\varrho - \widetilde{\varrho})u)\, \mathrm{d}x\right| \lesssim \|(-\Delta)^s \varpi\|_{L^2} \left(\||\nabla|^{1-2s}(\widetilde{\varrho} w)\|_{L^2} + \||\nabla|^{1-2s}[(\varrho - \widetilde{\varrho})u]\|_{L^2}\right)$$
$$\lesssim \|(\widetilde{\varrho}, u)\|_{W^{1,\infty}} \|(\varrho - \widetilde{\varrho}, w)\|_{H^1}^2.$$

For $k \geq 1$, invoking $(4.2)_3$ and (4.5), we reformulate $Q_7$ as

$$Q_7 = -\frac{1}{2}\frac{d}{dt}\int_{\mathbb{R}^3} \frac{1}{\widetilde{\rho}} |\nabla^k |\nabla|^s \varpi|^2 \, \mathrm{d}x + \sum_{j=1}^{6} Q_{7j},$$

where

$$Q_{71} := -\int_{\mathbb{R}^3} \frac{1}{\widetilde{\rho}} \left[\nabla^k |\nabla|^{-s}, \widetilde{\rho}\right] \operatorname{div} w \, \nabla^k |\nabla|^s \varpi \, \mathrm{d}x,$$
$$Q_{72} := -\frac{1}{2}\int_{\mathbb{R}^3} \frac{1}{\widetilde{\rho}^2} \partial_t \widetilde{\rho} |\nabla^k |\nabla|^s \varpi|^2 \, \mathrm{d}x,$$
$$Q_{73} := -\int_{\mathbb{R}^3} \frac{1}{\widetilde{\rho}} \nabla^k |\nabla|^{-s} \left(w \cdot \nabla\widetilde{\rho}\right) \nabla^k |\nabla|^s \varpi \, \mathrm{d}x,$$
$$Q_{74} := -\int_{\mathbb{R}^3} \frac{1}{\widetilde{\rho}} \nabla^k |\nabla|^{-s} \left((\varrho - \widetilde{\varrho}) \operatorname{div} u\right) \nabla^k |\nabla|^s \varpi \, \mathrm{d}x,$$
$$Q_{75} := -\int_{\mathbb{R}^3} \frac{1}{\widetilde{\rho}} \left[\nabla^k |\nabla|^{-s}, u\right] \cdot \nabla(\varrho - \widetilde{\varrho}) \, \nabla^k |\nabla|^s \varpi \, \mathrm{d}x,$$
$$Q_{76} := -\frac{1}{2}\int_{\mathbb{R}^3} \frac{1}{\widetilde{\rho}} u \cdot \nabla |\nabla^k |\nabla|^s \varpi|^2 \, \mathrm{d}x.$$

Rather than estimating these terms separately, we group them in a systematic way. With the aid of Lemma B.4 and $(C.1)_1$, the terms $Q_{72}$ and $Q_{76}$ can be estimated directly as

$$|Q_{72}| + |Q_{76}| \lesssim \|(\nabla\widetilde{\varrho}, \nabla u, \nabla\widetilde{u})\|_{L^\infty} \|\nabla^k |\nabla|^s \varpi\|_{L^2}^2.$$

Meanwhile, applying Lemma B.4 and $(C.1)_1$ to the remaining terms yields

$$\sum_{j \in \{1,3,4,5\}} |Q_{7j}| \lesssim \|(\nabla\widetilde{\varrho}, \nabla u, \nabla |\nabla|^s \widetilde{\phi})\|_{W^{k,\infty}} \|(\varrho - \widetilde{\varrho}, w, |\nabla|^s \varpi)\|_{H^k}^2.$$

Collecting the estimates for $Q_1$–$Q_7$ and taking the summation over $k \leq N_0'$, we arrive at

$$\frac{1}{2}\frac{d}{dt}\mathscr{E}_{N_0'}(t) + 2\varepsilon \|\nabla w\|_{H^{N_0'}}^2 \lesssim \|(\widetilde{\varrho}, \mathfrak{c}, \widetilde{\mathfrak{c}}, |\nabla|^s \widetilde{\phi}, u, \widetilde{u})\|_{W^{N_0'+1,\infty}} \mathscr{E}_{N_0'}(t) + \varepsilon\delta_0 \sqrt{\mathscr{E}_{N_0'}(t)}.$$

Here, we utilized the norm equivalence stated in (4.4). Furthermore, by Theorem 1.3, Lemma C.1, and the Sobolev embedding, the above differential inequality simplifies to

$$\frac{1}{2}\frac{d}{dt}\mathscr{E}_{N_0'}(t) + 2\varepsilon \|\nabla w\|_{H^{N_0'}}^2 \lesssim \delta_0 (1+t)^{-\beta_s} \mathscr{E}_{N_0'}(t) + \varepsilon\delta_0 \sqrt{\mathscr{E}_{N_0'}(t)}.$$

Dropping the non-negative dissipation term $2\varepsilon \|\nabla w\|_{H^{N_0'}}^2$ and dividing both sides by $\sqrt{\mathscr{E}_{N_0'}(t)}$, we obtain

$$\frac{d}{dt}\sqrt{\mathscr{E}_{N_0'}(t)} \leq C\delta_0 (1+t)^{-\beta_s} \sqrt{\mathscr{E}_{N_0'}(t)} + C\varepsilon\delta_0.$$

Since the viscous and inviscid systems share identical initial data, we have $\mathscr{E}_{N_0'}(0) = 0$. We define

$$\theta_s := \begin{cases} \dfrac{2p}{5p-6} & \text{when } 0 < s \le 1/2, \\ \dfrac{3p}{7p-8} & \text{when } 1/2 < s < 1. \end{cases}$$

Thus, by applying Gronwall's inequality, (4.4), and choosing $\varepsilon > 0$ sufficiently small and $T_\varepsilon := \varepsilon^{-\theta_s}$, we have $\mathscr{E}_{N_0'}(t) \lesssim \varepsilon^{2\lambda_s} < \delta_2^2/2$ on the interval $[0, T_\varepsilon]$. Moreover, with the help of the Sobolev embedding $H^{N_0'}(\mathbb{R}^3) \hookrightarrow W^{N_2,p}(\mathbb{R}^3)$ with $N_2 \le N_0' - 3/2(1-2/p)$, we deduce that

$$\sup_{0 \le t \le T_\varepsilon} \|(\rho - \widetilde{\rho}, |\nabla|^s(\phi - \widetilde{\phi}), u - \widetilde{u})(t)\|_{W^{N_2,p}} \lesssim \varepsilon^{\lambda_s},$$

where $\lambda_s$ is defined by (1.11). In addition, by Theorem 1.3 and Lemma C.1, we also have

$$\sup_{t > T_\varepsilon} \|(\rho - \widetilde{\rho}, |\nabla|^s(\phi - \widetilde{\phi}), u - \widetilde{u})(t)\|_{W^{N_2,p}} \lesssim (1+t)^{-\beta_s} \lesssim \varepsilon^{\lambda_s}.$$

Therefore, we conclude the proof of Corollary 4.2. □

## Appendix A. Properties of phase functions

In this section, we first study the following phase function,

$$b(r) := \sqrt{r^2 + r^{2-2s} - \varepsilon^2 r^4}. \tag{A.1}$$

**Lemma A.1.** *Let $0 < s < 1$ and $0 < \varepsilon \le 1$. There exists a sufficiently small constant $\kappa_0 > 0$ such that the following estimates hold uniformly for all $r > 0$ satisfying $\varepsilon r^2 < \kappa_0$. Then $b(r) > 0$ and*

$$b'(r) = \frac{r^{2s} + 1 - s - 2\varepsilon^2 r^{2+2s}}{r^s\big(r^{2s} + 1 - \varepsilon^2 r^{2+2s}\big)^{1/2}} > 0, \tag{A.2}$$

$$b''(r) = \frac{(2s^2 - s)r^{2s} - s + s^2 - \varepsilon^2 r^{2+2s}\big(3r^{2s} + 2s^2 + s + 3\big) + 2\varepsilon^4 r^{4+4s}}{r^{1+s}\big(r^{2s} + 1 - \varepsilon^2 r^{2+2s}\big)^{3/2}}, \tag{A.3}$$

$$b'''(r) = \frac{(-4s^3 + s)r^{4s} - 2(s^3 - s)r^{2s} - s^3 + s - \varepsilon^2 r^{2+2s}\big(3r^{4s} + c_{1,s} r^{2s} + c_{2,s}\big) - c_{3,s}\varepsilon^4 r^{4+4s}}{r^{2+s}\big(r^{2s} + 1 - \varepsilon^2 r^{2+2s}\big)^{5/2}}, \tag{A.4}$$

*where*

$$c_{1,s} = -8s^3 - 6s^2 + 8s + 6, \quad c_{2,s} = -2s^3 + 3s^2 + 8s + 3, \quad c_{3,s} = 4s^3 + 6s^2 + 2s.$$

*Moreover,*

$$b(r) \simeq \begin{cases} r^{1-s}, & 0 < r \le 1, \\ r, & r \ge 1, \end{cases} \qquad b'(r) \simeq \begin{cases} r^{-s}, & 0 < r \le 1, \\ 1, & r \ge 1. \end{cases} \tag{A.5}$$

*For every integer $k \ge 0$, one has the symbol-type bounds*

$$\big|\partial_r^k b(r)\big| \lesssim \begin{cases} r^{1-s-k}, & 0 < r \le 1, \\ r, & r \ge 1,\ k = 0, \\ 1, & r \ge 1,\ k = 1, \\ r^{1-2s-k}, & r \ge 1,\ k \ge 2,\ s \ne 1/2, \\ r^{-1-k}, & r \ge 1,\ k \ge 2,\ s = 1/2. \end{cases} \tag{A.6}$$

*In addition, the second derivative has the following more precise structure.*

*(1) If $0 < s \le 1/2$, then $b''(r) < 0$. More precisely,*

$$-b''(r) \simeq \begin{cases} r^{-1-s}, & 0 < r \le 1, \\ r^{-1-2s}, & r \ge 1,\ 0 < s < 1/2, \\ r^{-3}, & r \ge 1,\ s = 1/2. \end{cases} \tag{A.7}$$

*(2) If $1/2 < s < 1$, then $b''$ has a unique zero provided $\varepsilon r_*^2 < \kappa_0$ where $r_* := (1-s)^{\frac{1}{2s}}(2s-1)^{-\frac{1}{2s}}$. We denote the zero by $r_0$. Then*

$$r_0 \simeq 1, \qquad |b'''(r_0)| \simeq 1, \tag{A.8}$$

*and $b''$ changes sign from negative to positive at $r_0$, namely*

$$b''(r) < 0 \quad \text{for } r < r_0, \qquad b''(r) > 0 \quad \text{for } r > r_0.$$

*If instead $\varepsilon r_*^2 \ge \kappa_0$, then*

$$b''(r) < 0, \qquad \varepsilon r^2 < \kappa_0,$$

*and hence $b''$ has no zero in this region.*

*Proof.* Since $\varepsilon r^2 < \kappa_0$ and $0 < \varepsilon \le 1$, by choosing $\kappa_0 > 0$ sufficiently small we have

$$1 + (1-\kappa_0) r^{2s} \le r^{2s} + 1 - \varepsilon^2 r^{2+2s} \le 1 + r^{2s}. \tag{A.9}$$

Hence

$$b(r) = r^{1-s}\big(r^{2s} + 1 - \varepsilon^2 r^{2+2s}\big)^{1/2} \simeq r^{1-s}(1 + r^{2s})^{1/2},$$

which gives the estimate for $b$ in (A.5). The formula (A.2) for $b'$ follows from direct differentiation. Moreover, its numerator satisfies

$$r^{2s} + 1 - s - 2\varepsilon^2 r^{2+2s} \ge (1-s) + (1-2\kappa_0) r^{2s} > 0,$$

provided $\kappa_0$ is sufficiently small. This proves $b'(r) > 0$, and (A.5) for $b'$ follows.

The formulas (A.3) and (A.4) are obtained by direct differentiation of (A.1).

We now prove (A.6). The estimates for $k = 0, 1, 2$ have been proved as above. For $k \ge 2$, differentiating the identity

$$b^2(r) = r^2 + r^{2-2s} - \varepsilon^2 r^4$$

$k$ times gives the exact relation

$$\partial_r^k b(r) = \frac{1}{2b(r)}\left[\partial_r^k\big(r^2 + r^{2-2s} - \varepsilon^2 r^4\big) - \sum_{i=1}^{k-1}\binom{k}{i}\partial_r^i b(r)\partial_r^{k-i} b(r)\right]. \tag{A.10}$$

For $0 < r \le 1$, one has

$$\left|\partial_r^k\big(r^2 + r^{2-2s} - \varepsilon^2 r^4\big)\right| \lesssim r^{2-2s-k}.$$

Assuming the estimate (A.6) holds for all derivatives of order strictly less than $k$, each product in the sum of (A.10) is bounded by $Cr^{2-2s-k}$. Since $b(r) \simeq r^{1-s}$, (A.10) yields

$$|\partial_r^k b(r)| \lesssim r^{1-s-k}.$$

Let $r \ge 1$. If $s \ne 1/2$, the estimate for $k = 2$ follows from (A.3). For $k \ge 3$,

$$\left|\partial_r^k\big(r^2 + r^{2-2s} - \varepsilon^2 r^4\big)\right| \lesssim r^{2-2s-k}.$$

By the induction hypothesis, the terms in the sum of (A.10) are also bounded by $Cr^{2-2s-k}$. Since $b(r) \simeq r$, we obtain

$$|\partial_r^k b(r)| \lesssim r^{1-2s-k}, \qquad r \ge 1, \quad s \ne \frac{1}{2}, \quad k \ge 2.$$

It remains to handle the special case $s = 1/2$. The estimate for $k = 2$ again follows from (A.3). For $k \geq 3$, the part $r^2 + r$ contributes no derivative of order $k \geq 3$, while using $\varepsilon r^2 < \kappa_0$,

$$\left|\partial_r^k(\varepsilon^2 r^4)\right| \lesssim r^{-k}.$$

The induction hypothesis gives the same bound for the product terms in (A.10). Dividing by $b(r) \simeq r$, we conclude

$$|\partial_r^k b(r)| \lesssim r^{-1-k}, \qquad r \geq 1, \quad s = \frac{1}{2}, \quad k \geq 2.$$

This proves (A.6).

We next prove the assertions on $b''$. We write the numerator of (A.3) as

$$N_2(r) = (2s^2 - s)r^{2s} - s + s^2 - \varepsilon^2 r^{2+2s}\big(3r^{2s} + 2s^2 + s + 3\big) + 2\varepsilon^4 r^{4+4s}.$$

When $0 < s \leq 1/2$, the principal part

$$PN_2 := (2s^2 - s)r^{2s} - s + s^2 = s(2s-1)r^{2s} - s(1-s) < 0.$$

while the other terms can be made smaller than a fixed fraction of the principal negative part by choosing $\kappa_0$ sufficiently small. Therefore $N_2(r) < 0$, hence $b''(r) < 0$. The quantitative estimate (A.7) then follows from (A.3) together with (A.9).

Finally, consider $1/2 < s < 1$. The principal part of $N_2$ is

$$PN_2 := s(2s-1)r^{2s} - s(1-s),$$

which has the unique zero $r_*$. The remaining terms in $N_2$ are a small perturbation in the region $\varepsilon r^2 < \kappa_0$. Hence, inside that neighborhood of $r_*$, the implicit function theorem applies, because for $\varepsilon = 0$,

$$b''(r_*) = 0, \quad b'''(r_*) \neq 0.$$

More precisely,

$$b'''(r_*) = \frac{2s^3(1-s)}{(2s-1)r_*^{2+s}(1+r_*^{2s})^{5/2}} \neq 0 \quad \text{when } \varepsilon = 0.$$

Thus $b''$ has a unique zero in the region $\varepsilon r^2 < \kappa_0$ if $\varepsilon r_*^2 < \kappa_0$. We denote the zero by $r_0$, then $r_0 \simeq 1$, $|b'''(r_0)| \simeq 1$, and $b''$ changes sign from negative to positive at $r_0$. It remains to consider the case $\varepsilon r_*^2 \geq \kappa_0$. For any $r > 0$ satisfying $\varepsilon r^2 < \kappa_0$, we have

$$r < \left(\frac{\kappa_0}{\varepsilon}\right)^{1/2} \leq r_*.$$

Hence the principal part of the numerator of $b''$ satisfies $PN_2 < 0$. Moreover, choosing $\kappa_0 > 0$ sufficiently small, the remaining part of $N_2$ is non-positive. Therefore $N_2(r) < 0$. Since the denominator in (A.3) is positive by (A.9), we obtain

$$b''(r) < 0, \qquad \varepsilon r^2 < \kappa_0.$$

Thus, when $\varepsilon r_*^2 \geq \kappa_0$, the function $b''$ has no zero in the region $\{\varepsilon r^2 < \kappa_0\}$. This proves (A.8) and completes the proof. □

For our purpose, we further study the phase function $\Phi(\xi, \eta) = b(|\xi|) - b(|\xi - \eta|) - b(|\eta|)$, which will be used in demonstrating that the kernel $\mathfrak{M}$ belongs to the function spaces $\mathcal{M}^{k,\infty}_{\xi,\eta}$ and $\mathcal{M}^{k,\infty}_{\eta,\xi}$.

**Lemma A.2.** *Let* $0 < s < 1$, $\mathbf{k} = \mathbf{m} + \mathbf{n} \in \mathbb{R}^3$, *and* $|\mathbf{n}| \leq \min\{|\mathbf{k}|, |\mathbf{m}|\}$. *Assume that there exists a sufficiently small constant* $\kappa_0 > 0$ *such that* $\max\{\varepsilon|\mathbf{k}|^2, \varepsilon|\mathbf{m}|^2, \varepsilon|\mathbf{n}|^2\} \leq \kappa_0$. *Then, we have*

$$|b(\mathbf{k}) - b(\mathbf{m}) - b(\mathbf{n})| \gtrsim \frac{|\mathbf{n}|}{1 + (|\mathbf{k}||\mathbf{n}|)^s} + |\mathbf{n}|(1 - \cos[\mathbf{k}, \mathbf{n}] + 1 - \cos[\mathbf{k}, \mathbf{m}]) \tag{A.11}$$

*for* $|\mathbf{m}| \geq 1$*. Moreover, for* $|\mathbf{n}| < 1$*, we obtain*

$$|b(\mathbf{k}) - b(\mathbf{m}) - b(\mathbf{n})| \gtrsim \frac{|\mathbf{n}|^{1-s}}{(1+|\mathbf{m}|^{2s})^{1/2}}. \tag{A.12}$$

*Here and below,* $[\cdot,\cdot]$ *represents the angle between two vectors.*

*Proof.* We focus exclusively on estimating (A.11), as it provides a sharper exponent in the denominator than that in [12, Lemma 4.1]. The proof for (A.12) is omitted, as it follows from arguments strictly analogous to those in [12, Lemma 4.1] due to the shared properties of $b(r)$ and $q(r)$.

First, note that $|\mathbf{k}| \leq |\mathbf{m}|$ implies $b(\mathbf{k}) \leq b(\mathbf{m})$. Since $b(r)$ is a non-negative and increasing function by Lemma A.1, it follows that

$$|b(\mathbf{k}) - b(\mathbf{m}) - b(\mathbf{n})| \geq b(\mathbf{n}) \gtrsim (|\mathbf{n}|^2 + |\mathbf{n}|^{2-2s})^{1/2},$$

which yields the desired inequalities.

We now consider the case $|\mathbf{k}| > |\mathbf{m}|$. Let

$$q(r) = \frac{b(r)}{r} = \left(\frac{1 + r^{2s} - \varepsilon^2 r^{2s+2}}{r^{2s}}\right)^{1/2}.$$

One can observe that $q(r)$ is a non-negative, decreasing function satisfying

$$q(r) > 1, \quad q'(r) = -\frac{s + \varepsilon^2 r^{2+2s}}{r^{1+s}\sqrt{1 + r^{2s} - \varepsilon^2 r^{2+2s}}}, \quad -q'(r) \simeq \frac{s}{r^{1+2s}} \text{ as } r \to \infty,$$

$$q''(r) = \frac{s(s+1) + s(2s+1)r^{2s} - (2s+1)(s+1)\varepsilon^2 r^{2s+2} - \varepsilon^2 r^{4s+2}}{r^{s+2}(1 + r^{2s} - \varepsilon^2 r^{2s+2})^{3/2}} > 0.$$

Utilizing the relation $|\mathbf{k}| = |\mathbf{m}|\cos[\mathbf{k},\mathbf{m}] + |\mathbf{n}|\cos[\mathbf{k},\mathbf{n}]$, we decompose the phase difference as

$$\begin{aligned} b(\mathbf{m}) + b(\mathbf{n}) - b(\mathbf{k}) &= (|\mathbf{m}| + |\mathbf{n}| - |\mathbf{k}|)q(\mathbf{k}) + |\mathbf{m}|(q(\mathbf{m}) - q(\mathbf{k})) + |\mathbf{n}|(q(\mathbf{n}) - q(\mathbf{k})) \\ &= \{|\mathbf{m}|(1 - \cos[\mathbf{k},\mathbf{m}]) + |\mathbf{n}|(1 - \cos[\mathbf{k},\mathbf{n}])\}\, q(\mathbf{k}) \\ &\quad + |\mathbf{m}|(q(\mathbf{m}) - q(\mathbf{k})) + |\mathbf{n}|(q(\mathbf{n}) - q(\mathbf{k})). \end{aligned} \tag{A.13}$$

Since $q(r)$ is decreasing, each term on the right-hand side is positive.

To prove (A.11), we first assume $\cos[\mathbf{k},\mathbf{n}] \leq 9/10$. Then (A.13) implies

$$\begin{aligned} b(\mathbf{m}) + b(\mathbf{n}) - b(\mathbf{k}) &= |\mathbf{m}|(q(\mathbf{m}) - \cos[\mathbf{k},\mathbf{m}]q(\mathbf{k})) + |\mathbf{n}|(q(\mathbf{n}) - \cos[\mathbf{k},\mathbf{n}]q(\mathbf{k})) \\ &\geq |\mathbf{n}|(q(\mathbf{n}) - \cos[\mathbf{k},\mathbf{n}]q(\mathbf{k})) \\ &\geq \frac{|\mathbf{n}|q(\mathbf{n})}{10} + \frac{9|\mathbf{n}|}{10}(q(\mathbf{n}) - q(\mathbf{k})) \\ &\geq \frac{b(\mathbf{n})}{10} \end{aligned}$$

due to the decreasing property of $q(r)$. This establishes (A.11) for the case $\cos[\mathbf{k},\mathbf{n}] \leq 9/10$.

Next, we note that $|\mathbf{n}| \leq |\mathbf{m}|$ implies $\cos[\mathbf{k},\mathbf{n}] \leq \cos[\mathbf{k},\mathbf{m}]$, which leads to $|\mathbf{k}| \geq 2|\mathbf{n}|\cos[\mathbf{k},\mathbf{n}]$. When $\cos[\mathbf{k},\mathbf{n}] > 9/10$, it follows that

$$|\mathbf{k}| \geq \frac{4|\mathbf{n}|}{3}.$$

For $1 < |\mathbf{k}| < \eta^{-1}$ for $\eta > 0$ small enough, we have

$$|\mathbf{k}| \geq \frac{|\mathbf{k}|}{4} + |\mathbf{n}| \geq \frac{1}{4} + |\mathbf{n}|.$$

Consequently, we have

$$q(\mathbf{n}) - q(\mathbf{k}) = -\int_{|\mathbf{n}|}^{|\mathbf{k}|} q'(r)\,\mathrm{d}r \geq -q'(\eta^{-1}) \int_{|\mathbf{n}|}^{|\mathbf{k}|} \mathrm{d}r \gtrsim \eta^{2s+1} \int_{|\mathbf{n}|}^{|\mathbf{n}|+\frac{1}{4}} \mathrm{d}r \gtrsim 1,$$

which yields the estimate

$$|\mathbf{n}|(q(\mathbf{n}) - q(\mathbf{k})) \gtrsim |\mathbf{n}| \gtrsim \frac{|\mathbf{n}|}{1 + |\mathbf{k}|^s |\mathbf{n}|^s}. \tag{A.14}$$

For the larger regime $|\mathbf{k}| \geq \eta^{-1}$, the monotonicity of $-q'(r) > 0$ implies the situation simplifies to $|\mathbf{n}| > \eta^{-1}$. We then have

$$q(\mathbf{n}) - q(\mathbf{k}) = -\int_{|\mathbf{n}|}^{|\mathbf{k}|} q'(r)\,\mathrm{d}r \simeq \int_{|\mathbf{n}|}^{|\mathbf{k}|} \frac{s}{r^{1+2s}}\,\mathrm{d}r = \frac{|\mathbf{k}|^{2s} - |\mathbf{n}|^{2s}}{2|\mathbf{k}|^{2s}|\mathbf{n}|^{2s}} \gtrsim \frac{|\mathbf{k}|^s |\mathbf{n}|^s}{1 + (|\mathbf{k}||\mathbf{n}|)^{2s}} \gtrsim \frac{1}{1 + (|\mathbf{k}||\mathbf{n}|)^s}$$

due to $|\mathbf{k}| \geq 4|\mathbf{n}|/3$ and $|\mathbf{k}| > 1$. This yields

$$|\mathbf{n}|(q(\mathbf{n}) - q(\mathbf{k})) \gtrsim \frac{|\mathbf{n}|}{1 + |\mathbf{k}|^s |\mathbf{n}|^s}. \tag{A.15}$$

By combining (A.14) and (A.15) with the representation (A.13), we conclude that (A.11) holds for $\cos[\mathbf{k}, \mathbf{n}] > 9/10$. □

To prove Proposition B.12, we require precise estimates for the first- and second-order derivatives of the phase function $\Phi$, which are explicitly stated in Lemma A.3 below. We omit the detailed proof of these estimates, as it follows directly from the analytical framework established in [12, Lemma 4.1]. Specifically, because our function $b(r)$ shares analogous properties with the function $p(r)$ introduced in that work, the desired bounds can be recovered simply by substituting their parameter $\sigma$ with $2s$.

**Lemma A.3.** *Let $0 < s < 1$. Assume that there exists a sufficiently small constant $\kappa_0 > 0$ such that $\max\{\varepsilon|\xi|^2, \varepsilon|\xi - \eta|^2, \varepsilon|\eta|^2\} \leq \kappa_0$. Then the following estimates hold:*

$$|\nabla_\xi \Phi| \lesssim \begin{cases} |\eta| + 2\sin\dfrac{[\xi, \xi-\eta]}{2}, & \min\{|\xi|, |\xi-\eta|\} \geq 1, \\ \dfrac{1}{\min\{|\xi|, |\xi-\eta|\}^s} + 2\max\left\{1, \dfrac{1}{|\xi-\eta|^s}\right\} \sin\dfrac{[\xi, \xi-\eta]}{2}, & \min\{|\xi|, |\xi-\eta|\} < 1, \end{cases}$$

$$|\nabla_\eta \Phi| \lesssim \begin{cases} |\xi| + 2\sin\dfrac{[\eta, \xi-\eta]}{2}, & \min\{|\eta|, |\xi-\eta|\} \geq 1, \\ \dfrac{1}{\min\{|\eta|, |\xi-\eta|\}^s} + 2\max\left\{1, \dfrac{1}{|\xi-\eta|^s}\right\} \sin\dfrac{[\eta, \xi-\eta]}{2}, & \min\{|\eta|, |\xi-\eta|\} < 1, \end{cases}$$

$$|\Delta_\xi \Phi| \lesssim \begin{cases} |\eta|, & \min\{|\xi|, |\xi-\eta|\} \geq 1, \\ \dfrac{1}{\min\{|\xi|, |\xi-\eta|\}^{1+s}}, & \min\{|\xi|, |\xi-\eta|\} < 1, \end{cases}$$

$$|\Delta_\eta \Phi| \lesssim \begin{cases} |\xi|, & \min\{|\eta|, |\xi-\eta|\} \geq 1, \\ \dfrac{1}{\min\{|\eta|, |\xi-\eta|\}^{1+s}}, & \min\{|\eta|, |\xi-\eta|\} < 1. \end{cases}$$

*Moreover, if $\min\{|\xi|, |\xi-\eta|\} < 1$, then*

$$|\nabla_\xi \Phi| \lesssim \frac{|\eta|}{\min\{|\xi|, |\xi-\eta|\}^{1+s}} + 2\max\left\{1, \frac{1}{|\xi-\eta|^s}\right\} \sin\frac{[\xi, \xi-\eta]}{2},$$

$$|\Delta_\xi \Phi| \lesssim \frac{|\eta|}{\min\{|\xi|, |\xi-\eta|\}^{2+s}} + \frac{|\eta|^{1+s}}{\min\{|\xi|, |\xi-\eta|\}^{2+2s}}.$$

*Similarly, if* $\min\{|\eta|,|\xi-\eta|\}<1$*, then*

$$|\nabla_\eta\Phi|\lesssim\frac{|\xi|}{\min\{|\eta|,|\xi-\eta|\}^{1+s}}+2\max\left\{1,\frac{1}{|\xi-\eta|^s}\right\}\sin\frac{[\eta,\xi-\eta]}{2},$$

$$|\Delta_\eta\Phi|\lesssim\frac{|\xi|}{\min\{|\eta|,|\xi-\eta|\}^{2+s}}+\frac{|\xi|^{1+s}}{\min\{|\eta|,|\xi-\eta|\}^{2+2s}}.$$

## Appendix B. Some useful inequalities

In this section, we discuss several analytic inequalities that will be used frequently throughout this paper.

We begin with the Van der Corput lemma from [32].

**Lemma B.1** ([32]). *Suppose* $\phi$ *is real-valued and smooth in* $(a,b)$*, and that* $|\frac{d^k\phi(x)}{dx^k}|\gtrsim 1$ *for all* $x\in(a,b)$*. Then*

$$\left|\int_a^b e^{\mathrm{i}\lambda\phi(x)}\psi(x)\,\mathrm{d}x\right|\lesssim\lambda^{-\frac{1}{k}}\left(|\psi(b)|+\int_a^b|\psi'(x)|\,\mathrm{d}x\right),$$

*holds when* $k\geq 2$*, or* $k=1$ *and* $\phi'(x)$ *is monotonic.*

We then recall the Bernstein inequalities from [1].

**Lemma B.2** ([1]). *Let* $j\in\mathbb{Z}$*,* $\alpha\geq 0$*, and* $1\leq q\leq p\leq\infty$*, then*

$$\left\|\varphi\left(\frac{D}{2^j}\right)f\right\|_{L^p}\lesssim 2^{3j(\frac{1}{q}-\frac{1}{p})}\left\|\varphi\left(\frac{D}{2^j}\right)f\right\|_{L^q},$$

$$\left\||\nabla|^{\pm\alpha}\varphi\left(\frac{D}{2^j}\right)f\right\|_{L^p}\simeq 2^{\pm j\alpha}\left\|\varphi(\frac{D}{2^j})f\right\|_{L^p}\lesssim 2^{\pm j\alpha}\|f\|_{L^p},$$

*where* $\varphi(D/2^j)$ *is defined by* (1.12) *and* $\varphi$ *is given by* (1.18)*.*

The Hardy-Littlewood-Sobolev inequality from [1] will be used to handle negative derivatives.

**Lemma B.3** ([1]). *Let* $0<\vartheta<3$*,* $1<p<q<\infty$*,* $1/q+\vartheta/3=1/p$*, it holds that*

$$\||\nabla|^{-\vartheta}f\|_{L^q}\lesssim\|f\|_{L^p}.$$

For product and commutator estimates, we use the following Kato-Ponce inequalities.

**Lemma B.4** ([23, 25]). *Let* $s>0$ *and* $1<p<\infty$*. Assume that* $q_j,r_j\in(1,\infty]$ *for* $j=1,2$ *such that*

$$\frac{1}{p}=\frac{1}{q_1}+\frac{1}{r_1}=\frac{1}{q_2}+\frac{1}{r_2}.$$

*Then, if* $f\in W^{s,q_1}\cap L^{r_2}$ *and* $g\in W^{s,q_2}\cap L^{r_1}$*, we have*

$$\|fg\|_{\dot{W}^{s,p}}\lesssim\|f\|_{\dot{W}^{s,q_1}}\|g\|_{L^{r_1}}+\|f\|_{L^{r_2}}\|g\|_{\dot{W}^{s,q_2}},$$

$$\|fg\|_{W^{s,p}}\lesssim\|f\|_{W^{s,q_1}}\|g\|_{L^{r_1}}+\|f\|_{L^{r_2}}\|g\|_{W^{s,q_2}}.$$

*Furthermore, if* $f\in\dot{W}^{s,q_1}\cap\dot{W}^{1,r_2}$ *and* $g\in\dot{W}^{s-1,q_2}\cap L^{r_1}$*, we have*

$$\|[f,|\nabla|^s]\,g\|_{L^p}\lesssim\|f\|_{\dot{W}^{s,q_1}}\|g\|_{L^{r_1}}+\|\nabla f\|_{L^{r_2}}\|g\|_{\dot{W}^{s-1,q_2}}.$$

We use the following Moser-type composition estimate from [29].

**Lemma B.5** ([29]). *Suppose* $F\colon\mathbb{R}\to\mathbb{R}$ *is a smooth function with the condition* $F(0)=0$*. Then for any function* $u$ *that belongs to* $L^\infty\cap W^{k,p}$ *(*$k\geq 0$ *is an integer and* $1\leq p\leq+\infty$*), we have*

$$\|F(u)\|_{\dot{W}^{k,p}}\lesssim C(\|u\|_{L^\infty})\|u\|_{\dot{W}^{k,p}}.$$

*Furthermore, if* $k\geq 1$*, the condition* $F(0)=0$ *can be dropped.*

We now prove several multiplier bounds needed for the $W^{k,p}$ estimates.

Let us define the operators

$$n_1(D) := |\nabla| \quad \text{or} \quad \mathrm{i}b(D)\widetilde{\chi}_{\varepsilon,\kappa_0}(D)\mathcal{R}^* \quad \text{or} \quad \varepsilon\Delta\widetilde{\chi}_{\varepsilon,\kappa_0}(D)\mathcal{R}^*, \tag{B.1}$$

$$n_2(D) := \frac{|\nabla|}{a(D)} \quad \text{or} \quad \widetilde{\chi}_{\varepsilon,\kappa_0}(D)\mathcal{R}\frac{b(D)}{a(D)} \quad \text{or} \quad \widetilde{\chi}_{\varepsilon,\kappa_0}(D)\mathcal{R}\frac{-\varepsilon\Delta}{a(D)}, \tag{B.2}$$

$$n_3(D) := \frac{|\nabla|}{a(D)} \quad \text{or} \quad \mathcal{R}, \tag{B.3}$$

where $\mathcal{R}$ and $\mathcal{R}^*$ are the Riesz transform and its adjoint defined in (2.8), the differential operators $a(D)$ and $b(D)$ are defined in (2.2) and (2.11), respectively, and the differential operator $\chi_{\varepsilon,\kappa_0}$ is defined by (1.12), (1.16), and (1.15).

First, we recall the Hörmander–Mikhlin multiplier theorem.

**Lemma B.6** ([17]). *Let $1 < p < +\infty$. Assume that $m \in L^\infty(\mathbb{R}^3) \cap C^2(\mathbb{R}^3 \setminus \{0\})$ satisfies the Mikhlin condition*

$$[m]_{\mathcal{M}} := \sum_{|\alpha|\le 2} \sup_{\xi\in\mathbb{R}^3\setminus\{0\}} |\xi|^{|\alpha|}\,|\partial_\xi^\alpha m(\xi)| < +\infty.$$

*Then the Fourier multiplier $m(D)$, defined as in* (1.12)*, is bounded on $L^p(\mathbb{R}^3)$. More precisely,*

$$\|m(D)f\|_{L^p} \le C_p[m]_{\mathcal{M}}\|f\|_{L^p},$$

*where $C_p > 0$ depends only on $p$.*

**Lemma B.7.** *For any $1 < p < +\infty$, it holds that for $n_1$, $n_2$, and $n_3$ defined in* (B.1)*,* (B.2)*, and* (B.3)*, respectively,*

$$\|n_1(D)f\|_{L^p} \lesssim \|f\|_{W^{1,p}}, \quad \|n_2(D)f\|_{L^p} \lesssim \|f\|_{L^p}, \quad \|n_3(D)f\|_{L^p} \lesssim \|f\|_{L^p}.$$

*Proof.* We can apply the Hörmander-Mikhlin Theorem (Lemma B.6). Indeed, one can easily verify that the corresponding symbols $\langle\xi\rangle^{-1}n_1(\xi)$, $n_2(\xi)$, and $n_3(\xi)$ satisfy the Hörmander-Mikhlin multiplier conditions of order zero, uniformly in $\varepsilon \in (0, 1]$. □

As a consequence of the multiplier bounds and the definition of $Q(D)$ and $Q^{-1}(D)$ in (2.71), we obtain the following property for $Q(D)$ and $Q^{-1}(D)$.

**Corollary B.8.** *For any $1 < p < +\infty$, the operators $\chi_{\varepsilon,\kappa_0}(D)Q(D)$ and $\chi_{\varepsilon,\kappa_0}(D)Q^{-1}(D)$ are continuous on $L^p(\mathbb{R}^3)$ uniformly with respect to $\varepsilon \in (0, 1]$. More precisely, for any $F \in L^p(\mathbb{R}^3)$, it holds that*

$$\|\chi_{\varepsilon,\kappa_0}(D)Q(D)F\|_{L^p} \lesssim \|F\|_{L^p}, \qquad \|\chi_{\varepsilon,\kappa_0}(D)Q^{-1}(D)F\|_{L^p} \lesssim \|F\|_{L^p},$$

*where $\chi_{\varepsilon,\kappa_0}(\xi) = \chi\big(\sqrt{\tfrac{\varepsilon}{\kappa_0}}\xi\big)$, $\chi$ is defined in* (1.14)*.*

We also use the heat semigroup estimate from [29].

**Lemma B.9** ([29]). *Let $k \in \mathbb{N}^*$ and $1 < q < +\infty$. Then we have*

$$\|e^{\varepsilon t\Delta}(\varepsilon\Delta)^k\chi_{\varepsilon,\kappa_0}(D)f\|_{L^q} \lesssim (1+t)^{-k}\|f\|_{L^q},$$

*where $\chi_{\varepsilon,\kappa_0}(\xi) = \chi\big(\sqrt{\tfrac{\varepsilon}{\kappa_0}}\xi\big)$, $\chi$ is defined in* (1.14)*.*

The next lemma gives the $W^{k,p}$ estimates for the nonlinear terms in the system (2.3).

**Lemma B.10.** *Assume that $\|V\|_{H^{\frac{3}{2}+\delta}} \le \delta_1$ for $\delta_1 > 0$ sufficiently small. Then for every $1 < p < +\infty$, $\frac{1}{p} = \frac{1}{q_i} + \frac{1}{r_i}$, $1 < q_i, r_i < +\infty$, we have that for any $k \ge 0$, $\mathcal{B}$ defined in* (2.6) *and $\mathcal{Q}$ defined in* (2.7) *satisfy*

$$\begin{aligned}
\|\mathcal{B}(V,V)\|_{W^{k,p}} &\lesssim \|V\|_{W^{k+1,q_1}} \|V\|_{L^{r_1}} + \|V\|_{L^{r_2}} \|V\|_{W^{k+1,q_2}} ,\\
\|\mathcal{B}(V,V)\|_{W^{k,p}} &\lesssim \|n_3(D)V\|_{L^\infty} \|V\|_{W^{k+1,p}} ,\\
\|\mathcal{Q}(V)\|_{W^{k,p}} &\lesssim \|n_3(D)V\|_{L^\infty} \left(\|V\|_{W^{k+1,q_1}} \|V\|_{L^{r_1}} + \|V\|_{L^{r_2}} \|V\|_{W^{k+1,q_2}}\right) ,\\
\|\mathcal{Q}(V)\|_{W^{k,p}} &\lesssim \|n_3(D)V\|_{L^\infty}^2 \|V\|_{W^{k+1,p}} ,
\end{aligned}$$

*where $n_3$ is defined in* (B.3).

*Proof.* For simplicity, we only estimate the terms $\mathcal{B}_1(V,V)$ and $\mathcal{Q}(V)$, since $\mathcal{B}_2(V,V)$ can be treated in a similar manner. By the definition of $\mathcal{B}(V,V)$, the boundedness of the Riesz transform in $L^r$ $(1 < r < +\infty)$, Lemma B.4, Lemma B.7 and the fact

$$\|a(D)f\|_{L^r} \lesssim \|f\|_{W^{1,r}},$$

we have

$$\begin{aligned}
\|\mathcal{B}_1(V,V)\|_{W^{k,p}} &\lesssim \left\|\frac{|\nabla|}{a(D)} h \mathcal{R} q\right\|_{W^{k+1,p}}\\
&\lesssim \|n_3(D)h\|_{W^{k+1,q_1}} \|n_3(D)q\|_{L^{r_1}} + \|n_3(D)h\|_{L^{r_2}} \|n_3(D)q\|_{W^{k+1,q_2}}\\
&\lesssim \|V\|_{W^{k+1,q_1}} \|V\|_{L^{r_1}} + \|V\|_{L^{r_2}} \|V\|_{W^{k+1,q_2}},
\end{aligned}$$

and

$$\begin{aligned}
\|\mathcal{B}_1(V,V)\|_{W^{k,p}} &\lesssim \left\|\frac{|\nabla|}{a(D)} h \mathcal{R} q\right\|_{W^{k+1,p}}\\
&\lesssim \|n_3(D)h\|_{W^{k+1,p}} \|n_3(D)q\|_{L^\infty} + \|n_3(D)h\|_{L^\infty} \|n_3(D)q\|_{W^{k+1,p}}\\
&\lesssim \|V\|_{W^{k+1,p}} \|n_3(D)V\|_{L^\infty}.
\end{aligned}$$

We denote $g := \frac{|\nabla|}{a(D)} h$ and

$$\mathcal{H}(g) := \frac{1}{\gamma-1}(1 + \frac{|\nabla|}{a(D)} h)^{\gamma-1} - \frac{1}{\gamma-1} - \frac{|\nabla|}{a(D)} h - \frac{\gamma-2}{2} |\frac{|\nabla|}{a(D)} h|^2.$$

Then

$$\mathcal{H}(g) = \frac{(\gamma-2)(\gamma-3)}{2} \int_0^1 (1-\theta)^2 (1+\theta g)^{\gamma-4} g^3 \, \mathrm{d}\theta := \Psi(g) g^3,$$

together with the boundedness of the Riesz transform in $L^r$ $(1 < r < +\infty)$, Lemma B.4 and Lemma B.5, we have for any $j \le k$

$$\begin{aligned}
\|\mathcal{R}^* \nabla [\mathcal{H}(g)]\|_{\dot{W}^{j,p}} &\lesssim \left\|g^3\right\|_{\dot{W}^{j+1,p}} \|\Psi(g)\|_{L^\infty} + \left\|g^3\right\|_{L^{r_2}} \|\Psi(g)\|_{\dot{W}^{j+1,q_2}}\\
&\lesssim \left\|g^2\right\|_{\dot{W}^{j+1,p}} \|g\|_{L^\infty} + \left\|g^2\right\|_{L^{r_2}} \|g\|_{\dot{W}^{j+1,q_2}}\\
&\lesssim \|g\|_{L^\infty} \left(\|g\|_{\dot{W}^{j+1,q_1}} \|g\|_{L^{r_1}} + \|g\|_{L^{r_2}} \|g\|_{\dot{W}^{j+1,q_2}}\right)\\
&\lesssim \|n_3(D)V\|_{L^\infty} \left(\|V\|_{\dot{W}^{j+1,q_1}} \|V\|_{L^{r_1}} + \|V\|_{L^{r_2}} \|V\|_{\dot{W}^{j+1,q_2}}\right),
\end{aligned}$$

and

$$\begin{aligned}
\|\mathcal{R}^*\nabla\left[\mathcal{H}(g)\right]\|_{\dot{W}^{j,p}} &\lesssim \left\|g^3\right\|_{\dot{W}^{j+1,p}}\|\Psi(g)\|_{L^\infty}+\left\|g^3\right\|_{L^\infty}\|\Psi(g)\|_{\dot{W}^{j+1,p}}\\
&\lesssim \left\|g^2\right\|_{\dot{W}^{j+1,p}}\|g\|_{L^\infty}+\left\|g^2\right\|_{L^\infty}\|g\|_{\dot{W}^{j+1,p}}\\
&\lesssim \|g\|_{L^\infty}^2\|g\|_{\dot{W}^{j+1,p}}\\
&\lesssim \|n_3(D)V\|_{L^\infty}^2\|V\|_{\dot{W}^{j+1,p}},
\end{aligned}$$

where we have used the fact

$$\|g\|_{L^\infty}\lesssim\|g\|_{H^{\frac{3}{2}+\delta}}\lesssim\|h\|_{H^{\frac{3}{2}+\delta}}\le C\delta_1.$$

Summing the above estimates over $j\le k$ immediately yields the $W^{k,p}$ estimates for $\mathcal{Q}(V)$. This completes the proof. □

The same argument gives the corresponding estimates after applying $|\nabla|^{-(1-s)}$.

**Lemma B.11.** *Assume that $\|V\|_{H^{\frac{3}{2}+\delta}}\le\delta_1$ for $\delta_1>0$ sufficiently small. Then for every $1<p<+\infty$, $\frac{1}{p}=\frac{1}{q_i}+\frac{1}{r_i}$, $1<q_i,r_i<+\infty$, we have that for any $k\ge0$, $\mathcal{B}$ defined in* (2.6) *and $\mathcal{Q}$ defined in* (2.7) *satisfy*

$$\begin{aligned}
\||\nabla|^{-(1-s)}\mathcal{B}(V,V)\|_{W^{k,p}} &\lesssim \|V\|_{W^{k+s,q_1}}\|V\|_{L^{r_1}}+\|V\|_{L^{r_2}}\|V\|_{W^{k+s,q_2}},\\
\||\nabla|^{-(1-s)}\mathcal{B}(V,V)\|_{W^{k,p}} &\lesssim \|n_3(D)V\|_{L^\infty}\|V\|_{W^{k+s,p}},\\
\||\nabla|^{-(1-s)}\mathcal{Q}(V)\|_{W^{k,p}} &\lesssim \|n_3(D)V\|_{L^\infty}\left(\|V\|_{W^{k+s,q_1}}\|V\|_{L^{r_1}}+\|V\|_{L^{r_2}}\|V\|_{W^{k+s,q_2}}\right),\\
\||\nabla|^{-(1-s)}\mathcal{Q}(V)\|_{W^{k,p}} &\lesssim \|n_3(D)V\|_{L^\infty}^2\|V\|_{W^{k+s,p}}.
\end{aligned}$$

*where $n_3$ is defined in* (B.3).

For the normal form estimates, we introduce the following bilinear multiplier

$$\mathfrak{M}_{jk}:=\frac{|\xi|^{1-s}|\eta|^{1-s}|\xi-\eta|^{1-s}}{\langle\xi-\eta\rangle^{2\lambda}\langle\eta\rangle^{2\lambda}\Phi_{jk}}\widetilde{\chi}_{\varepsilon,\kappa_0}(\xi-\eta)\widetilde{\chi}_{\varepsilon,\kappa_0}(\eta)\widetilde{\chi}_{\varepsilon,\kappa_0}(\xi)\frac{a(\xi)}{2ib(\xi)}. \tag{B.4}$$

The next proposition verifies the required symbol bounds for this multiplier.

**Proposition B.12.** *Let $0<s<1$ and $j,k=1,2$. For $0\le\ell\le3/2-\delta$, we have*

$$\|\mathfrak{M}_{jk}\|_{\mathcal{M}^{\ell,\infty}_{\xi,\eta}}+\|\mathfrak{M}_{jk}\|_{\mathcal{M}^{\ell,\infty}_{\eta,\xi}}\lesssim1,$$

*where $\mathfrak{M}_{jk}$ is defined in* (B.4) *with $\lambda>\frac{7}{4}+s\ell-\frac{1}{2}s$.*

*Proof.* We only present the proof for the most singular phase

$$\Phi=-\Phi_{11}=b(|\xi|)-b(|\xi-\eta|)-b(|\eta|)$$

and for the corresponding kernel $\mathfrak{M}=-\mathfrak{M}_{11}$. The proof follows the same strategy as in [12, Section 4.2 and Appendix A], with the exponent $\sigma/2$ replaced by $s$.

The remaining sign combinations can be handled in the same way. Indeed, the phases $\Phi_{12}$ and $\Phi_{21}$ have the same singular structure as $\Phi$ after a harmless relabeling of the variables $\xi$, $\eta$, and $\xi-\eta$, possibly combined with a change of the overall sign. On the other hand, the phase $\Phi_{22}$ is non-resonant because $|\Phi_{22}(\xi,\eta)|=b(|\xi|)+b(|\xi-\eta|)+b(|\eta|)$, so that it is bounded from below by the elementary positivity and monotonicity properties of $b$.

It remains to prove the desired estimate for $\mathfrak{M}=-\mathfrak{M}_{11}$. Since the estimates are obtained by the same frequency decompositions as in [12, Section 4.2 and Appendix A], we only record the corresponding frequency regimes and the restrictions imposed on $\lambda$:

- (Case A) $|\xi| \le |\eta|/2$.
  In this regime, the condition for the kernel to belong to the desired space is $\lambda > \frac{7}{8} - \frac{3}{4}s$.

- (Case B) $|\eta| \le |\xi|/2$.
  For this interaction, the requirement is $\lambda > \frac{7}{4} + s\ell - \frac{1}{2}s$.

- (Case C) $1/3 < |\xi|/|\eta| < 3$.
  In this balanced frequency region, we require $\lambda > \max\left\{1 + s\ell - \frac{1}{2}s + \frac{1}{2}\ell, \frac{7}{4} + s\ell - \frac{1}{2}s\right\}$.

Analogously, to demonstrate that $\mathfrak{M} \in M^{\ell,\infty}_{\xi,\eta}$ with $0 \le \ell \le 3/2 - \delta$, we employ a similar classification, yielding the following constraints:

- (Case A) $|\xi| \le |\eta|/2$.
  We require $\lambda > \max\left\{\frac{5}{4} - s, \frac{5}{4} - s + \frac{1}{2}s\ell - \frac{1}{4}\ell\right\}$.

- (Case B) $|\eta| \le |\xi|/2$.
  We require $\lambda > 1 + s\ell - \frac{1}{2}s + \frac{\ell}{2}$.

- (Case C) $1/3 < |\xi|/|\eta| < 3$.
  We require $\lambda > \max\left\{\frac{5}{4} - s, \frac{5}{4} - s + \frac{1}{2}s\ell - \frac{1}{4}\ell, 1 + s\ell - \frac{1}{2}s + \frac{1}{2}\ell\right\}$.

In summary, by synthesizing the constraints from all aforementioned cases and considering $0 < s < 1$ and $0 \le \ell \le 3/2 - \delta$, we conclude that the desired result holds provided that

$$\lambda > \frac{7}{4} + s\ell - \frac{1}{2}s.$$

□

The two lemmas below give the bilinear estimates developed by Gustafson, Nakanishi and Tsai [16], which we use for the normal form analysis.

**Lemma B.13** ([20]). *Let $0 \le \ell < 3/2$, and suppose $\|\mathfrak{M}\|_{\mathcal{M}^{\ell,\infty}_{\eta,\xi}}, \|\mathfrak{M}\|_{\mathcal{M}^{\ell,\infty}_{\xi,\eta}} < \infty$. Then, we have*

$$\|T_{\mathfrak{M}}(f,g)\|_{L^{\frac{l_1}{l_1-1}}} \lesssim \|\mathfrak{M}\|_{\mathcal{M}^{\ell,\infty}_{\eta,\xi}} \|f\|_{L^2} \|g\|_{L^{l_2}},$$

*and*

$$\|T_{\mathfrak{M}}(f,g)\|_{L^2} \lesssim \|\mathfrak{M}\|_{\mathcal{M}^{\ell,\infty}_{\xi,\eta}} \|f\|_{L^{l_1}} \|g\|_{L^{l_2}},$$

*for $f, g \in \mathcal{S}$ where*

$$2 \le l_1, l_2 \le \frac{6}{3-2\ell} \quad \text{and} \quad \frac{1}{l_1} + \frac{1}{l_2} = 1 - \frac{\ell}{3}.$$

*Here, $\mathcal{S}$ is the Schwartz class.*

**Lemma B.14** ([12]). *Under the conditions of Lemma B.13, we have*

$$\|\langle\nabla\rangle^k T_{\mathfrak{M}}(f,g)\|_{L^{\frac{l_1}{l_1-1}}} \lesssim \|f\|_{H^k} \|g\|_{L^{l_2}} + \|f\|_{L^{l_2}} \|g\|_{H^k},$$

*and*

$$\|\langle\nabla\rangle^k T_{\mathfrak{M}}(f,g)\|_{L^2} \lesssim \|f\|_{W^{k,l_1}} \|g\|_{L^{l_2}} + \|f\|_{L^{l_2}} \|g\|_{W^{k,l_1}}.$$

We now apply these bilinear estimates to the multiplier defined in (B.4) to obtain the corresponding bounds.

**Lemma B.15.** *Define the phase function*

$$-\Phi_{jk}(\xi,\eta)=b(\xi)+(-1)^j b(\xi-\eta)+(-1)^k b(\eta),\quad j,k=1,2$$

*and the multiplier*

$$m(\xi,\eta)=\widetilde{\chi}_{\varepsilon,\kappa_0}(\xi-\eta)\widetilde{\chi}_{\varepsilon,\kappa_0}(\eta)\widetilde{\chi}_{\varepsilon,\kappa_0}(\xi)\frac{a(\xi)}{2\mathrm{i}b(\xi)}n_1(\xi)n_2(\xi-\eta)n_2(\eta).$$

*Then, the following estimates hold uniformly for $\varepsilon\in(0,1]$:*

$$\begin{aligned}\|T_{\frac{m}{\Phi_{jk}}}(f,g)\|_{W^{\sigma,\frac{p}{p-1}}}&\lesssim \||\nabla|^{-(1-s)}f\|_{H^{\sigma+s+2\lambda}}\||\nabla|^{-(1-s)}g\|_{W^{2\lambda,r}}\\&\quad+\||\nabla|^{-(1-s)}f\|_{W^{2\lambda,r}}\||\nabla|^{-(1-s)}g\|_{H^{\sigma+s+2\lambda}}\\ \|T_{\frac{m}{\Phi_{jk}}}(f,g)\|_{H^{\sigma}}&\lesssim \||\nabla|^{-(1-s)}f\|_{W^{\sigma+s+2\lambda,p}}\||\nabla|^{-(1-s)}g\|_{W^{2\lambda,r}}\\&\quad+\||\nabla|^{-(1-s)}f\|_{W^{2\lambda,r}}\||\nabla|^{-(1-s)}g\|_{W^{\sigma+s+2\lambda,p}},\end{aligned}$$

*where we choose $\lambda>\frac{7}{4}+s\ell-\frac{1}{2}s$ and $0\le\ell\le\frac{3}{2}-\delta$ such that*

$$2\le p,r\le\frac{6}{3-2\ell}\quad\text{and}\quad\frac{1}{p}+\frac{1}{r}=1-\frac{\ell}{3},$$

*the bilinear operator $T_{\frac{m}{\Phi_{jk}}}$ is defined by* (1.13)*, $n_1$ and $n_2$ are given in* (B.1) *and* (B.2)*.*

*Proof.* We begin by observing that

$$\frac{m}{\Phi_{jk}}=\mathfrak{M}_{jk}\frac{n_1(\xi)}{|\xi|^{1-s}}\frac{n_2(\xi-\eta)\langle\xi-\eta\rangle^{2\lambda}}{|\xi-\eta|^{1-s}}\frac{n_2(\eta)\langle\eta\rangle^{2\lambda}}{|\eta|^{1-s}}.$$

Next, we define the functions $\mathfrak{f}$ and $\mathfrak{g}$ via their Fourier transforms:

$$\widehat{\mathfrak{f}}(\xi-\eta):=\frac{n_2(\xi-\eta)\langle\xi-\eta\rangle^{2\lambda}}{|\xi-\eta|^{1-s}}\widehat{f}(\xi-\eta),\quad \widehat{\mathfrak{g}}(\eta):=\frac{n_2(\eta)\langle\eta\rangle^{2\lambda}}{|\eta|^{1-s}}\widehat{g}(\eta),$$

along with the multiplier

$$\widetilde{q}(\xi):=\frac{n_1(\xi)}{|\xi|^{1-s}\langle\xi\rangle^{s}}\widetilde{\chi}_{\varepsilon,\kappa_0}(\xi).$$

By Lemma B.7 and the Hörmander-Mikhlin theorem (Lemma B.6), the operator $\widetilde{q}(D)$ is bounded on $L^q$ $(1<q<+\infty)$, which, combined with Lemma B.14, yields

$$\begin{aligned}&\|T_{\frac{m}{\Phi_{jk}}}(f,g)\|_{W^{\sigma,\frac{p}{p-1}}}\\ \lesssim\ &\|q(D)\langle\nabla\rangle^{\sigma+s}T_{\mathfrak{M}_{jk}}(\mathfrak{f},\mathfrak{g})\|_{L^{\frac{p}{p-1}}}\\ \lesssim\ &\|\langle\nabla\rangle^{\sigma+s}T_{\mathfrak{M}_{jk}}(\mathfrak{f},\mathfrak{g})\|_{L^{\frac{p}{p-1}}}\\ \lesssim\ &\|\mathfrak{f}\|_{H^{\sigma+s}}\|\mathfrak{g}\|_{L^r}+\|\mathfrak{f}\|_{L^r}\|\mathfrak{g}\|_{H^{\sigma+s}}\\ \lesssim\ &\||\nabla|^{-(1-s)}f\|_{H^{\sigma+s+2\lambda}}\||\nabla|^{-(1-s)}g\|_{W^{2\lambda,r}}+\||\nabla|^{-(1-s)}f\|_{W^{2\lambda,r}}\||\nabla|^{-(1-s)}g\|_{H^{\sigma+s+2\lambda}}.\end{aligned}$$

Similarly, for the $H^\sigma$ estimate, we obtain

$$\begin{aligned}&\|T_{\frac{m}{\Phi_{jk}}}(f,g)\|_{H^{\sigma}}\\ \lesssim\ &\|\widetilde{q}(D)\langle\nabla\rangle^{\sigma+s}T_{\mathfrak{M}_{jk}}(\mathfrak{f},\mathfrak{g})\|_{L^2}\\ \lesssim\ &\|\langle\nabla\rangle^{\sigma+s}T_{\mathfrak{M}_{jk}}(\mathfrak{f},\mathfrak{g})\|_{L^2}\\ \lesssim\ &\|\mathfrak{f}\|_{W^{\sigma+s,p}}\|\mathfrak{g}\|_{L^r}+\|\mathfrak{f}\|_{L^r}\|\mathfrak{g}\|_{W^{\sigma+s,p}}\\ \lesssim\ &\||\nabla|^{-(1-s)}f\|_{W^{\sigma+s+2\lambda,p}}\||\nabla|^{-(1-s)}g\|_{W^{2\lambda,r}}+\||\nabla|^{-(1-s)}f\|_{W^{2\lambda,r}}\||\nabla|^{-(1-s)}g\|_{W^{\sigma+s+2\lambda,p}}.\end{aligned}$$

This completes the proof. ☐

Finally, we prove a time integral estimate used to close the decay bounds.

**Lemma B.16.** *Let $\varepsilon \in (0,1]$ and $m > 0$ be given parameters. Assume $\alpha \geq 0$, $\beta \geq 0$, and $\alpha + \beta > 1$. Then we have for any $t \geq 0$,*

$$\int_0^t (1+\varepsilon^m \tau)^{-\alpha}(1+\tau)^{-\beta}\,\mathrm{d}\tau \lesssim \begin{cases} 1, & \text{if } \beta > 1, \\ \ln(1+\varepsilon^{-m}), & \text{if } \beta = 1, \\ \varepsilon^{-m(1-\beta)}, & \text{if } \beta < 1. \end{cases}$$

*Proof.* For the case $\beta > 1$, using $(1+\varepsilon^m\tau)^{-\alpha} \leq 1$ and the integrability of $(1+\tau)^{-\beta}$ over $(0,\infty)$, we directly obtain $\mathcal{I}(t) \lesssim 1$.

For the case $\beta \leq 1$, we decompose the integral at $\tau = \varepsilon^{-m}$

$$\mathcal{I}(t) \leq \int_0^{\varepsilon^{-m}} (1+\tau)^{-\beta}\,\mathrm{d}\tau + \int_{\varepsilon^{-m}}^{\infty} (\varepsilon^m\tau)^{-\alpha}\tau^{-\beta}\,\mathrm{d}\tau := \mathcal{I}_1 + \mathcal{I}_2.$$

Direct integration yields $\mathcal{I}_1 \lesssim \ln(1+\varepsilon^{-m})$ for $\beta = 1$, and $\mathcal{I}_1 \lesssim \varepsilon^{-m(1-\beta)}$ for $\beta < 1$. For the term $\mathcal{I}_2$, since $\alpha + \beta > 1$, we evaluate it as follows

$$\mathcal{I}_2 = \varepsilon^{-m\alpha} \int_{\varepsilon^{-m}}^{\infty} \tau^{-(\alpha+\beta)}\,\mathrm{d}\tau \lesssim \varepsilon^{-m\alpha}(\varepsilon^{-m})^{1-\alpha-\beta} = \varepsilon^{-m(1-\beta)}.$$

When $\beta = 1$, it follows that $\mathcal{I}(t) \leq \mathcal{I}_1 + \mathcal{I}_2 \lesssim \ln(1+\varepsilon^{-m}) + 1 \lesssim \ln(1+\varepsilon^{-m})$. When $\beta < 1$, we see that both $\mathcal{I}_1$ and $\mathcal{I}_2$ are bounded by $\varepsilon^{-m(1-\beta)}$. This finishes the proof. ☐

## Appendix C. Well-posedness of the compressible Euler-Riesz system

In [12], Choi, Jung, and Lee studied the following problems for the compressible Euler-Riesz system,

$$\begin{cases} \partial_t \widetilde{\rho} + \operatorname{div}(\widetilde{\rho}\widetilde{u}) = 0, \\ \partial_t \widetilde{u} + \widetilde{u}\cdot\nabla\widetilde{u} + \dfrac{1}{\gamma-1}\nabla(\widetilde{\rho}^{\gamma-1}) + \nabla\widetilde{\phi} = 0, \\ (-\Delta)^s \widetilde{\phi} = \widetilde{\rho} - \bar{\rho}, \\ (\widetilde{\rho},\widetilde{u})(0,x) = (\widetilde{\rho}_0,\widetilde{u}_0)(x), \quad \lim\limits_{|x|\to\infty} (\widetilde{\rho},\widetilde{u}) = (\bar{\rho},0). \end{cases} \tag{C.1}$$

We summarize the following lemma from their result [12, Theorem 1.1], which states the global existence and properties of smooth solutions to the system (C.1).

**Lemma C.1.** *For $0 < s < 1$, let $p$ be chosen according to*

$$\begin{cases} 6 < p < \infty & \text{when } 0 < s \leq 1/2, \\ 8 < p < \infty & \text{when } 1/2 < s < 1. \end{cases}$$

*There exists $\varepsilon_0 > 0$ such that for any initial perturbation $(\rho_0, u_0)$ satisfying*

$$\||\nabla|^{-(1-s)}(\widetilde{\rho}_0 - \bar{\rho}, \widetilde{u}_0)\|_{H^N} + \|(\widetilde{\rho}_0 - \bar{\rho}, \widetilde{u}_0)\|_{W^{N_0+3(1-\frac{2}{p}),\frac{p}{p-1}}} \leq \varepsilon_0 \quad \text{and} \quad \nabla\times\widetilde{u}_0 = 0,$$

*with $N_0 \geq 6$ and $N \geq N_0 + 14$, the system admits a global solution $(\widetilde{\rho}, \widetilde{u})$ with*

$$\sup_{0\leq t<+\infty} \{\||\nabla|^{-(1-s)}(\widetilde{\rho}-\bar{\rho},\widetilde{u})(t)\|_{H^N} + (1+t)^{\beta_s}\|(\widetilde{\rho}-\bar{\rho},\widetilde{u})(t)\|_{W^{N_0,p}}\} \leq 2\varepsilon_0,$$

*where $\beta_s$ is defined by*

$$\beta_s = \begin{cases} \frac{3}{2}\left(1 - \frac{2}{p}\right) & \text{when } 0 < s \le 1/2, \\ \frac{4}{3}\left(1 - \frac{2}{p}\right) & \text{when } 1/2 < s < 1. \end{cases}$$

**Remark C.2.** *It should be noted that although Lemma C.1 is essentially taken from* [12, Theorem 1.1]*, the result requires less regularity for initial data in $W^{N_0,p}$ and $H^N$ due to the refined estimates of the phase function $\Phi$ as shown in Lemma A.2.*

**Acknowledgments:** A.Ś.-G. and J.H.Z. are supported in part by National Science Centre, Poland, grant number 2023/50/A/ST1/00447. Y.X. and J.H.Z. would like to thank Prof. Fucai Li and Prof. Lvqiao Liu for helpful comments.

**Conflict of Interest:** The authors declare that they have no conflict of interest.

**Data availability:** Data sharing not applicable to this article as no datasets were generated or analyzed during the current study.

University of Warsaw, Institute of Applied Mathematics and Mechanics, Banacha 2, 02-097 Warszawa, Poland

*Email address*: `aswiercz@mimuw.edu.pl`

Central China Normal University, School of Mathematics and Statistics, Wuhan, 430079, China

*Email address*: `yuanxu2019@163.com`

University of Warsaw, Institute of Applied Mathematics and Mechanics, Banacha 2, 02-097 Warszawa, Poland

*Email address*: `jhzhang@mimuw.edu.pl`